\documentclass[a4paper]{easychair}
\usepackage{bussproofs}
\usepackage{amsmath,amssymb}
\usepackage{mathrsfs}
\usepackage{enumerate}
\usepackage{stmaryrd}
\usepackage{setspace}
\newtheorem{theorem}{Theorem}[section]
\newtheorem{lemma}[theorem]{Lemma}
\newtheorem{proposition}[theorem]{Proposition}
\newtheorem{corollary}[theorem]{Corollary}
\newtheorem{fact}[theorem]{Fact}
\newtheorem{remark}[theorem]{Remark}
\theoremstyle{plain}
\newtheorem{definition}[theorem]{Definition}
\newtheorem{example}[theorem]{Example}
\renewcommand{\phi}{\varphi}

\newcommand{\bd}{\blacklozenge}
\newcommand{\bb}{\blacksquare}
\newcommand{\sub}{\subseteq}

\newcommand{\comb}[1]{(#1)^{c}}
\newcommand{\ol}{\overline}

\newcommand{\ve}{\varnothing}
\newcommand{\D}{\Diamond}
\newcommand{\B}{\Box}
\newcommand{\ua}{\uparrow}

\newcommand{\cset}[1]{{\{ #1 \}}}
\newcommand{\tup}[1]{\langle #1 \rangle}

\newcommand{\rsto}{{\upharpoonright}}
\newcommand{\iso}{\cong}

\newcommand{\R}{R_\sharp}
\renewcommand{\S}{S_\sharp}

\newcommand{\gf}{\mathbb{F}}
\renewcommand{\gg}{\mathbb{G}}
\newcommand{\F}{\mathfrak{F}}
\newcommand{\G}{\mathfrak{G}}
\renewcommand{\H}{\mathfrak{H}}

\newcommand{\Z}{\mathfrak{Z}}
\newcommand{\M}{\mathfrak{M}}

\newcommand{\C}{\mathfrak{C}}
\renewcommand{\L}{ \mathscr{L}_t}

\newcommand{\J}{\mathcal{J}}

\newcommand{\md}{\models}
\newcommand{\NExt}{\mathsf{NExt}}
\newcommand{\Ext}{\mathsf{Ext}}
\usepackage{tikz}
\tikzset{shorten <>/.style={shorten >=#1, shorten <=#1}}
\usetikzlibrary{patterns,arrows}
\title{Pretabular Tense Logics over $\mathsf{S4}_t$}
\author{
    Qian Chen
}

\institute{The Tsinghua-UvA JRC for Logic, Department of Philosophy, Tsinghua University, China  \and  Institute for Logic, Language and Computation, University of Amsterdam,  The Netherlands}

\authorrunning{~}

\titlerunning{~}

\begin{document}

\maketitle

\begin{abstract}
    A logic $L$ is called tabular if it is the logic of some finite frame and $L$ is pretabular if it is not tabular while all of its proper consistent extensions are tabular. In this work, we study pretabular tense logics in the lattice $\NExt(\mathsf{S4}_t)$ of all extensions of $\mathsf{S4}_t$, tense $\mathsf{S4}$. For all $0<n,m,k,l\leq\aleph_0$, we define the tense logic $\mathsf{S4BP}_{n,m}^{k,l}$ with respectively bounded width, depth and z-degree. We give a full characterization of the set $\mathsf{PTAB}(\mathsf{S4.3}_t)$ of all pretabular logics extending $\mathsf{S4.3}_t$, which entails that there are exactly 5 pretabular logics in $\NExt(\mathsf{S4.3}_t)$. Moreover, by providing a full characterization of $\mathsf{PTAB}(\mathsf{S4BP}^{2,\omega}_{2,2})$ and proving that $|\mathsf{PTAB}(\mathsf{S4BP}^{2,\omega}_{2,3})|=2^{\aleph_0}$, we show the anti-dichotomy theorem for cardinality of pretabular extensions in $\NExt(\mathsf{S4}_t)$: for all cardinal $\kappa$ such that $\kappa\leq{\aleph_0}$ or $\kappa=2^{\aleph_0}$, $|\mathsf{PTAB}(L)|=\kappa$ for some $L\in\NExt(\mathsf{S4}_t)$. It follows that $|\mathsf{PTAB}(\mathsf{S4}_t)|=2^{\aleph_0}$, which answers an open problem concerning the cardinality of $\mathsf{PTAB}(\mathsf{S4}_t)$ raised by Rautenberg in 1979.
\end{abstract}

\section{Introduction}

A logic is called \textit{tabular} if it is the logic of some finite frame or algebra. Tabular modal logics have been well-investigated, see \cite{Blok1980,Chagrov.Zakharyaschev1997}. A logic is called \textit{pretabular} if it is not tabular but all of its proper consistent extensions are tabular. The concept of pretabularity was introduced by Kuznetsov \cite{Kuznetsov1971} in order to study tabularity. 
In general, given a lattice $\mathcal{L}$ of logics, tabularity in $\mathcal{L}$ is decidable if all of the following holds: (i) every non-tabular logic has a pretabular extension; (ii) there are finitely many pretabular logics in $\Ext(\mathsf{IPC})$; and (iii) every pretabular extension of $\mathsf{IPC}$ is decidable. Consider the intuitionistic proposition logic $\mathsf{IPC}$ the lattice $\Ext(\mathsf{IPC})$ of all intermediate logics. Kuznetsov \cite{Kuznetsov1971} showed that (i) holds for $\Ext(\mathsf{IPC})$ and every pretabular extension of $\mathsf{IPC}$ has the FMP. Maksimova \cite{Maksimova1972} provided a full characterization of pretabular intermediate logics, which entails $|\mathsf{PTAB}(\mathsf{IPC})|=3$ and every pretabular intermediate logic is finitely axiomatizable. Combining the results above, we see that tabularity in $\Ext(\mathsf{IPC})$ is decidable: given any formula $\phi$, there is a effective method to decide whether $\mathsf{IPC}\oplus\phi$ is tabular.

A number of by-now famous results on pretabular modal logics have been obtained since the 1970s. Let $\mathsf{PTAB}(L)$ denote the set of pretabular logics in the lattice $\NExt(L)$ of normal extensions of a normal modal logic $L$. Esakia and Meskhi \cite{Esakia.Meskhi1977} and Maksimova~\cite{Maksimova1975} provided full characterizations of $\mathsf{PTAB}(\mathsf{S4})$ and showed that $|\mathsf{PTAB}(\mathsf{S4})|=5$ independently, where $\mathsf{S4}$ is the modal logic of pre-orders. Bellissima \cite{Bellissima1988} studied the finitely alternative modal logics $\mathsf{KAlt}_n$ and proved that $|\mathsf{PTAB}(\mathsf{KAlt}_1)|=1$. Note that every pretabular logic mentioned above is decidable, therefore tabularity in $\NExt(\mathsf{S4})$ and $\NExt(\mathsf{KAlt}_1)$ are decidable.
On the other hand, Blok \cite{Blok1980b} showed that there exist $2^{\aleph_0}$ pretabular logics extending the transitive modal logic $\mathsf{K4}$, which means that the lattice $\NExt(\mathsf{K4})$ of all extensions of $\mathsf{K4}$ is highly complex and decidability of tabularity cannot be established in the same way. In fact, Blok \cite{Blok1980a} proved that tabularity in $\NExt(\mathsf{K})$ is undecidable and the decidability of tabularity in $\NExt(\mathsf{K4})$ remains unknown \cite[Section~5]{Rautenberg.Zakharyaschev.ea2006}. 

A \emph{tense logic} is a bi-modal logic which includes a future-looking necessity modality $\B$ and a past-looking possibility modality $\bd$. From an algebraic perspective, these two modalities are adjoint in the sense that, for any tense logic $L$ and formulas $\phi$ and $\psi$, we have $\bd\phi\to\psi\in L$ if and only if $\phi\to\B\psi\in L$. The addition of a past-looking modality makes the lattices of tense logics significantly more complex than those of unary modal logics. Makinson \cite{Makinson1971} showed that $\NExt(\mathsf{K})$ has only 2 co-atoms, while \cite{Chen.Ma2024a} showed that there are continuum many co-atoms in $\NExt(\mathsf{K}_t)$. These differences are apparent even in very small sublattices. For example, consider the modal logic $\mathsf{S4.3}$ and the tense logic $\mathsf{S4.3}_t$ of linear frames. It was proved by Bull \cite{Bull1966} and Fine \cite{Fine1971} that every modal logic in the lattice $\NExt(\mathsf{S4.3})$ enjoys the finite model property (FMP), but there are countably many tense logics in $\NExt(\mathsf{S4.3}_t)$ which do not have the FMP (see \cite{Wolter1995}). For more on the differences and similarities between the lattices of modal logics and tense logics, we refer the readers to \cite{Chen.Ma2024a,Rautenberg1979,Kracht1992,Thomason1972,Ma.Chen2021,Ma.Chen2023}.

The aim of this paper is to develop a better understanding of the lattices of tense logics, with a particular focus on tabularity and pretabularity. Tabular tense logics were studied in \cite{Chagrov.Shehtman1995} and \cite{Chen.Ma2024}, where different characterizations of tabular tense logics are provided. For pretabularity, Ma and Chen \cite{Ma.Chen2021} investigated the tense analogs $T_{n,m}$ of alternative modal logics $\mathsf{KAlt}_n$ and showed that $|\mathsf{PTAB}(T_{1,1})|=1$ and $|\mathsf{PTAB}(T_{1,2})|=|\mathsf{PTAB}(T_{2,1})|\geq\aleph_0$. However, unlike the modal case, the pretabularity in the lattice $\NExt(\mathsf{S4}_t)$ has not been fully understood yet, where $\mathsf{S4}_t$ is the tense logic of pre-orders. Kracht \cite{Kracht1992} defined a pretabular tense logic $\mathsf{Ga}\in\NExt(\mathsf{S4}_t)$ of depth 2, whose modal fragment is tabular. Rautenberg \cite[Page 23]{Rautenberg1979} posed the question: how many pretabular tense logics exist in the lattice $\NExt(\mathsf{S4}_t)$? As a partial answer, Rautenberg claimed in the same paper that there are infinitely many pretabular tense logics extending $\mathsf{S4}_t$. However, no proof was provided, and the only explanation for this claim was that the tense logic $\mathsf{Ga}$ is pretabular and has dimension 3. Since Kracht \cite{Kracht1992} later demonstrated that $\mathsf{Ga}$ does not have dimension 3, we believe determining the exact cardinality of $\mathsf{PTAB}(\mathsf{S4}_t)$ remains an open problem.

In this paper, we investigate pretabular tense logics in $\NExt(\mathsf{S4}_t)$, and our main goal is to determine the cardinality of $\mathsf{PTAB}(\mathsf{S4}_t)$. As indicated by Kracht \cite{Kracht1992}, the lattice $\NExt(\mathsf{S4}_t)$ is far more complex than its modal counterpart $\NExt(\mathsf{S4})$. We therefore begin by studying sublattices of $\NExt(\mathsf{S4}_t)$. For modal logics, transitive modal logics with bounded forth-width and depth have been thoroughly studied, see \cite{Fine1974,Fine1985,Chagrov.Zakharyaschev1997}. In the tense setting, however, the presence of a past-looking modality motivates to study logics not only with bounded depth and forth-width, but also with back-width and \emph{z-degree}. Intuitively, a tense logic $L$ is of \emph{z-degree} $l$ if every rooted frames for $L$ can be generated within $l$ back-and-forth steps from any point. As we shall see later, the z-degree is a crucial parameter for tense logics. 

For all $n,m,k,l\in\mathbb{Z}^+\cup\cset{\omega}$, we define $\mathsf{S4BP}_{n,m}^{k,l}$ to be the tense logic of all frames with forth-width, back-width, depth and z-degree no more than $n,m,k$ and $l$, respectively. We give first full characterizations of $\mathsf{PTAB}(\mathsf{S4BP}_{n,m}^{k,l})$ where $n,m,k,l$ are finite. Then we move to lattices $\NExt(\mathsf{S4BP}_{n,m}^{k,l})$ where some of the parameters are infinite. The first logic we consider is $\mathsf{S4BP}_{1,1}^{\omega,1}$, where the depth is unbounded and all other parameters are bounded by $1$. Note that $\mathsf{S4BP}_{1,1}^{\omega,1}$ is exactly the tense logic $\mathsf{S4.3}_t$ of all linear frames. We give a full characterization of $\mathsf{PTAB}(\mathsf{S4.3}_t)$ and prove that $|\mathsf{PTAB}(\mathsf{S4.3}_t)|=5$. Moreover, we show that every logic in $\mathsf{PTAB}(\mathsf{S4.3}_t)$ is finitely axiomatizable and has the FMP. As a corollary, tabularity in $\mathsf{NExt}(\mathsf{S4.3}_t)$ is decidable. Then we turn our attention to the lattice $\NExt(\mathsf{S4BP}^{2,\omega}_{2,2})$, where the logic $\mathsf{Ga}$ is proved to be the unique pretabular logic of infinite z-degree. We provide a characterization of the rooted frames of $\mathsf{Ga}$ and identify an error in the characterization of $\NExt(\mathsf{Ga})$ given in \cite{Kracht1992}. Moreover, we provide a full characterization of $\mathsf{PTAB}(\mathsf{S4BP}_{2,2}^{2,\omega})$ and show that every logic in $\mathsf{PTAB}(\mathsf{S4BP}_{2,2}^{2,\omega})$ has the FMP. It follows from the characterization that $|\mathsf{PTAB}(\mathsf{S4BP}_{2,2}^{2,\omega})|=\aleph_0$.

Based on the above, we can already observe a key difference between $\mathsf{PTAB}(\mathsf{S4})$ and $\mathsf{PTAB}(\mathsf{S4}_t)$: the former is finite, whereas the latter is infinite. Moreover, there exist pretabular tense logic, such as $\mathsf{Ga}$, whose modal fragment are not pretabular. Next, in order to address the problem of determining the cardinality of $\mathsf{PTAB}(\mathsf{S4}_t)$, we move from the lattice $\NExt(\mathsf{S4BP}^{2,\omega}_{2,2})$ to $\NExt(\mathsf{S4BP}^{2,\omega}_{2,3})$, where the bound of back-width is increased from $2$ to $3$. The main result we obtain is that $|\mathsf{PTAB}(\mathsf{S4BP}_{2,3}^{2,\omega})|=2^{\aleph_0}$. To prove this result, we introduce the generalized Thue-Morse sequences, generalized Jankov formulas and local t-morphisms. It follows that the following anti-dichotomy theorem for cardinality of pretabular extensions in $\NExt(\mathsf{S4}_t)$ holds: 
\begin{center}
    For all $\kappa\leq{\aleph_0}$ or $\kappa=2^{\aleph_0}$, there exists $L\in\NExt(\mathsf{S4}_t)$ such that $|\mathsf{PTAB}(L)|=\kappa$.
\end{center}
Consequently, we determine that the cardinality of $\mathsf{PTAB}(\mathsf{S4}_t)$ is $2^{\aleph_0}$. This gives a full solution to the problem presented in \cite{Rautenberg1979}. At the same time, this result also indicates that decidability of tabularity in $\NExt(\mathsf{S4}_t)$ cannot be established via pretabular logics.

The paper is structured as follows. Section 2 recalls the preliminaries of tense logic. Section 3 introduces generalized Jankov formulas and local t-morphisms. Section 4 presents tense logics with bounded parameters $\mathsf{S4BP}_{n,m}^{k,l}$ and studies their basic properties. Section 5 gives a characterization of $\mathsf{PTAB}(\mathsf{S4BP}_{n,m}^{k,l})$ for all $n,m,k,l\in\mathbb{Z}^+$. Section 6 provides a characterization of $\mathsf{PTAB}(\mathsf{S4.3}_t)$. Section 7, gives a characterization of $\mathsf{PTAB}(\mathsf{S4BP}_{2,2}^{2,\omega})$ and proves the anti-dichotomy theorem for the cardinality of pretabular logics extending $\mathsf{S4BP}_{2,2}^{2,\omega}$. It also highlights the connection between Kracht's work on the lattice $\NExt(\mathsf{Ga})$ in \cite[Section 4]{Kracht1992} and our results. Section 8 introduces generalized Thue-Morse sequences and defines a continual family of pretabular logics in $\NExt(\mathsf{S4BP}_{2,3}^{2,\omega})$, which shows that $|\mathsf{PTAB}(\mathsf{S4}_t)|=2^{\aleph_0}$. Section 9 concludes the paper.

\section{Preliminaries}

Basic notions on modal and tense logic can be found in e.g. \cite{Ma.Chen2021,Chagrov.Zakharyaschev1997,Blackburn.deRijke.ea2001}. Let $\mathbb{N}$, $\mathbb{Z}$ and $\mathbb{Z}^+$ be sets of all natural numbers, integers and positive integers respectively. We use $\mathbb{O}$ and $\mathbb{E}$ for sets of odd numbers and even numbers in $\mathbb{N}$ respectively. Let $\mathbb{Z}^*=\mathbb{Z}^+\cup\cset{\omega}$. The cardinal of a set $X$ is denoted by $|X|$. The power set of $X$ is denoted by $\mathcal{P}(X)$. The set of all finite tuples of elements of $X$ is denoted by $X^*$. We use Boolean operations $\cap$, $\cup$ and $\comb{\cdot}$ (complementation) on $\mathcal{P}(X)$. Let $X,X'$ be sets. We write $f:X\to X'$ if $f$ is a function from $X$ to $X'$. Moreover, $f\sub X\times X'$ is called a partial function from $X$ to $X'$ if $f:Y\to X'$ for some $Y\sub X$. Given any (partial) function $f:X\to X'$, we write $\mathsf{dom}(f)$ and $\mathsf{ran}(f)$ for its domain and range, respectively. In this paper, we sometimes write an ordinal $\lambda$ for the set $\cset{\alpha:\alpha<\lambda}$.

\begin{definition}
The {language of tense logic} consists of a denumerable set of variables $\mathsf{Prop}=\{p_i: i\in\mathbb{N}\}$, connectives $\bot$ and $\to$, and unary tense operators $\B$ and $\bd$. The set of all formulas $\L$ is defined as follows:
\[
\L\ni \phi ::= p \mid \bot \mid (\phi\to\phi) \mid \B\phi\mid \bd\phi,~\text{where $p \in \mathsf{Prop}$.}
\]
Formulas in $\mathsf{Prop}\cup\{\bot\}$ are called {\em atomic}. 
The connectives $\top,\neg,\wedge$ and $\vee$ are defined as usual. Let $\D\phi:=\neg\B\neg\phi$ and $\bb\phi:=\neg\bd\neg\phi$. Let $var(\phi)$ be the set of all variables in a formula $\phi$. For each $n\in\omega$, let $\mathsf{Prop}(n)=\cset{p_i:i<n}$ and $\L(n)=\cset{\phi\in\L:var(\phi)\sub\mathsf{Prop}(n)}$.  The {\em complexity} $c(\phi)$ of $\phi$ is defined by 
\begin{align*}
c(p) &= 0 = c(\bot),\\
c(\phi\to\psi)&=\max\{c(\phi),c(\psi)\}+1,\\
c(\B\phi)&=c(\phi)+1=c(\bd\phi).
\end{align*}
The {\em modal degree} $md(\phi)$ of $\phi$ is defined inductively as follows: 
\begin{align*}
md(p) &= 0 = md(\bot),\\
md(\phi\to\psi)&=\max\{md(\phi),md(\psi)\},\\
md(\B\phi)&=md(\phi)+1=md(\bd\phi).
\end{align*}
The set of all subformulas of $\phi$ is denoted by $Sub(\phi)$. A set of formulas $\Sigma$ is {\em closed under subformulas} if $\Sigma=\bigcup_{\phi\in\Sigma}Sub(\phi)$. A {\em substitution} is a homomorphism $(\cdot)^s:\L\to \L$ on the formula algebra $\L$.
\end{definition}

\begin{definition}
    A {\em frame} is a pair $\F =(X,R)$ where $X\neq\ve$ and $R\sub X\times X$. We write $Rxy$ if $\tup{x,y}\in R$. The {\em inverse of $R$} is defined as $\breve{R}=\{\tup{y,x}:Rxy\}$. The {\em inverse of $\F$} is defined as the frame $\breve{\F}=(X,\breve{R})$. For every $Y\sub X$, we define $R[Y]=\cset{z\in X:\exists{y\in Y}(Ryz)}$ and $\breve{R}[Y]=\cset{z\in X:\exists{y\in Y}(Rzy)}$. Let $x\in X$. We write $R[x]$ and $\breve{R}[x]$ for $R[\cset{x}]$ and $\breve{R}[\cset{x}]$, respectively.
    For all $n\geq 0$, we define $\R^n[x]$ by:
    \begin{center}
        $\R^0[x] = \{x\}$ and $\R^{k+1}[x] =\R^k[x]\cup R[\R^k[x]]\cup \breve{R}[\R^k[x]]$.
    \end{center}
    Let $\R[x] = \R^1[x]$ and $\R^\omega[x] = \bigcup_{k\geq 0}\R^k[x]$.
    Let $\mathsf{Fr}$ and $\mathsf{Fin}$ denote the class of all frames and finite frames, respectively.
\end{definition}

\begin{definition}
    A {\em general frame} is a triple $\gf=(X,R,A)$ where $(X,R)$ is a frame and $A\sub\mathcal{P}(X)$ is a set such that $\ve\in A$ and $A$ is closed under $\cap$, $(\cdot)^c$, $R[\cdot]$ and $\breve{R}[\cdot]$. We call $A\sub\mathcal{P}(X)$ {\em the set of internal sets in $\gf$}.

    Let $\mathsf{GF}$ denote the class of all general frames.
\end{definition}
It is not hard to verify that for every family $(\gf_i=(X,R,A_i))_{i\in I}$ of general frames, $\gf=(X,R,\bigcap_{i\in I}A_i)$ is again a general frame. Thus, for any frame $\F=(X,R)$ and $Q\sub\mathcal{P}(X)$, there exists the smallest subset $[Q]$ of $\mathcal{P}(X)$ such that $Q\sub[Q]$ and $(\F,[Q])$ is a general frame. A general frame $\gf=(X,R,A)$ is called {\em finitely generated}, if $A=[Q]$ for some finite $Q\sub\mathcal{P}(X)$.

Let $X$ be a set and $A\sub\mathcal{P}(X)$. We say that $A$ has the {\em finite intersection property (FIP)}, if $\bigcap B\neq\ve$ for any non-empty finite subset $B$ of $A$.

\begin{definition}
    Let $\gf=(X,R,A)$ be a general frame. Then $\gf$ is called
    \begin{enumerate}[(i)]
        \item {\em differentiated}, if for all distinct $x,y\in X$, $x\in U$ and $y\not\in U$ for some $U\in A$;
        \item {\em tight}, if for all $x,y\in X$ such that $y\not\in R[x]$, there are internal sets $U,V\in A$ such that $x\in U\setminus\breve{R}[V]$ and $y\in V\setminus R[U]$;
        \item {\em compact}, if $\bigcap B\neq\ve$ for any $B\sub A$ which has the FIP;
        \item {\em refined}, if $\gf$ is both differentiated and tight;
        \item {\em descriptive}, if $\gf$ is both refined and compact.
    \end{enumerate}
    Let $\mathsf{RF}$ and $\mathsf{DF}$ denote the class of all refined and descriptive frames, respectively.
\end{definition}

For each general frame $\gf=(X,R,A)$, we call $\kappa\gf=(X,R)$ the {\em underlying frame of $\gf$}. Moreover, if a general frame is of the form $\gf=(X,R,\mathcal{P}(X))$, then we identify $\gf$ with its underlying frame $\kappa\gf$. In this sense, we see that $\mathsf{Fr}\sub\mathsf{RF}\sub\mathsf{GF}$.

\begin{definition}
    Let $\gf=(X,R,A)$ be a general frame. Then a map $V:\mathsf{Prop}\to A$ is called a {\em valuation in $\gf$}. A valuation $V$ is extended to $V:\L\to A$ as follows:
    \begin{align*}
        V(\bot) &=\ve,
        &V(\phi\to\psi) &= (V(\phi))^c\cup V(\psi),\\
        V(\bd\phi) &= R[V(\phi)],
        &V(\B\phi) &= (\breve{R}[V(\phi)^c])^c.
    \end{align*}
    A {\em model} is a pair $\M=(\gf,V)$ where $\gf\in\mathsf{GF}$ and $V$ a valuation in $\gf$. Let $\phi$ be a formula and $w\in X$. Then (i) $\phi$ is {\em true} at $w$ in $\M$ (notation: $\M,w\models\phi$) if $w\in V(\phi)$; (ii) $\phi$ is {\em valid} at $w$ in $\gf$ (notation: $\gf,w\models\phi$) if $w\in V(\phi)$ for every valuation $V$ in $\gf$; (iii) $\phi$ is {\em valid} in $\gf$ (notation: $\gf\models\phi$) if $\gf,w\models\phi$ for every $w\in X$; and (iv) $\phi$ is {\em valid} in a class of general frame $\mathcal{K}$ (notation: $\mathcal{K}\models\phi$) if $\gf\models\phi$ for every $\gf\in\mathcal{K}$. For all sets $\Sigma\sub\L$ of formulas and classes $\mathcal{K}\sub\mathsf{GF}$ of general frames, let 
    \begin{center}
        $\mathcal{K}(\Sigma)=\cset{\gf\in\mathcal{K}:\gf\md\Sigma}$ and $\mathsf{Log}(\mathcal{K})=\cset{\phi:\mathcal{K}\md\phi}$.
    \end{center}
    We call $\mathsf{Log}(\mathcal{K})$ the {\em tense logic of $\mathcal{K}$}.
\end{definition}
We often write $(X,R,V)$ for the model $(\gf,V)$. Let $\M=(X,R,V)$ be a model and $\Sigma$ a set of formulas. Then for all $x,y\in X$, we say $x$ and $y$ are {\em equivalent with respect to $\Sigma$} (notation: $x\equiv_\Sigma y$) if $\cset{\phi\in\L:\M,x\md\phi}=\cset{\phi\in\L:\M,y\md\phi}$.

\begin{fact}
    Let $\M=(X,R,V)$ be a model and $x\in X$. Then for all $\phi\in\L$,
    \begin{enumerate}[(1)]
        \item $\M,x\md\bd\phi$ if and only if $\M,y\md\phi$ for some $y\in\breve{R}[x]$, and
        \item $\M,x\md\B\phi$ if and only if $\M,y\md\phi$ for all $y\in R[x]$.
    \end{enumerate}
\end{fact}

\begin{definition}
    Let $\gf=(X,R,A)$ be a general frame. For every $Y\sub X$, the {\em subframe of $\gf$ induced by $Y$} is defined as $\gf\rsto Y=(Y,R\rsto Y,A\rsto Y)$, where $R\rsto Y=R\cap(Y\times Y)$ and $A\rsto Y=\cset{U\cap Y:U\in A}$. For all $x\in X$, let $\gf_x$ denote the frame $\gf\rsto \R^\omega[x]$ and we call $\gf_x$ the {\em subframe of $\gf$ generated by $x$}. We call $\gf$ a {\em rooted frame} if $\gf=\gf\rsto \R^\omega[x]$ for some $x\in X$. For each class $\mathcal{K}$ of general frames, we write $\mathcal{K}_r$ for the subclass $\cset{\gf_x:\gf\in\mathcal{K}\text{ and }x\in\gf}$ of all rooted elements in $\mathcal{K}$.
    
    Let $\mathcal{F}=(\gf_i=(X_i,R_i,A_i))_{i\in I}$ be a family of general frames. Then the {\em disjoint union of $\mathcal{F}$} is defined as $\bigoplus_{i\in I}\gf_i=(X,R,A)$, where $X=\bigcup_{i\in I}X_i\times\cset{i}$, $R=\cset{\tup{\tup{x,i},\tup{y,i}}:R_ixy\text{ and }i\in I}$ and $A=\cset{U\sub X:\forall_{i\in I}(\cset{y\in X_i:\tup{y,i}\in U}\in A_i)}$. 

    Let $\gf=(X,R,A)$ and $\gf'=(X',R',A')$ be general frames. A map $f:X\to X'$ is said to be a {\em t-morphism from $\gf$ to $\gf'$} (notation: $f:\gf\to\gf'$), if 
    \begin{itemize}
        \item $f^{-1}[Y']\in A$ for any $Y'\in A'$; and
        \item for all $x\in X$, $f[R[x]]= R'[f(x)]$ and $f[\breve{R}[x]]=\breve{R'}[f(x)]$.
    \end{itemize}
    We write $f:\gf\rightarrowtail\gf'$, $f:\gf\twoheadrightarrow\gf'$ and $f:\gf\cong\gf'$ if $f$ is an injective, surjective and bijective t-morphism from $\gf$ to $\gf'$, respectively. Moreover, $\gf'$ is called (i) a {\em t-morphic image of $\gf$ (notation: $\gf\twoheadrightarrow\gf'$)} if there exists $f:\gf\twoheadrightarrow\gf'$; and (ii) an {\em isomorphic image of $\gf$ (notation: $\gf\cong\gf'$)} if there exists $f:\gf\cong\gf'$. For each class $\mathcal{K}$ of general frames, let $\mathsf{TM}(\mathcal{K})$ and $\mathsf{IM}(\mathcal{K})$ denote the class of all t-morphic images and isomorphic image of frames in $\mathcal{K}$, respectively.
\end{definition}

\begin{fact}\label{fact:inj-tmorphism}
    Let $\F=(X,R)$, $\F'=(X',R')$ be frames and $f:\F\rightarrowtail\F'$. Then for all $x,y\in X$, $Rxy$ if and only if $R'f(x)f(y)$. As a corollary, $f:\F\cong\F'\rsto{\mathsf{ran}(f)}$.
\end{fact}

We recall the definition and some basic results about ultraproducts of general frames (c.f. \cite{Kracht1999l}). Given a family $\mathcal{X}=(X_i)_{i\in I}$ of sets, the {\em direct product $\prod_{i\in I}X_i$ of $\mathcal{X}$} is defined as follows:
$$
    \prod_{i\in I}X_i=\cset{f:I\to\bigcup_{i\in I}X_i: \forall{i\in I}(f(i)\in X_i)}.
$$

Let $I$ be a set and $F\sub\mathcal{P}(I)$. Then $F$ is called a {\em filter over $I$} if (i) $I\in F$; (ii) for all $J,K\in F$, $J\cap K\in F$ and (iii) for all $J\in F$ and $K\in\mathcal{P}(I)$, $J\sub K$ implies $K\in F$.
A filter $F$ over $I$ is called {\em proper} if $\ve\not\in F$. Moreover, $F$ is called an {\em ultrafilter} if $F$ is proper and for all $J\in\mathcal{P}(I)$, exactly one of $J\in F$ and $J^c\in F$ holds.

\begin{definition}
    Let $(\gf_i=(X_i,R_i,A_i))_{i\in I}$ be a family of general frames and $U$ an ultrafilter over $I$. Then the {\em ultraproduct $\prod_U\gf_i=(\prod_UX_i,R,A)$ of $(\gf_i)_{i\in I}$} is defined as follows:
    \begin{itemize}
        \item $\prod_UX_i=\cset{[x]:x\in\prod_{i\in I}X_i}$, where 
        \begin{center}
          $[x]=\cset{y\in\prod_{i\in I}X_i:\cset{i\in I:x(i)=y(i)}\in U}$;
        \end{center}
        \item the set $A$ of internal sets is $\cset{[P]:P\in\prod_{i\in I}A_i}$, where
        \begin{center}
            $[P]=\cset{[x]\in \prod_UX:\cset{i\in I:x(i)\in P(i)}\in U}$.
        \end{center}
        \item for all $x,y\in\prod_{i\in I}X_i$, $R[x][y]$ if and only if $\cset{i\in I:R_ix(i)y(i)}\in U$.
    \end{itemize}
\end{definition}

\begin{proposition}\label{prop:UltraproductOfGF}
    Let $(\gf_i=(X_i,R_i,A_i))_{i\in I}$ be a family of general frames, $U$ an ultrafilter over $I$ and $\prod_U\gf_i$ the ultraproduct of $(\gf_i)_{i\in I}$. Then for all $P,Q\in\prod_{i\in I}A_i$
    \begin{align*}
        ([P])^c&=[\tup{P(i)^c:i\in I}],
        &[P]\cap[Q]&=[\tup{P(i)\cap Q(i):i\in I}],\\
        {R}[[P]]&=[\tup{{R_i}[P(i)]:i\in I}],
        &\breve{R}[[P]]&=[\tup{\breve{R_i}[P(i)]:i\in I}].
    \end{align*}
    As a corollary, $\prod_U\gf_i$ is a general frame. Moreover, if $\gf_i$ is differentiated (tight) for each $i\in I$, then $\prod_U\gf_i$ is differentiated (tight).
\end{proposition}
\begin{proof}
    The fact that $\prod_U\gf_i$ is a general frame follows from \cite[Proposition 5.7.1]{Kracht1999l}. Suppose $\gf_i$ is differentiated for each $i\in I$. Take any distinct $[x],[y]\in X$. Then $J=\cset{j\in I:x(j)\neq y(j)}\in U$. For each $j\in J$, there exists $P_j\in A_j$ such that $x(j)\in P_j$ and $y(j)\not\in P_j$. Let $[P]\in A$ be such that $P(j)=P_j$ for all $j\in J$. Then clearly, $[x]\in[P]$ and $[y]\not\in[P]$, which entails $\prod_U\gf_i$ is differentiated. The case for tightness is similar.
\end{proof}

\begin{theorem}\label{thm:UltraproductOfGF}
    Let $(\M_i=(\gf_i,V_i))_{i\in I}$ be a family of models, $V=\tup{V_i:i\in I}$ and $U$ an ultrafilter over $I$. Let $[V]$ be the valuation in $\prod_U\gf_i$ such that $[V]:p\mapsto[\tup{V_i(p):i\in I}]$. Then for all $x=\tup{x_i:i\in I}$ and $\phi\in\L$,
    \begin{enumerate}[(1)]
        \item $\prod_U\gf_i,[V],[x]\md\phi$ if and only if $\cset{i\in I:\M_i,x_i\md\phi}\in U$.
        \item $\bigoplus_{i\in I}\gf_i\md\phi$ implies $\prod_U\gf_i\md\phi$.
    \end{enumerate}
\end{theorem}
\begin{proof}
    (1) follows from \cite[Theorem 5.7.2]{Kracht1999l}. For (2), suppose $\bigoplus_{i\in I}\gf_i\md\phi$. Take any $y=\tup{y_i:i\in I}$ and valuation $V'$ in $\prod_U\gf_i$. Note that for all $j\in\omega$, $V'(p_j)=[P_j]$ for some $P_j\in\prod_{i\in I}A_i$. For all $i\in I$, let $V'_i$ be the valuation in $\gf_i$ such that $V'_i(p_j)=P_j(i)$ for all $j\in\omega$. Then we see $V'=[\tup{V'_i:i\in I}]$. Since $\bigoplus_{i\in I}\gf_i\md\phi$, we have $\cset{i\in I:\gf_i,V'_i,y_i\md\phi}=I\in U$. Then by (1),  $\prod_U\gf_i,V',[y]\md\phi$. 
\end{proof}

\begin{remark}
    It is good to note that an ultraproduct of frames is not always a frame. Here is an example. Let $\mathsf{grz}=\B(\B(p\to\B p)\to p)\to p$. For each $n\in\mathbb{Z}^+$, let $\C_n$ be the reflexive-transitive chain of $n$ elements, that is, $\C_n=(n,\geq)$. Let $U$ be an ultrafilter over $\omega$ and $\gf=\prod_U\C_n$. Since $\C_n\md\mathsf{grz}$ for each $n\in\mathbb{Z}^+$, by Theorem~\ref{thm:UltraproductOfGF}(2), $\gf\md\mathsf{grz}$. However, it is not hard to verify that $\kappa\gf$ contains an infinite ascending chain, which entails that $\kappa\gf\not\md\mathsf{grz}$. Hence $\kappa\gf\neq\gf$.
\end{remark}

\begin{definition}
A (normal) {\em tense logic} is a set of formulas $L$ such that (i) all instances of classical propositional tautologies belong to $L$;
(ii) $\bd \phi\to \psi\in L$ if and only if $\phi\to \B\psi\in L$;
(iii) if $\phi,\phi\to\psi\in L$, then $\psi\in L$; (iv) if $\phi\in L$, then $\phi^s\in L$ for every substitution $s$. The least tense logic is denoted by $\mathsf{K}_t$.\footnote{Since tense logics are closed under arbitrary intersections, the least tense logic exists. The logic $\mathsf{K}_t$ is also known as the normal bi-modal logic with extra axioms $p\to\B\bd p$ and $p\to\bb\D p$, see \cite{Prior1967}.} 
\end{definition}

For every tense logic $L$ and set of formulas $\Sigma$, let $L\oplus\Sigma$ denote the smallest tense logic containing $L\cup\Sigma$.  A tense logic $L_1$ is a {\em sublogic} of $L_2$ (or $L_2$ is an {\em extension} of $L_1$) if $L_1\sub L_2$. $L_2$ is called a {\em proper extension of $L_1$} if $L_2\supsetneq L_1$. Let $\NExt(L)$ be the set of all extensions of $L$. A tense logic $L$ is {\em consistent} if $\bot\not\in L$. The only inconsistent tense logic is $\L$. Clearly $(\NExt(L), \cap,\oplus)$ is a complete lattice with top $\L$ and bottom $L$. 

\begin{fact}\label{fact:alg-completeness}
    For all tense logic $L$, $L=\mathsf{Log}(\mathsf{DF}(L))=\mathsf{Log}(\mathsf{RF}_r(L))$.
\end{fact}
\begin{proof}
    By \cite[Theorems~17 and 67]{Venema2007}, $L=\mathsf{Log}(\mathsf{DF}(L))$. Note that every descriptive frame is a disjoint union of rooted refined frames, we see $L=\mathsf{Log}(\mathsf{RF}_r(L))$.
\end{proof}

\begin{definition}
    Let $L$ be a tense logic. Then we say (i) $L$ is {\em Kripke complete} if $L=\mathsf{Log}(\mathsf{Fr}(L))$; (ii) $L$ has the {\em finite model property (FMP)} (notation: $\mathrm{FMP}(L)$) if $L=\mathsf{Log}(\mathsf{Fin}(L))$; (iii) $L$ is {\em tabular} if $L = \mathsf{Log}(\F)$ for some finite frame $\F$; (iv) $L$ is {\em pretabular} if $L$ is not tabular and all consistent proper extension of $L$ is tabular.
    Let $\mathsf{TAB}$ and $\mathsf{PTAB}$ denote the sets of all tabular and pretabular tense logics, respectively.
\end{definition}

Recall some results on tabularity of tense logics (c.f. \cite{Chen.Ma2024a}). For each $n\in\omega$ and $\phi\in\L$, we define the formula $\Delta^n\phi$ by:
\begin{center}
    $\Delta^0\phi=\phi$ and $\Delta^{k+1}\phi=\Delta^k\phi\vee\D\Delta^k\phi\vee\bd\Delta^k\phi$.
\end{center}
Let $\nabla^n\phi:=\neg\Delta^n\neg\phi$. For each $n\in\omega$, the formula $\mathbf{tab}^T_n$ is defined as 
\begin{align*}
\mathbf{tab}^T_n=\neg(\Delta^n\psi_0\wedge\cdots\wedge\Delta^n\psi_n),
\end{align*}
where $\psi_i=\neg p_0\wedge\cdots\wedge\neg p_{i-1}\wedge p_i$ for each $i\leq n$. Note that $\psi_0=p_0$.

For example, $\Delta p = p\vee\D p\vee\bd p$, $\Delta^2q=q\vee\D q\vee\bd q\vee\D^2 q\vee\D\bd q\vee\bd\D q\vee\bd^2 q$ and $\mathbf{tab}^T_1=\neg((p_0\vee\D p_0\vee\bd p_0)\wedge((\neg p_0\wedge p_1)\vee\D(\neg p_0\wedge p_1)\vee\bd(\neg p_0\wedge p_1)))$.
Semantically, one can check that for all model $\M=(X,R,V)$ and $x\in X$, 
\begin{center}
    $\M,x\md\Delta^n\phi$ if and only if $\M,y\md\phi$ for some $y\in \R^n[x]$.
\end{center}

\begin{fact}\label{fact:tabn}
    Let $\gf=(X,R,A)\in\mathsf{GF}$ and $x\in X$. Then for every $n\in\omega$,
    \begin{enumerate}[(1)]
        \item $\R^{n}[x]=\R^{n+1}[x]$ if and only if $\R^n[x]=\R^\omega[x]$.
        \item $\kappa\gf,x\md\mathbf{tab}^T_n$ if and only if $|\R^n[x]|\leq n$.
    \end{enumerate}
\end{fact}
\begin{proof}
    See \cite[Lemma 3.4 and Lemma 3.5]{Chen.Ma2024a}.
\end{proof}

\begin{theorem}[{\cite[Theorem 3.7]{Chen.Ma2024a}}]\label{thm:tabular}
    For every consistent tense logic $L\in\NExt(\mathsf{K}_t)$, $L\in\mathsf{TAB}$ if and only if $\mathbf{tab}^T_n\in L$ for some $n\geq 1$.
\end{theorem}

\begin{theorem}\label{thm:nontabular-infrootedframe}
    Let $L$ be a non-tabular tense logic. Then $\mathsf{Fin}_r(L)\neq\mathsf{RF}_r(L)$, i.e., there exists an infinite rooted refined frame $\gf$ such that $\gf\md L$.
\end{theorem}
\begin{proof}
    Since $L$ is non-tabular, by Theorem~\ref{thm:tabular}, $\mathbf{tab}^T_n\not\in L$ for any $n\in\omega$. By Fact~\ref{fact:alg-completeness}, $L=\mathsf{Log}(\mathsf{RF}_r(L))$. Thus for all $n\in\omega$, there exists $\gf_n=(X_n,R_n,A_n)\in\mathsf{RF}_r(L)$ and $x_n\in X_n$ such that $\gf_n,x_n\not\md\mathbf{tab}^T_n$. By Fact~\ref{fact:tabn}(2), we see $\gf_n,x_n\not\md\mathbf{tab}^T_m$ for all $m\geq n$. By Theorem~\ref{thm:UltraproductOfGF}(1), $\prod_U\gf_n,[x]\not\md\mathbf{tab}^T_n$ for any $n\in\omega$, where $U$ is a non-principal ultrafilter over $\omega$ and $x:n\mapsto x_n$ for all $n\in\omega$. Let $\gg=(\prod_U\gf_n)_{[x]}$. By Theorem~\ref{thm:UltraproductOfGF}, $\prod_U\gf_n\md L$, which implies $\gg\md L$. Note that $\gg\in\mathsf{RF}_r$ and $\gg,[x]\not\md\mathbf{tab}^T_n$ for any $n\in\omega$, by Fact~\ref{fact:tabn}(2), we see that $\gg\in\mathsf{RF}_r(L)\setminus\mathsf{Fin}_r(L)$. 
\end{proof}

\begin{corollary}
    If $L\in\mathsf{PTAB}$, then $L=\mathsf{Log}(\gf)$ for some rooted refined frame $\gf$.
\end{corollary}

\section{Generalized Jankov formulas and local t-morphisms}

Jankov formulas are widely used in the study of lattices of intermediate logics and modal logics (see \cite{Jankov1963,deJongh1968,Kracht1993,Rautenberg1980,Wolter1997}). Fine \cite{Fine1974} developed frame formulas for finite rooted $\mathsf{S4}$-frames, which are similar to Jankov formulas for finite subdirectly irreducible Heyting algebras. Such a formula for $\F$ is refuted by some $\mathsf{S4}$-frame $\G$ if and only if $\F$ is a p-morphic image of a generated subframe of $\G$. In this section, we introduce the generalized Jankov formulas for image-finite Kripke frames and local t-morphisms, which will be the main tools we use in the following sections. Given any rooted image-finite frame $\F=(X,R)$, $x\in X$ and $k\in\mathbb{Z}^+$, we associate the generalized Jankov formula $\J^k(\F,x)$ to it. We show that for any frame $\G=(Y,S)$ and $y\in Y$, $(\G,y)$ validates $\neg\J^k(\F,x)$ if and only if there exists no $k$-t-morphism $f:(\G,y)\to^k(\F,x)$. This generalizes the Jankov formulas for finite algebras or finite rooted frames defined in \cite{deJongh1968,Fine1974}.

\begin{definition}
    Let $\gf=(X,R,A)\in\mathsf{GF}$. Then $\gf$ is said to be {\em image-finite} if $|\R[x]|<\aleph_0$ for all $x\in X$.\footnote{$\aleph_0$ is the least infinite cardinal.} A tense logic is call {\em image-finite} if every refined frame of $L$ is image-finite.
\end{definition}

\begin{example}
    Consider the finitely alternative logics $T_{n,m}=\mathsf{K}_t\oplus\cset{\mathbf{alt}^+_n,\mathbf{alt}^-_m}$ defined in \cite{Ma.Chen2021}, where the formulas $\mathbf{alt}^+_n$ and $\mathbf{alt}^-_n$ are defined as follows:
    \begin{align*}
        \tag*{($\mathbf{alt}^+_n$)} \B p_0 \vee \B(p_0 \to p_1) \vee \cdots \vee \B(p_0 \wedge \cdots \wedge p_{n-1} \to p_{n})\\
        \tag*{($\mathbf{alt}^-_n$)} \bb p_0 \vee \bb(p_0 \to p_1) \vee \cdots \vee \bb(p_0 \wedge \cdots \wedge p_{n-1} \to p_{n})
    \end{align*}
    Refined frames $\gf=(X,R,A)$ for $T_{n,m}$ are exactly those such that $|R[x]|\leq n$ and $|\breve{R}[x]|\leq m$ for all $x\in X$. Clearly, these logics are all image-finite.
\end{example}

\begin{lemma}\label{lem:pointwise-finite-kappaF}
    Let $\gf=(X,R,A)\in\mathsf{GF}$ be image-finite. Then 
    \begin{enumerate}[(1)]
        \item $|\R^n[x]|<\aleph_0$ for all $x\in X$ and $n\in\omega$.
        \item If $\gf$ is differentiated, then $\mathsf{Log}(\gf)=\mathsf{Log}(\kappa\gf)$.
    \end{enumerate}
\end{lemma}
\begin{proof}
    For (1), we prove by induction on $n$. The case $n=0$ is trivial. Suppose $n>0$. Then by induction hypothesis, we have 
    \begin{center}
        $|\R^n[x]|\leq\sum_{y\in R^{n-1}_\sharp[x]}|\R[y]|<\aleph_0$.
    \end{center}
    For (2), it suffices to show $\mathsf{Log}(\gf)\sub\mathsf{Log}(\kappa\gf)$. Take any $\phi\not\in\mathsf{Log}(\kappa\gf)$. Then there exists $x\in X$ and a valuation $V$ in $\kappa\gf$ such that $\kappa\gf,V,x\not\md\phi$. Let $d(\phi)=m$. By~(1), $\R^m[x]$ is finite. Since $\gf$ is differentiated, there exists a valuation $V'$ in $\gf$ such that $V(p)\cap \R^m[x]=V'(p)\cap \R^m[x]$ for all $p\in\mathsf{Prop}$. Thus $(\kappa\gf,V)\rsto \R^m[x]\cong(\gf,V')\rsto \R^m[x]$. Note that $d(\phi)=m$, we have $\gf,V',x\not\md\phi$. Thus $\phi\not\in\mathsf{Log}(\gf)$.
\end{proof}

\begin{definition}
    Let $\F=(X,R)$ be a image-finite frame and $x\in X$. Let $k\in\mathbb{Z}^+$ and $\tup{x_i:i\in n}$ be an enumeration of $\R^{k}[x]$ where $x=x_0$. Then the formula $\J^k(\F,x)$ is defined to be the conjunction of the following formulas:
    \begin{enumerate}[(1)]
        \item $p_0\wedge\nabla^k(p_0\vee\cdots\vee p_{n-1})$
        \item $\nabla^k(p_i\to\neg p_j)$, for all $i\neq j$
        \item $\nabla^{k-1}((p_i\to\D p_j)\wedge(p_j\to\bd p_i))$, for all $Rx_ix_j$
        \item $\nabla^{k-1}((p_i\to\neg\D p_j)\wedge(p_j\to\neg\bd p_i))$, for all $x_j\not\in R[x_i]$
    \end{enumerate}
    $\J^k(\F,x)$ is called the {\em Jankov formula of $(\F,x)$ with degree $k$}.
\end{definition}

By Lemma~\ref{lem:pointwise-finite-kappaF}(1), $\R^k[x]$ is finite for all $k\in\omega$. Thus $\J^k(\F,x)$ is well defined. Let $V$ be a valuation in $\F$ such that $V(p_i)=\cset{x_i}$ for all $i\in n$. Then it is not hard to check that $\F,V,x\md\J^k(\F,x)$. Thus $\F,x\not\md\neg\J^k(\F,x)$ for any $k\in\mathbb{Z}^+$. 

Intuitively, $\J^k(\F,x)$ is a formula that describes the subframe $\F\rsto\R^k[x]$ of an image-finite frame $\F$, while the classical Jankov formula describes its corresponding finite frames. If $X=\R^k[x]$, then $\F$ is finite and the formula $\J^k(\F,x)$ captures all the information about $\F$. In this sense, our formulas generalize the Jankov formulas.

\begin{definition}
    Let $\F=(X,R)$ and $\F'=(X',R')$ be frames, $x\in X$ and $x'\in X'$. For each $k\in\mathbb{Z}^+$, a partial function $f:X\to X'$ is called a {\em $k$-t-morphism from $(\F,x)$ to $(\F',x')$ (notation: $f:(\F,x)\to^k(\F',x')$)}, if $\mathsf{dom}(f)\supseteq\R^k[x]$, $f(x)=x'$ and
    \begin{center}
        for all $y\in \R^{k-1}[x]$, $f[R[y]]= R'[f(y)]$ and $f[\breve{R}[y]]=\breve{R'}[f(y)]$.
    \end{center} 
    A $k$-t-morphism is also called a {\em local t-morphism}. Moreover, $f:(\F,x)\to^k(\F',x')$ is said to be {\em full} if $\mathsf{ran}(f)\supseteq(R'_\sharp)^\omega[x']$.
\end{definition}

Intuitively, $f:(\F,x)\to^k(\F',x')$ is a `t-morphism for points in $\R^{k-1}[x]$'. Clearly, $f:\F\to\F'$ if $\R^{k-1}[x]=X$ and $f$ is full. It is also not hard to verify that if $f:\F\to\F'$ and $f(x)=x'$, then $f:(\F,x)\to^k(\F',x')$ for all $k\in\mathbb{Z}^+$.

\begin{lemma}\label{lem:k-t-morphism}
    Let $\F=(X,R)$ and $\F'=(X',R')$ be frames, $x\in X$ and $x'\in X'$. Suppose $k\in\mathbb{Z}^+$ and $f:(\F,x)\to^{k}(\F',x')$. Then for all formula $\phi$ with $md(\phi)\leq k$,
    \begin{center}
        $\F,x\md\phi$ implies $\F',x'\md\phi$.
    \end{center}
\end{lemma}
\begin{proof}
    Suppose $\F',V',x'\not\md\phi$. Let $V$ be a valuation in $\F$ such that $V(p)=f^{-1}[V'(p)]$ for all $p\in\mathsf{Prop}$. It suffices to prove the following claim:
    
    \noindent\textbf{Claim:} For all $l\leq k$, $y\in\R^{k-l}[x]$ and  $\psi\in\L$ with $\mathsf{md}(\psi)\leq l$
    \begin{center}
        $\F,V,y\md\psi$ if and only if $\F',V',f(y)\md\psi$.
    \end{center}
    \textbf{Proof of Claim:} By induction on $l$. The case $l=0$ follows from the definition of $V$ immediately. Let $l>0$. By induction hypothesis, $y$ and $f(y)$ agree on all formulas with modal depth at most $l-1$. Take any $\gamma\in\L$ with $\mathsf{md}(\gamma)\leq l-1$. Suppose $\F,V,y\md\bd\gamma$. Then $\F,V,z\md\gamma$ for some $z\in\breve{R}[y]$. Note that $z\in\R^{k-(l-1)}[x]$, by induction hypothesis, $\F',V',f(z)\md\gamma$. Since $f(z)\in f[\breve{R}[y]]=\breve{R'}[f(y)]$, we see $\F',V',f(y)\md\bd\gamma$. Suppose $\F',V',f(y)\md\bd\gamma$. Then $\F',V',z'\md\gamma$ for some $z'\in\breve{R'}[f(y)]$. Since $\breve{R'}[f(y)]=f[\breve{R}[y]]$, there exists $z\in\breve{R}[y]$ such that $f(z)=z'$. Note that $z\in\R^{k-(l-1)}[x]$, by induction hypothesis, $\F,V,z\md\gamma$, which entails $\F,V,y\md\bd\gamma$. Since $f[R[y]]= R'[f(y)]$, by a similar proof, we see $\F,V,y\md\B\gamma$ if and only if $\F',V',f(y)\md\B\gamma$. \hfill$\dashv$

    By the claim above, take $l=k$ and $\psi=\phi$, we see $\F,x\not\md\phi$.
\end{proof}

Recall that the refutation of a Jankov formula means the existence of particular p-morphisms. Similarly, the following lemma shows that the refutation of a generalized Jankov formula $\J^k(\F,x)$ amounts to the existence of particular $k$-t-morphisms.

\begin{lemma}\label{lem:JankovLemma-k}
    Let $\F=(X,R)$ be a frame, $\G=(Y,S)$ a image-finite frame, $x\in X$ and $y\in Y$. Then for all $k\in\mathbb{Z}^+$,
    \begin{center}
        $\F,x\not\md\neg\J^k(\G,y)$ if and only if there exists $f:(\F,x)\to^{k}(\G,y)$.
    \end{center}
\end{lemma}
\begin{proof}
    Suppose $f:(\F,x)\to^{k}(\G,y)$. Since $\G,y\not\md\neg\J^k(\G,y)$ and $\mathsf{md}(\J^k(\G,y))=k$, by Lemma~\ref{lem:k-t-morphism}, we have $\F,x\not\md\neg\J^k(\G,y)$. Suppose $\F,V,x\not\md\neg\J^k(\G,y)$ for some valuation $V$ in $\F$. Let $\tup{y_i:i\in n}$ be an enumeration of $\S^k[y]$ used in the definition of $\J^k(\G,y)$. We define the function $f:\R^k[x]\to Y$ as follows:
    \begin{center}
        for all $z\in\R^k[x]$, $f(z)=y_i$ if and only if $z\md p_i$.
    \end{center}
    Since $x\md\nabla^k(p_0\vee\cdots\vee p_{n-1})$ and $x\md\nabla^k(p_i\to\neg p_j)$ for all $i\neq j$, we see $\R^k[x]\sub\bigcup_{i\in n}V(p_i)$ and $\R^k[x]\cap V(p_i)\cap V(p_j)=\ve$ for all $i\neq j$. Thus $f$ is well-defined. It suffices to show that $f:(\F,x)\to^k(\G,y)$. Since $y=y_0$ and $x\md p_0$, $f(x)=y$. Take any $z\in\R^{k-1}[x]$. Then $f(z)=y_i$ and $z\md p_i$ for some $i\in n$. For all $y_j\in Y$, we have
    \begin{center}
        $y_j\in S[f(z)]\Longrightarrow x\md\nabla^{k-1}(p_i\to\D p_j)\Longrightarrow z\md\D p_j\Longrightarrow y_j\in f[R[z]]$.
    \end{center}
    and
    \begin{center}
        $y_j\not\in S[f(z)]\Longrightarrow x\md\nabla^{k-1}(p_i\to\neg\D p_j)\Longrightarrow z\md\neg\D p_j\Longrightarrow y_j\not\in f[R[z]]$.
    \end{center}
    Thus $f[R[z]]= S[f(z)]$. Similarly, $f[\breve{R}[z]]= \breve{S}[f(z)]$. Hence $(\F,x)\to^{k}(\G,y)$.
\end{proof}

Let us take a closer look at local t-morphisms. For the purpose of this paper, we are interested in full local t-morphisms. In what follows, we provide a sufficient condition for a local t-morphism $f:(\F,x)\to^k(\F',x')$ to be full, which requires only information about $f$. 

\begin{definition}
    Let $f:(\F,x)\to^{k}(\F',x')$ be a local t-morphism. A subset $Y\sub\R^{k-1}[x]$ is called {\em sufficient} if $f[\R[y]]\sub f[Y]$ for all $y\in Y$. We say that $f$ is {\em sufficient} if there is a nonempty sufficient set $Y\sub\mathsf{dom}(f)$.
\end{definition}

\begin{fact}\label{fact:boundary-sufficient}
    Let $f:(\F,x)\to^{k}(\F',x')$ and $Y\sub\R^{k-1}[x]$. Then $Y$ is sufficient if one of the following holds:
    \begin{enumerate}[(1)]
        \item for all $y\in Y$, there exist $u,v\in Y$ with $f(y)=f(u)=f(v)$ and $R[u]\cup\breve{R}[v]\sub Y$.
        \item there exists $Z\sub Y$ such that $\R[Z]\sub Y$ and $f[Y]=f[Z]$.
    \end{enumerate}
\end{fact}
\begin{proof}
    For (1), take any $y\in Y$. Then we see $f[R[y]]=R'[f(y)]=R'[f(u)]=f[R[u]]\sub f[Y]$. Similarly, $f[\breve{R}[y]]=f[\breve{R}[v]]\sub f[Y]$. Thus $f[\R[y]]\sub f[Y]$. Hence (1) holds. It is clear that (2) follows from (1).
\end{proof}

\begin{lemma}\label{lem:sufficient-full}
    Let $\F=(X,R)$ and $\F'=(X',R')$ be frames, $x\in X$ and $x'\in X'$. Suppose $k\in\mathbb{Z}^+$ and $f:(\F,x)\to^{k}(\F',x')$ is sufficient. Then $f$ is full.
\end{lemma}
\begin{proof}
    Suppose $Y\sub\mathsf{dom}(f)$ is sufficient and $y\in Y$. We show by induction on $n$ that $(R')_\sharp^n[f(y)]\sub f[Y]$ for all $n\in\omega$. The case $n=0$ is trivial. Let $n>0$. Take any $z'\in(R')_\sharp^{n-1}[f(y)]$. By induction hypothesis, $f(z)=z'$ for some $z\in Y$. Since $Y$ is sufficient, we see $(R'_\sharp)[z']=f[\R[z]]\sub f[Y]$. Thus $\bigcup_{z'\in(R')_\sharp^{n-1}[f(y)]}(R')_\sharp[z']\sub f[Y]$, which entails $(R')_\sharp^n[f(y)]\sub f[Y]$. Hence $(R')_\sharp^\omega[f(y)]\sub f[Y]\sub\mathsf{ran}(f)$.
\end{proof}

\begin{corollary}\label{coro:tmorphism-ktmorphism}
    Let $\F=(X,R)$ and $\F'=(X',R')$ be rooted frames, $x\in X$ and $x'\in X'$. Suppose $k\in\mathbb{Z}^+$, $f:(\F,x)\to^{k}(\F',x')$ and $\R^{k-1}[x]=X$. Then $f:\F\twoheadrightarrow\F'$.
\end{corollary}
\begin{proof}
    Since $f:(\F,x)\to^{k}(\F',x')$ and $\R^{k-1}[x]=X$, $f:X\to X'$ is a t-morphism. Note that $X\sub\R^{k-1}[x]$ is sufficient, by Lemma~\ref{lem:sufficient-full}, $f$ is full. Thus $f:\F\twoheadrightarrow\F'$.
\end{proof}

By Lemma~\ref{lem:JankovLemma-k} and Corollary~\ref{coro:tmorphism-ktmorphism}, we have

\begin{theorem}\label{thm:JankovLemma}
    Let $\F=(X,R)$ be a rooted frame, $\G=(Y,S)$ a rooted image-finite frame and $y\in Y$. Let $k\in\mathbb{Z}^+$. Suppose $X=\R^{k-1}[x]$ for all $x\in X$. Then
    \begin{center}
        $\F\twoheadrightarrow\G$ if and only if $\F\not\md\neg\J^k(\G,y)$.
    \end{center}
\end{theorem}

\section{Tense logics over $\mathsf{S4}_t$ with bounded parameters}

Recall that $\mathsf{S4}_t=\mathsf{K}_t\oplus\cset{\mathbf{4},\mathbf{T}}$, where $\mathbf{T}=\B p\to p$ and $\mathbf{4}=\B p\to\B\B p$. 
From this section, we focus on tense logics in $\NExt(\mathsf{S4}_t)$. Unless otherwise specified, frames are always assumed to be reflexive and transitive.

To study a lattice with a complex structure, it is useful to investigate its sublattices. For example, a better understanding of the lattice $\NExt(\mathsf{K4})$ of transitive modal logics can be obtained by studying modal logics with bounded depth and width. Segerberg \cite{Segerberg1971} proved that every modal logic of finite depth enjoys the FMP and Fine \cite{Fine1974b} showed that every modal logic of finite width is Kripke complete. 

In this section, we introduce tense logics in $\NExt(\mathsf{S4}_t)$ with bounded depth, forth-width, back-width and z-degrees. We call extensions of them tense logics with bounded parameters and investigate their basic properties. 

\begin{definition}
    Let $\gf=(X,R,A)\in\mathsf{GF}$ be a general frame, $\alpha$ an ordinal and $\mathcal{Y}=\tup{y_i\in X:i<\alpha}$. Then 
    \begin{itemize}
        \item $\mathcal{Y}$ is called a {\em chain} in $\F$ if $Rx_\lambda x_\gamma$ for all $\lambda<\gamma<\alpha$;
        \item $\mathcal{Y}$ is called a {\em strict chain} in $\F$ if it is a chain and $x_\lambda\not\in R(x_\gamma)$ for all $\lambda<\gamma<\alpha$;
        \item $\mathcal{Y}$ is called a {\em (strict) co-chain} in $\F$ if it is a (strict) chain in $\breve{\F}$;
        \item $\cset{y_i\in X:i<\alpha}$ is called an {\em anti-chain} in $\F$ if $x_\lambda\not\in R(x_\gamma)$ for all $\lambda\neq\gamma<\alpha$.
    \end{itemize}
    The {\em length} $l(\mathcal{Y})$ of a strict chain $\mathcal{Y}=\tup{y_i\in X:i<\alpha}$ is defined to be $\alpha$.
    We say that $x\in X$ is {\em of depth $n$ (notation: $\mathrm{dep}(x)=n$)}, if there exists a strict chain $\mathcal{Y}$ in $\gf\rsto R[x]$ with $l(\mathcal{Y})=n$ and there is no strict chain of greater length. Otherwise $x$ is said to be of infinite depth and we write $\mathrm{dep}(x)=\aleph_0$. We define the {\em depth $\mathrm{dep}(\gf)$ of $\gf$} by $\mathrm{dep}(\gf)=\mathrm{sup}\cset{\mathrm{dep}(x):x\in X}$. For a tense logic $L$, we define the {\em depth $\mathrm{dep}(L)$ of $L$} by $\mathrm{dep}(L)=\mathrm{sup}\cset{\mathrm{dep}(\gf):\mathsf{RF}(L)}$.

    Let $n\in\mathbb{Z}^+$ and $x\in X$. We say that $x$ is {\em of forth-width $n$ (notation: $\mathrm{wid}^+(x)=n$)}, if there exists an anti-chain $Y\sub R[x]$ with $|Y|=n$ and there is no anti-chain in $R[x]$ with greater size. Otherwise we write $\mathrm{wid}^+(\gf)=\aleph_0$. We say that $\gf$ is {\em of forth-width $n$ (notation: $\mathrm{wid}^+(\gf)=n$)}, if $\mathrm{sup}\cset{\mathrm{wid}^+(x):x\in X}=n$. Again, for a tense logic $L$, we define the {\em forth-width $\mathrm{wid}^+(L)$ of $L$} by $\mathrm{wid}^+(L)=\mathrm{sup}\cset{\mathrm{wid}^+(\gf):\mathsf{RF}(L)}$.
    {\em Back-width} is defined dually. We write $\mathrm{wid}^-(x)=n$, $\mathrm{wid}^-(\gf)=n$ and $\mathrm{wid}^-(L)=n$ if $x$, $\gf$ and $L$ is of back-width $n$, respectively.
\end{definition}

\begin{definition}
    Let $\gf=(X,R,A)\in\mathsf{GF}$ be a frame and $x\in X$. Let $k\in\mathbb{Z}^+$. Then we say {\em $x$ is of z-degree $k$ (notation: $\mathrm{zdg}(x)=k$)}, if $\R^{k-1}[x]\neq \R^{k}[x]=\R^\omega[x]$. Specially, $\mathrm{zdg}(x)=0$ if $\R^\omega[x]=\cset{x}$ and $\mathrm{zdg}(x)=\aleph_0$ if $\R^k[x]\neq\R^{k+1}$ for any $k\in\omega$. We define the {\em z-degree $\mathrm{zdg}(\gf)$ of $\gf$} by $\mathrm{zdg}(\gf)=\mathrm{sup}\cset{\mathrm{zdg}(x):x\in X}$. For a tense logic $L$, we define the {\em z-degree $\mathrm{zdg}(L)$ of $L$} by $\mathrm{zdg}(L)=\mathrm{sup}\cset{\mathrm{zdg}(\gf):\gf\in\mathsf{RF}(L)}$.
\end{definition}

\begin{definition}
    For each $n\in\mathbb{Z}^+$, let $(\mathbf{bz}_n)$, $(\mathbf{bw}^+_n)$ and $(\mathbf{bw}^-_n)$ denote the following formulas respectively:
    \begin{align*}
        \tag{$\mathbf{bz}_n$} \Delta^{n+1}p&\to\Delta^{n}p\\
        \tag{$\mathbf{bw}^+_n$} \bigwedge_{i\leq n}\D p_i&\to\bigvee_{i\neq j\leq n}\D(p_i\wedge(p_j\vee\D p_j))\\
        \tag{$\mathbf{bw}^-_n$} \bigwedge_{i\leq n}\bd p_i&\to\bigvee_{i\neq j\leq n}\bd(p_i\wedge(p_j\vee\bd p_j))
    \end{align*}
    Moreover, we define the formula $(\mathbf{bd}_n)$ for each $n\in\mathbb{Z}^+$ as follow:
    \begin{align*}
        \mathbf{bd}_1 &= \D\B p_0\to p_0\\
        \mathbf{bd}_{k+1} &= \D(\B p_{k}\wedge\neg\mathbf{bd}_k)\to p_{k}
    \end{align*}
    Specially, we define $\mathbf{bd}_\omega=\mathbf{bz}_\omega=\mathbf{bw}^+_\omega=\mathbf{bw}^-_\omega=\top$.
    Let $k,l,n,m\in\mathbb{Z}^*$. We define the tense logics $\mathsf{S4BP}^{k,l}_{n,m}$ by
    \begin{center}
        $\mathsf{S4BP}^{k,l}_{n,m}=\mathsf{S4}_t\oplus\cset{\mathbf{bd}_k,\mathbf{bz}_l,\mathbf{bw}^+_n,\mathbf{bw}^-_m}$.
    \end{center}
    We say a tense logic $L$ is {\em with bounded parameters} if $L\supseteq\mathsf{S4BP}^{k,l}_{n,m}$ for some $k,l,n,m\in\mathbb{Z}^*$ such that $\min\cset{k,l,n,m}<\omega$.
\end{definition}

Lots of well-studied logics are instants of tense logics with bounded parameters. For example, the tense logic $\mathsf{S4.3}_t$ ($\mathsf{Lin}_t\oplus\mathsf{Dens}_1$ in \cite{Wolter1996b}) of linear frames turns out to be equal to $\mathsf{S4BP}^{\omega,1}_{1,1}$. Moreover, the following fact holds.

\begin{fact}\label{fact:bounds}
    Let $\gf=(X,R,A)\in\mathsf{RF}$, $x\in X$ and $n\in\mathbb{Z}^+$. Then
    \begin{enumerate}[(1)]
        \item $\gf,x\md\mathbf{bd}_n$ if and only if $\mathrm{dep}(x)\leq n$.
        \item $\gf,x\md\mathbf{bw}^+_n$ if and only if $\mathrm{wid}^+(x)\leq n$.
        \item $\gf,x\md\mathbf{bw}^-_n$ if and only if $\mathrm{wid}^-(x)\leq n$.
        \item $\gf,x\md\mathbf{bz}_n$ if and only if $\mathrm{zdg}(x)\leq n$.
    \end{enumerate}
\end{fact}
\begin{proof}
    The proofs of (1) - (3) are analogous to those of \cite[Propositions 3.42 and 3.44]{Chagrov.Zakharyaschev1997}, while the proof of (4) is similar to that of \cite[Corollary 3.35]{Chagrov.Zakharyaschev1997}.
\end{proof}

\begin{definition}
    Let $\gf=(X,R,A)\in\mathsf{GF}$ be and $x\in W$. The {\em cluster generated by $x$}, denoted by $C(x)$, is defined as follows:
    \[
    C(x)=R[x]\cap \breve{R}[x]
    \]
    A subset $C\sub X$ is called a {\em cluster in $\gf$} if $C=C(x)$ for some $x\in X$.
    Let $n\in\mathbb{Z}^+$. We say $\gf$ is {\em of the girth $n$ (notation: $\mathrm{gir}(\gf)=n$)} if there exists a cluster $C$ in $\gf$ such that $|C|=n$ and there is no cluster in $\gf$ of larger size. We write $\mathrm{gir}(\gf)=\aleph_0$ if for all $k\in\mathbb{Z}^+$, there exists a cluster $C$ in $\gf$ such that $|C|>k$.
\end{definition}

\begin{lemma}\label{lem:cluster-n-formula}
    Let $\M=(X,R,V)$ be a model and $C\sub X$ be a cluster in $\F$. Suppose $n\in\omega$, $x,y\in C$ and $x\equiv_{\mathsf{Prop}(n)} y$. Then $x\equiv_{\L(n)} y$.
\end{lemma}
\begin{proof}
    The proof proceeds by induction on the complexity of $\phi$. The case $\phi\in\mathsf{Prop}(n)$ follows from $x\equiv_{\mathsf{Prop}(n)} y$ immediately and the Boolean cases are standard. Consider the case $\phi=\B\psi$. Suppose $\M,x\md\phi$. Then $\M,z\md\psi$ for all $z\in R[y]$. Since $C$ is a cluster and $x,y\in C$, we see $R[x]= R[y]$. By induction hypothesis, $\M,z\md\psi$ for all $z\in R[y]$. Thus $\M,y\md\phi$. Symmetrically, $\M,y\md\phi$ implies $\M,x\md\phi$. Note that $\breve{R}[x]=\breve{R}[y]$, the proof for the case $\phi=\bd\psi$ is similar.
\end{proof}

\begin{proposition}\label{prop:finitecluster}
    Let $\gf=(X,R,A)\in\mathsf{RF}$ be finitely generated. Then $\mathrm{gir}(\gf)<\aleph_0$. 
\end{proposition}
\begin{proof}
    Let $\gf$ be generated by $U_0,\cdots,U_{n-1}\in A$ for some $n\in\omega$. Let $V$ be a valuation in $\gf$ such that $V(p_i)=U_i$ for all $i<n$. Then $A=\cset{V(\phi):\phi\in\L(n)}$. Let $\M=(\gf,V)$. Since $\gf$ is differentiated, we have $x\not\equiv_{\L(n)}y$ for any different $x,y\in X$. Take any cluster $C$ in $\gf$. By Lemma~\ref{lem:cluster-n-formula}, $x\not\equiv_{\mathsf{Prop}(n)}y$ for any different $x,y\in X$. Thus $|C|\leq 2^n<\aleph_0$. Since $C$ is chosen arbitrarily, $\mathrm{gir}(\gf)\leq 2^n<\aleph_0$.
\end{proof}

\begin{lemma}\label{lem:nstep-fin}
    Let $n,m,k\in\mathbb{Z}^+$ and $\gf=(X,R,A)$ be a finitely generated refined frame for $\mathsf{S4BP}^{k,\omega}_{n,m}$. Then $\gf$ is image-finite.
\end{lemma}
\begin{proof}
    Let $x\in X$. By Fact~\ref{fact:bounds}, $R[x]$ and $\breve{R}[x]$ contains only finitely many distinct chains and each of them is of finite depth. Since $\gf$ is finitely generated, by Proposition~\ref{prop:finitecluster}, every cluster in $\gf$ is finite. Thus $\R[x]=\cset{x}\cup R[x]\cup\breve{R}[x]$ is finite.
\end{proof}

\begin{theorem}\label{thm:BP-KC-FMP}
    Let $n,m,k,l\in\mathbb{Z}^+$ and $L\in\NExt(\mathsf{S4BP}^{k,\omega}_{n,m})$. Then $L$ is Kripke complete. Moreover, if $L\in\NExt(\mathsf{S4BP}^{k,l}_{n,m})$, then $L$ has the finite model property.
\end{theorem}
\begin{proof}
    Let $\mathcal{K}$ be the class of all rooted finitely generated refined frames of $L$. By Lemma~\ref{lem:nstep-fin}, every frame in $\mathcal{K}$ is image-finite. By Lemma~\ref{lem:pointwise-finite-kappaF}(2), we see $L=\mathsf{Log}(\mathcal{K})=\mathsf{Log}(\cset{\kappa\gf:\gf\in\mathcal{K}})$, which entails that $L$ is Kripke complete. 
    
    Suppose $L\in\NExt(\mathsf{S4BP}^{k,l}_{n,m})$. Take any $\gf=(X,R,A)\in\mathcal{K}$ and $x\in X$. Since $\gf$ is rooted, by Lemma~\ref{lem:pointwise-finite-kappaF}(1) and Fact~\ref{fact:bounds}, we see that $X=\R^l[x]$ is finite. Thus every frame in $\mathcal{K}$ is finite, which entails $L$ has the FMP immediately.
\end{proof}

\begin{remark}\label{rem:K4BP}
    The readers can see that reflexivity of frames plays no role in the proofs above. In fact, it is also natural to drop axiom $\mathbf{T}$ and define tense logics $\mathsf{K4BP}^{k,l}_{n,m}$ with bounded parameters. Basic properties above can be easily generalized.
\end{remark}

\section{Pretabular tense logics over $\NExt(\mathsf{S4BP}^{k,l}_{n,m})$}\label{sec:S4BP}

In this section, we study the pretabular logics in $\NExt(\mathsf{S4BP}^{k,l}_{n,m})$ where $k,l,n,m\in\mathbb{Z}^+$. We obtain a full characterization of pretabular logics in $\NExt(\mathsf{S4BP}^{k,l}_{n,m})$. By Theorem~\ref{thm:BP-KC-FMP}, every logic in $\NExt(\mathsf{S4BP}^{k,l}_{n,m})$ is Kripke complete. %Thus in this section, we consider only Kripke frame. 
It turns out that a tense logic $L\in\NExt(\mathsf{S4BP}^{k,l}_{n,m})$ is pretabular if and only if $L=\mathsf{Log}(\F)$ for some rooted frame $\F$ with certain conditions.

\begin{definition}
    Let $\F=(X,R)$ be a frame. Let $W^S=\cset{C(x):x\in X}$ and $R^S=\cset{\tup{C(x),C(y)}:Rxy}$. Then $\F^S=(W^S,R^S)$ is called the {\em skeleton of $\F$}. We call $\F$ a skeleton if $\F\iso\F^S$.
\end{definition}

\begin{definition}
    Let $\F=(X,R)$ be a frame. Then for each $\lambda\leq\omega$ and $x\in X$, we define the frame $\F^x_\lambda=(X^x_\lambda,R^x_\lambda)$ as follows:
    \begin{itemize}
        \item $X^x_\lambda=X\uplus N^x_\lambda$, where $N^x_\lambda=\cset{x_i:0<i\leq\lambda}$;
        \item $R^x_\lambda=R\cup(C^x_\lambda\times R[x])\cup(\breve{R}[x]\times C^x_\lambda)\cup(C^x_\lambda\times C^x_\lambda)$, where $C^x_\lambda=N^x_\lambda\cup\cset{x}$.
    \end{itemize}
    A frame $\G$ is called a {\em $\lambda$-pre-skeleton} if $\G\iso\F^x_\lambda$ for some skeleton $\F$ and $\lambda>0$.
\end{definition}
Intuitively, $\F^x_\lambda$ is the frame obtained from $\F$ by replacing one reflexive point $x$ in $\F$ by a cluster with $1+\lambda$ points. The readers can readily verify that pre-skeletons are those frames containing exactly 1 proper cluster. It should be clear that $\F$ and $\F^x_\lambda$ share the same skeleton and have the same z-degree, width and depth. In what follows, without loss of generality, we always assume that $X\cap N^x_\lambda=\ve$.

\begin{example}
    Concrete examples of skeletons and pre-skeletons is provided in Figure~\ref{fig:ske-preske}, where $\F$ is a skeleton and the others are pre-skeletons.
    \begin{figure}        
    \begin{center}
            \begin{tikzpicture}
                \draw (-2.5,0.5) node{$\F$};
                \draw (-.5,0) node{$\circ$};
                \draw (-.5,0) node[right]{\small $x_4$};
                
                \draw (-1,1) node{$\circ$};
                \draw [->] (-.55,.05) -- (-0.95,.95);
                \draw (-2,1) node{$\circ$};
                \draw (-1.5,0) node{$\circ$};
                \draw [->] (-1.55,.05) -- (-1.95,.95);
                \draw [->] (-1.45,.05) -- (-1.05,.95);
                \draw (-1,1) node[right]{\small $x_3$};
                \draw (-1.5,0) node[right]{\small $x_1$};
                \draw (-2,1) node[right]{\small $x_0$};
            \end{tikzpicture}
            \qquad\qquad
            \begin{tikzpicture}
                \draw (-2.5,0.5) node{$\F^{x_4}_{6}$};
                \draw (-.5,0) node{$7$};
                \draw (-.5,0) ellipse (.2 and .2);
                
                \draw (-1,1) node{$\circ$};
                \draw [->] (-.6,.2) -- (-0.95,.95);
                \draw (-2,1) node{$\circ$};
                \draw (-1.5,0) node{$\circ$};
                \draw [->] (-1.55,.05) -- (-1.95,.95);
                \draw [->] (-1.45,.05) -- (-1.05,.95);
                \draw (-1,1) node[right]{\small $x_3$};
                \draw (-1.5,0) node[right]{\small $x_1$};
                \draw (-2,1) node[right]{\small $x_0$};
            \end{tikzpicture}
            \qquad\qquad
            \begin{tikzpicture}
                \draw (-2.5,0.5) node{$\F^{x_4}_{\omega}$};
                \draw (-.5,0) node{$\omega$};
                \draw (-.5,0) ellipse (.2 and .2);
                
                \draw (-1,1) node{$\circ$};
                \draw [->] (-.6,.2) -- (-0.95,.95);
                \draw (-2,1) node{$\circ$};
                \draw (-1.5,0) node{$\circ$};
                \draw [->] (-1.55,.05) -- (-1.95,.95);
                \draw [->] (-1.45,.05) -- (-1.05,.95);
                \draw (-1,1) node[right]{\small $x_3$};
                \draw (-1.5,0) node[right]{\small $x_1$};
                \draw (-2,1) node[right]{\small $x_0$};
            \end{tikzpicture}
    \end{center}
        \caption{Skeleton and pre-skeleton}
        \label{fig:ske-preske}
    \end{figure}
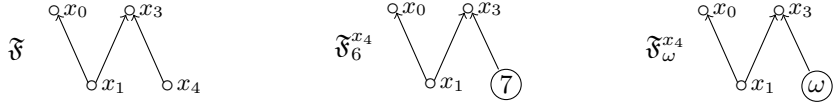
\end{example}

\begin{lemma}\label{lem:cluster-clopen-S4BP}
    Let $\gf=(X,R,A)\in\mathsf{RF}_r(\mathsf{S4BP}^{k,l}_{n,m})$ for some $k,l,n,m\in\mathbb{Z}^+$. Then 
    \begin{enumerate}[(1)]
        \item $\kappa\gf^S$ is finite.
        \item $C\in A$ for every cluster $C\sub W$.
    \end{enumerate}    
\end{lemma}
\begin{proof}
    For (1), by Fact~\ref{fact:bounds}, $\kappa\gf\in\mathsf{RF}_r(\mathsf{S4BP}^{k,l}_{n,m})$. Note that $\kappa\gf$ contains no proper cluster, $\gf$ is image-finite and (1) follows from Lemma~\ref{lem:pointwise-finite-kappaF}(1). For (2), let $\tup{C_i}_{i\leq s}$ be an enumeration of clusters in $\gf$ such that $C=C_s$. Let $c\in C$. For each $i<s$, take a point $c_i\in C_i$. Suppose $c_i\not\in R[c]$. Since $\gf$ is tight, $c_i\in U_i$ and $c\not\in\D U_i$ for some $U_i\in A$. Since $\gf$ is differentiated, $c_i\in U_i'$ and $c\not\in U_i'$ for some $U_i'\in A$. Let $V_i=(U_i\cap U_i')\cup\breve{R}[U_i\cap U_i']$. Then $c_i\in V_i$ and $c\not\in V_i$. Note that $V_i\sub\breve{R}[V_i]$, we have $C_i\sub V_i$ and $C\cap V_i=\ve$. Suppose $c_i\in R[c]$. Since $c_i\not\in C$, we see $c_i\not\in\breve{R}[c]$. By a similar argument, there exists $V_i\in A$ with $C_i\sub V_i$ and $C\cap V_i=\ve$. It is not hard to see that $C=X\setminus\bigcup_{i<s}V_i\in A$.
\end{proof}

\begin{lemma}\label{lem:pre-skeleton-inftofin}
    Let $\F=(X,R)$ be a frame and $x\in X$. Suppose $\gf=(\F^x_\lambda,A)\in\mathsf{RF}_r(\mathsf{S4BP}^{k,l}_{n,m})$ for some $k,l,n,m\in\mathbb{Z}^+$ and $\lambda\geq\aleph_0$. Then $\gf\twoheadrightarrow\F^x_s$ for all $s\in\omega$.
\end{lemma}
\begin{proof}
    Let $s\in\omega$ and $C^x_\lambda$ denote the cluster in $\gf$ generated by $x$. By Lemma~\ref{lem:cluster-clopen-S4BP}, $C^x_\lambda\in A$. Note that $C^x_\lambda$ is infinite and $\gf$ is differentiated, there are pairwise disjoint $U_0,\cdots,U_{s}\in A$ such that $\bigcup_{i\leq s}U_i=C^x_\lambda$. Let $D=\cset{d_0,\cdots,d_{s}}=C^x_s$ denote the cluster in $\F^x_s$ generated by $x$. We define the map $f:X^x_\lambda\to X^x_s$ by
    \begin{align*}
        f(x)=
        \begin{cases}
            d_i &\text{ if }x\in U_i,\\
            x &\text{ otherwise.}
        \end{cases}
    \end{align*}
    It is easy to check that $f:\gf\twoheadrightarrow\F^x_s$.
\end{proof}

\begin{definition}
    Let $\M=(X,R,V)$ be a model and $C\sub X$ a cluster. We call $D\sub C$ a {\em $\mathsf{Prop}(n)$-approximation of $C$} if for all $c\in C$, $c\equiv_{\mathsf{Prop}(n)}d$ for some $d\in D$.
\end{definition}

\begin{lemma}\label{lem:selection-in-cluster}
    Let $\M=(X,R,V)$ be a model and $C\sub X$ be a cluster in $\F$. Let $Y$ be a subset of $X$ such that $X\setminus C\sub Y$ and $Y\cap C$ is a $\mathsf{Prop}(n)$-approximation of $C$. Then for all $y\in Y$ and $\phi\in\L(n)$, 
    \begin{center}
        $\M,y\md\phi$ if and only if $\M\rsto Y,y\md\phi$.
    \end{center}
\end{lemma}
\begin{proof}
    Let $\M\rsto Y=(Y,S)$. The proof proceeds by induction on the complexity of $\phi$. The case $\phi\in\mathsf{Prop}(n)$ is trivial and the Boolean cases are standard. Consider the case $\phi=\B\psi$. Suppose $\M,y\md\phi$. Then $\M,z\md\psi$ for all $z\in R[y]$. Since $S[y]\sub R[y]$, by induction hypothesis, $\M\rsto Y,z\md\psi$ for all $z\in S[y]$. Thus $\M\rsto Y,y\md\phi$. Suppose $\M,y\not\md\phi$. Then $\M,z\not\md\psi$ for some $z\in R[y]$. Assume $z\not\in C$. Since $X\setminus C\sub Y$, we see $z\in Y$. By induction hypothesis, $\M\rsto Y,z\not\md\psi$ and so $\M\rsto Y,y\not\md\phi$. Assume $z\in C$. Since $Y\cap C$ is a $\mathsf{Prop}(n)$-approximation of $C$, there is $z'\in C\cap Y$ such that $z\equiv_{\mathsf{Prop}(n)}z'$. By Lemma~\ref{lem:cluster-n-formula}, $\M,z'\not\md\psi$. By induction hypothesis, $\M\rsto Y,z'\not\md\psi$ and so $\M\rsto Y,y\not\md\phi$. The case $\phi=\bd\psi$ can be proved in a similar way.
\end{proof}

\begin{lemma}\label{lem:ThFomega=CapThFn}
    Let $\F=(X,R)$ be a skeleton. Then for all $x\in X$,
    \begin{center}
        $\mathsf{Log}(\F^x_\omega)=\bigcap\limits_{n\in\omega}\mathsf{Log}(\F^x_n)$.
    \end{center}
\end{lemma}
\begin{proof}
    Note that $\F^x_\omega\twoheadrightarrow\F^x_n$ for each $n\in\omega$, we have $\mathsf{Log}(\F^x_\omega)\sub\bigcap_{n\in\omega}\mathsf{Log}(\F^x_n)$. Suppose $\phi(p_0,\cdots,p_{m-1})\not\in\mathsf{Log}(\F^x_\omega)$. Then $\F^x_\omega,V\not\md\phi$ for some valuation $V$ in $\F^x_\omega$. Let $\M=(\F^x_\omega,V)$. Since $\mathsf{Prop}(m)$ is finite, there exists a finite $\mathsf{Prop}(m)$-approximation $C$ of $C^x_\omega$. Let $Y=X_{-x}\cup C$. Clearly, $\F^x_{|C|+1}\iso\F^x_\omega\rsto Y$. Note that $X^x_\omega\setminus C^x_\omega\sub Y$, by Lemma~\ref{lem:selection-in-cluster}, $\M\rsto Y\not\md\phi$, which entails $\phi\not\in\mathsf{Log}(\F^x_{|C|+1})$ and so $\phi\not\in\bigcap_{n\in\omega}\mathsf{Log}(\F^x_n)$.
\end{proof}

\begin{lemma}\label{lem:Finr-FMP}
    Let $L$ be a tense logic with $\mathbf{bz}_k\in L$ and $\cset{\F_i:i\in I}$ a family of rooted frames. Suppose $L=\bigcap\cset{\mathsf{Log}(\F_i):i\in I}$. Then $\mathsf{Fin}_{r}(L)\sub\bigcup_{i\in I}\mathsf{TM}(\F_i)$.
\end{lemma}
\begin{proof}
    Let $L=\mathsf{Log}(\F^x_\omega)$, $\G\in\mathsf{Fin}_r(L)$ and $y\in\G$. Note that $\G\not\md\neg\J^{k+1}(\G,y)$, we see $\neg\J^{k+1}(\G,y)\not\in L$. Since $L=\bigcap\cset{\mathsf{Log}(\F_i):i\in I}$, $\F_i\not\md\neg\J^{k+1}(\G,y)$ for some $i\in I$. By Theorem~\ref{thm:JankovLemma}, $\F_i\to^{k+1}\G$. Since $\mathbf{bz}_k\in L$ and $\F_i\md L$, we see $\mathsf{zdg}(\F_i)\leq k$. By Corollary~\ref{coro:tmorphism-ktmorphism}, we have $\F_i\twoheadrightarrow\G$. Hence $\G\in\mathsf{TM}(\F_i)$.
\end{proof}

\begin{lemma}\label{lem:finiteK4-ontopmorphism-cluster}
    Let $\F=(X,R)$ and $\G=(Y,S)$ be frames of finite depth and $x\in X$. Suppose $f:\F\twoheadrightarrow\G$ and $|C(f(x))|\geq 2$. Then 
    \begin{enumerate}[(1)]
        \item $|C(y)|\geq 2$ for some $y\in R[x]$.
        \item $|C(y')|\geq 2$ for some $y'\in\breve{R}[x]$.
    \end{enumerate}
    Moreover, if $\F=(\F')^x_n$ is a $n$-pre-skeleton for some $n\geq 1$, then 
    \begin{enumerate}
        \item[(3)] for all $z\in X$, $C(z)$ is proper if and only if $C(f(z))$ is proper.
        \item[(4)] $f[C(x)]=C(f(x))$.
    \end{enumerate}
\end{lemma}
\begin{proof}
    For (1), let $f(x)=y_0$ and $\mathrm{dep}(\F)=k$. Since $|C(y_0)|\geq 2$, there exists $y_1\in C(y_0)$ with $y_0\neq y_1$. By $f$ is a t-morphism and $Sf(x)y_1$, there exists $x_1\in X$ such that $f(x_1)=y_1$ and $Rxx_1$. Again, since $Sf(x_1)y_0$, there exists $x_2\in X$ such that $f(x_2)=y$ and $Rx_1x_2$. By repeating this construction, we get an $R$-chain $\tup{x_i\in X:i\leq 2k+1}$ such that $x=x_0$, $f(x_{2i})=y_0$ and $f(x_{2i+1})=y_1$ for all $i\leq k$. Since $\mathrm{dep}(\F)=k<2k+1$, we see $x_{m}=x_{l}$ for some $m<l\leq 2k+1$. By transitivity of $\F$, we see $x_{m+1}Rx_l=x_m$ and so $C(x_m)=C(x_{m+1})$. Note that $f(x_m)\neq f(x_{m+1})$, we have $x_m\neq x_{m+1}$ and so $|C(x_m)|\geq 2$. 
    (2) can be proved symmetrically.

    For (3), suppose $C(f(z))$ is proper. Then $|C(f(z))|\geq 2$. By (1) and (2), there are $z_1\in R[z]$ and $z_2\in\breve{R}[z]$ such that $|C(z_1)|\geq 2$ and $|C(z_2)|\geq 2$. Since $\F$ is a pre-skeleton, there is exactly 1 proper cluster $C$ in $\F$. Thus $C(z_1)=C(z_2)=C$. Since $z_2RzRz_2$, we have $z\in C$ and so $C(z)=C$, which implies $C(z)$ is proper. Suppose $C(z)$ is proper. Then we have $z\in C(x)$, which entails $C(f(z))=C(f(x))$ is proper.
    
    For (4), it is clear that $f(y)\in C(f(x))$ for all $y\in C(x)$, which entails $f[C(x)]\sub C(f(x))$. Let $u\in C(f(x))$. Since $f$ is onto, $f(z)=u$ for some $u\in X$. Since $C(u)=C(f(x))$ is proper, by (3), $C(z)$ is proper. Note that $C(x)$ is the only proper cluster in $\F$, $z\in C(x)$ and so $u\in f[C(x)]$. Thus $f[C(x)]=C(f(x))$.
\end{proof}

\begin{lemma}\label{lem:preskeleton-tmorphism}
    Let $\F=(X,R)$ and $\G=(Y,S)$ be finite skeletons. Then the following are equivalent:
    \begin{enumerate}[(1)]
        \item $\F^x_n\twoheadrightarrow\G^y_m$ for some $1\leq m\leq n\leq\omega$;
        \item $\F^x_k\twoheadrightarrow\G^y_l$ for all $l\leq k\leq\omega$.
    \end{enumerate}
\end{lemma}
\begin{proof}
    Clearly (1) follows from (2). Suppose $f:\F^x_n\twoheadrightarrow\G^y_m$. Let $X_0=X\setminus C(x)$ and $Y_0=Y\setminus C(y)$. By Lemma~\ref{lem:finiteK4-ontopmorphism-cluster}(3-4), $f[C^x_n]=C^y_m$ and $f[X_0]=Y_0$. For any $k\geq l$, there exists $g:C^x_k\twoheadrightarrow C^y_l$. Let $f'=f\rsto X_0$ and $h=g\cup f'$. Clearly, $h:\F^x_k\twoheadrightarrow\G^y_l$.
\end{proof}

\begin{definition}
    Let $\lambda>0$. Then a $\lambda$-pre-skeleton $\F^x_\lambda$ is called c-irreducible if
    \begin{center}
        $\mathsf{TM}(\F^x_\lambda)=\mathsf{IM}(\cset{\F^x_m:0<m\leq\lambda})\cup\mathsf{TM}(\F)$.
    \end{center}
    Otherwise, $\F^x_\lambda$ is c-reducible.
\end{definition}

Intuitively, a pre-skeleton $\F^x_\lambda$ is c-irreducible if and only if $\F^x_\lambda$ does not admit any proper t-morphic image with the same girth. 

\begin{example}
    Consider the frames $\F_1$ and $\F_2$ in Figure~\ref{fig:example-c-irreducible}. We see that $\F_1$ is c-reducible, since $\F_2$ is a t-morphic image of $\F_1$. However, $\F_2$ is c-irreducible, since for every non-injective t-morphism $f$, the $f$-image is always a skeleton.
    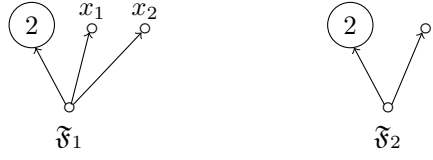
\begin{figure}[ht]
        \[
        \begin{tikzpicture}[scale=1]
            \draw (0.5,-1.2) node[below]{$\F_1$};
            \draw (0,0) node{$2$};
            \draw (0,0) ellipse (.3 and .3);
            \draw [->] (0.52,-1.05) -- (0.76,-.1);
            \draw (0.5,-1.1) node{$\circ$};
            \draw [->] (0.45,-1.05) -- (.05,-.3);
            \draw (.8,-.05) node{$\circ$};
            \draw (.8,-.05) node[above]{$x_1$};
            \draw [->] (0.55,-1.05) -- (1.45,-.1);
            \draw (1.5,-.05) node{$\circ$};
            \draw (1.5,-.05) node[above]{$x_2$};
        \end{tikzpicture}
        \quad\quad\quad\quad\quad\quad
        \begin{tikzpicture}[scale=1]
            \draw (0.5,-1.2) node[below]{$\F_2$};
            \draw (0,0) node{$2$};
            \draw (0,0) ellipse (.3 and .3);
            \draw [->] (0.55,-1.05) -- (0.95,-.1);
            \draw (0.5,-1.1) node{$\circ$};
            \draw [->] (0.45,-1.05) -- (.05,-.3);
            \draw (1,-.05) node{$\circ$};
        \end{tikzpicture}
        \]
        \caption{A c-reducible frame and a c-irreducible frame}
        \label{fig:example-c-irreducible}
        \end{figure}
\end{example}

\begin{lemma}\label{lem:c-irreducible-2-omega}
    Let $\F=(X,R)$ be a skeleton and $x\in X$ a reflexive point. Then 
    \begin{center}
        $\F^x_1$ is c-irreducible if and only if $\F^x_\omega$ is c-irreducible.
    \end{center}
\end{lemma}
\begin{proof}
    The right-to-left direction follows from Lemma~\ref{lem:finiteK4-ontopmorphism-cluster}(4) immediately. For the left-to-right direction, take any and $\H'\in\mathsf{TM}(\F^x_\omega)$. Let $h:\F^x_\omega\to\H'$. By Lemma~\ref{lem:finiteK4-ontopmorphism-cluster}, $\H'$ is of the form $\H^z_m$ for some $m\leq\omega$.
    If $m=0$, then $\H'\in\mathsf{TM}(\F)$ by Lemma~\ref{lem:preskeleton-tmorphism}. Suppose $m>0$. By Lemma~\ref{lem:preskeleton-tmorphism}, $\H^z_1\in\mathsf{TM}(\F^x_1)$. Note that $\mathsf{TM}(\F^x_1)=\mathsf{IM}(\cset{\F^x_1})\cup\mathsf{TM}(\F)$ and every frame in $\mathsf{TM}(\F)$ contains no proper cluster, $\H^z_1\iso\F^x_1$. Thus $\H'\iso\H^z_m\iso\F^x_m$.
\end{proof}

\begin{lemma}\label{lem:pretabularity-fin-skeleton}
    Let $\F=(X,R)$ be a rooted finite skeleton and $x\in X$. Then 
    \begin{center}
        $\mathsf{Log}(\F^x_\omega)$ is pretabular if and only if $\F^x_1$ is c-irreducible.
    \end{center}
\end{lemma}
\begin{proof}
    Suppose $\F^x_1$ is c-reducible. Then there exists a frame $\G_0$ such that $\G_0\not\iso\F^x_1$ and $\G_0\in\mathsf{TM}(\F^x_1)\setminus\mathsf{TM}(\F)$. Let $f:\F^x_1\twoheadrightarrow\G_0$. Note that $C^x_1$ is the unique proper cluster in $\F^x_1$.    
    We claim $|f[C^x_1]|>1$. Suppose $|f[C^x_1]|=1$. Then $f[C^x_1]=\cset{y_0}$ for some $y_0\in\G_0$. Then we define $f':\F\to\G_0$ by $f'(x)=y_0$ and $f'(z)=f(z)$ for all $z\in X\setminus\cset{x}$. It is clear that $f':\F\twoheadrightarrow\G_0$, which entails $\G_0\in\mathsf{TM}(\F)$ and leads to a contradiction. Thus $f[C^x_1]=C(f(x))$ is a proper cluster in $\G_0$. By Lemma~\ref{lem:finiteK4-ontopmorphism-cluster}(3), $\G_0$ is a 2-pre-skeleton and so $\G_0\cong\G^y_1$ for some $\G=(Y,S)$. It suffices to show $\mathsf{Log}(\G^y_\omega)\supsetneq\mathsf{Log}(\F^x_\omega)$. Since $\F^x_1\twoheadrightarrow\G^y_1$, by Lemma~\ref{lem:preskeleton-tmorphism}, $\F^x_\omega\twoheadrightarrow\G^y_\omega$ and so $\mathsf{Log}(\G^y_\omega)\supseteq\mathsf{Log}(\F^x_\omega)$. Consider the formula $\phi=\neg\J^k(\F^x_1,x)$ where $k=\mathsf{zdg}(\F)+1$. We show that $\phi\in\mathsf{Log}(\G^y_\omega)\setminus\mathsf{Log}(\F^x_\omega)$. Clearly, $\phi\not\in\mathsf{Log}(\F^x_\omega)$. Suppose $\phi\not\in\mathsf{Log(\G^y_\omega)}$. By Lemma~\ref{lem:ThFomega=CapThFn}, $\G^y_m\not\md\phi$ for some $m\in\omega$. By Theorem~\ref{thm:JankovLemma}, $\G^y_m\twoheadrightarrow\F^x_1$. By Lemma~\ref{lem:preskeleton-tmorphism}, $\G_0=\G^y_1\twoheadrightarrow\F^x_1$. Thus $|\F^x_1|=|\G_0|$ and so $f:\F^x_1\twoheadrightarrow\G_0$ is injective. By Fact~\ref{fact:inj-tmorphism}, $\G_0\iso\F^x_1$, which contradicts the assumption.

    Suppose $\F^x_1$ is c-irreducible. By Lemma~\ref{lem:c-irreducible-2-omega}, $\F^x_\omega$ is c-irreducible. Take any $L\supsetneq\mathsf{Log}(\F^x_\omega)$. By Theorem~\ref{thm:BP-KC-FMP}, $L=\mathsf{Log}(\mathsf{Fin}_r(L))$. By Lemma Lemma~\ref{lem:Finr-FMP}, we see
    \begin{center}
        $\mathsf{Fin}_r(L)\subsetneq\mathsf{Fin}_r(\mathsf{Log}(\F^x_\omega))\sub\mathsf{TM}(\F^x_\omega)=\mathsf{IM}(\cset{\F^x_{m}:0<m\leq\omega})\cup\mathsf{TM}(\F)$. 
    \end{center}
    Thus there exists $n\in\mathbb{Z}^+$ such that $\mathsf{Fin}_r(L)\sub\mathsf{IM}(\cset{\F^x_{m}:0<m\leq n})\cup\mathsf{TM}(\F)$, which entails that $L=\mathsf{Log}(\mathsf{Fin}_r(L))$ is tabular. Clearly, $\mathsf{Log}(\F^x_\omega)$ is non-tabular. Hence $\mathsf{Log}(\F^x_\omega)$ is pretabular.
\end{proof}

\begin{theorem}\label{thm:pretabular-S4BP-ch}
    Let $L\in\NExt(\mathsf{S4BP}^{k,l}_{n,m})$ for some $k,l,n,m\in\mathbb{Z}^+$. Then $L$ is pretabular if and only if $L=\mathsf{Log}(\F^x_\omega)$ for some c-irreducible rooted pre-skeleton $\F^x_\omega$.
\end{theorem}
\begin{proof}
    The right-to-left direction follows from Lemma~\ref{lem:pretabularity-fin-skeleton} immediately. For the other direction, suppose $L$ is pretabular. By Theorem~\ref{thm:nontabular-infrootedframe}, $L=\mathsf{Log}(\gf')$ for some infinite rooted refined frame $\gf'$. Since $\gf'\md\mathsf{S4BP}^{k,l}_{n,m}$, we see that $\kappa\gf'$ is finite and so $\kappa\gf'\iso\F^x_\lambda$ for some frame $\F^x_\lambda$ and $\lambda\geq\aleph_0$. By Lemma~\ref{lem:pre-skeleton-inftofin}, $\gf'\twoheadrightarrow\F^x_n$ for all $n\in\mathbb{Z}^+$. By Lemma~\ref{lem:ThFomega=CapThFn}, 
    \begin{center}
        $\mathsf{Log}(\F^x_\lambda)=\mathsf{Log}(\kappa\gf')\sub\mathsf{Log}(\gf')\sub\bigcap_{0<i<\omega}\mathsf{Log}(\F^x_i)=\mathsf{Log}(\F^x_\lambda)=\mathsf{Log}(\F^x_\omega)$.
    \end{center}
    Since $L=\mathsf{Log}(\gf')$ is pretabular, $L=\mathsf{Log}(\F^x_\omega)$. Let $\G=(\F^x_\omega)^S$. Consider the frame $\G^{C(x)}_\omega$. Note that $\G^{C(x)}_\omega$ is isomorphic to the frame obtained by collapsing all proper cluster except $C(x)$. Then $\F^x_\omega\twoheadrightarrow\G^{C(x)}_\omega$. Since $L$ is pretabular, $L=\mathsf{Log}(\G^{C(x)}_\omega)$. Note that $\G$ is a finite skeleton, by Lemma~\ref{lem:pretabularity-fin-skeleton}, $\G^{C(x)}_\omega$ is c-irreducible.
\end{proof}

\begin{lemma}\label{lem:irreducible-frame-same-logic-implies-iso}
    Let $\F^x_\omega$ and $\G^y_\omega$ be finite c-irreducible rooted pre-skeletons. Then 
    \begin{center}
        $\mathsf{Log}(\F^x_\omega)=\mathsf{Log}(\G^y_\omega)$ implies $\F^x_1\iso\G^y_1$.
    \end{center}    
\end{lemma}
\begin{proof}
    Let $k=\mathsf{zdg}(\F)+\mathsf{zdg}(\G)$. Since $\F^x_1\md\mathsf{Log}(\F^x_\omega)$ and $\F^x_1,x\not\md\neg\J^k(\F^x_1,x)$, we see $\neg\J^k(\F^x_1,x)\not\in\mathsf{Log}(\G^y_\omega)$ and so $\G^y_\omega\not\md\neg\J^k(\F^x_1,x)$. By Theorem~\ref{lem:JankovLemma-k} and Corollary~\ref{coro:tmorphism-ktmorphism}, $\G^y_\omega\twoheadrightarrow\F^x_1$. By Lemma~\ref{lem:preskeleton-tmorphism}, $\G^y_1\twoheadrightarrow\F^x_1$. Symmetrically, we can show that $\F^x_1\twoheadrightarrow\G^y_1$. Then $|\G^y_1|=|\F^x_1|$. By Fact~\ref{fact:inj-tmorphism}, $\F^x_1\iso\G^y_1$.
\end{proof}

\begin{theorem}
    For all $k,l,n,m\in\mathbb{Z}^+$, $|\mathsf{PTAB}(\mathsf{S4BP}^{k,l}_{n,m})|<\aleph_0$.
\end{theorem}
\begin{proof}
    Take any skeleton $\F\in\mathsf{Fr}_r(\mathsf{S4BP}^{k,l}_{n,m})$ and $x\in\F$. It is not hard to see that $|\F|\leq|\R^l[x]|<{(k(m+n))}^l$. Thus there are only finitely many skeletons validating $\mathsf{S4BP}^{k,l}_{n,m}$. By Theorem~\ref{thm:pretabular-S4BP-ch}, $|\mathsf{PTAB}(\mathsf{S4BP}^{k,l}_{n,m})|<\aleph_0$.
\end{proof}

\begin{theorem}
    For all $k,l,n,m\in\mathbb{Z}^+$ and $L\in\mathsf{PTAB}(\mathsf{S4BP}^{k,l}_{n,m})$, $L$ has the FMP.
\end{theorem}
\begin{proof}
    Follows from Lemma~\ref{lem:ThFomega=CapThFn} and Theorem~\ref{thm:pretabular-S4BP-ch} immediately.
\end{proof}

\section{Pretabular tense logics in $\NExt(\mathsf{S4.3}_t)$}

In Section~\ref{sec:S4BP}, pretabular tense logics of finite width, depth and z-degree are studied and a full characterization is obtained. In the present section, we turn our attention to tense logics with weaker restrictions on these parameters. As a representative case, we consider the tense logic $\mathsf{S4.3}_t$ of linear frames and its extensions. Recall that $\mathsf{S4.3}_t=\mathsf{S4BP}_{1,1}^{\omega,1}$, which indicates that it has an very strong bounds on width and z-degree, while imposing no restriction on depth.

Linearity is important in the research of tense logics. The tense logic $\mathsf{S4.3}_t$ and its modal fragment $\mathsf{S4.3}$ have been well-studied. Bull \cite{Bull1966} showed that every extension of $\mathsf{S4.3}$ has the FMP. By the characterization of pretabular modal logics over $\mathsf{S4}$ in \cite{Esakia.Meskhi1977,Maksimova1975}, we see $|\mathsf{PTAB}(\mathsf{S4.3})|=3$. Wolter \cite{Wolter1996b} proved that every extension of $\mathsf{S4.3}_t$ is finitely axiomatizable. The aim of this section is to show the following theorem:

\begin{theorem}
    There are exactly 5 pretabular tense logics in $\NExt(\mathsf{S4.3}_t)$.
\end{theorem}

Let $L^\ua$ denote the tense logic of all finite chains. To be precise, recall that $\C_n=(n,\geq)$ for each $n\in\mathbb{Z}^+$. We define $L^\ua=\bigcap_{n=1}^\omega\mathsf{Log}(\C_n)$. We first show that $L^\ua$ is the only pretabular logic in $\NExt(\mathsf{S4.3}_t)$ with infinite depth.

\begin{lemma}\label{lem:S4.3-infdep-morphism}
    Let $\gf\in\mathsf{RF}_r(\mathsf{S4.3}_t)$ and $n\in\mathbb{Z}^+$. If $\mathrm{dep}(\gf)>n$, then $\C_{n+1}\in\mathsf{TM}(\gf)$.
\end{lemma}
\begin{proof}
    Let $\gf=(X,R,A)$. Since $\mathrm{dep}(\gf)>n$, by Fact~\ref{fact:bounds}(1), $\gf\not\md\mathbf{bd}_n$. Then there exists a valuation $V$ and a co-chain $\tup{x_i:i\leq n}$ in $\gf$ such that $\gf,V,x_i\md\B p_i$ and $\gf,V,x_{i+1}\md\neg p_i$ for all $i<n$. We define the function $f:X\to n+1$ by:
    \begin{align*}
        f(x)=
        \begin{cases}
            i &\text{ if }x\in V(\B p_i)\setminus V(\B p_{i+1})\text{ and }i<n\\
            n &\text{ otherwise.}
        \end{cases}
    \end{align*}
    Clearly, $f^{-1}[D]\in P$ for all $D\sub n+1$. Take any $j<n$ and $x\in V(\B p_j)$. Note that $\gf$ is linear, by $x_{j+1}\md\B p_{j+1}\wedge\neg p_j$, we see $x\in R[x_{j+1}]$ and so $x\in V(\B p_{j+1})$. Thus $V(\B p_i)\subsetneq V(\B p_{i+1})$ for all $i<n$. Hence $f:\gf\twoheadrightarrow\C_{n+1}$.
\end{proof}

\begin{theorem}\label{thm:infchain-Lua}
    Let $L\in\NExt(\mathsf{S4.3}_t)$. Then $L\sub L^\ua$ if and only if $\mathrm{dep}(L)=\aleph_0$.
\end{theorem}
\begin{proof}
    The left-to-right direction is trivial. Suppose $\mathrm{dep}(L)=\aleph_0$. Then for each $n\in\mathbb{Z}^+$, there exists $\gf_n\in\mathsf{RF}_r(L)$ such that $\mathrm{dep}(\gf_n)>n$. By Lemma~\ref{lem:S4.3-infdep-morphism}, we have $\gf_n\twoheadrightarrow\C_{n+1}$, which entails $\C_{n+1}\md L$. Thus $L\sub\bigcap_{n=1}^\omega\mathsf{Log}(\C_n)=L^\ua$.
\end{proof}

\begin{lemma}\label{lem:grz.3-pretabular}
    $L^\ua$ is pretabular.
\end{lemma}
\begin{proof}
    Let $L\supsetneq L^\ua$. Note that $\mathbf{bz}_1\in L^\ua$, by Lemma~\ref{lem:Finr-FMP}, $\mathsf{Fin}_r(L^\ua)=\bigcup_{n\in\mathbb{Z}^+}\mathsf{TM}(\C_n)$. By Theorem~\ref{thm:infchain-Lua} and Fact~\ref{fact:bounds}(1), $\mathbf{bd}_n\in L$ for some $n\in\mathbb{Z}^+$. Thus $\mathsf{Fin}_r(L)=\mathsf{TM}(\cset{\C_n})$ for some $n\in\mathbb{Z}^+$. By Theorem~\ref{thm:BP-KC-FMP}, $L=\mathsf{Log}(\mathsf{Fin}_r(L))$, which entails $L=\mathsf{Log}(\C_n)$ is tabular. Hence $L^\ua\in\mathsf{PTAB}(\mathsf{S4.3}_t)$.
\end{proof}

Consider the lattice $\NExt(\mathsf{S4.3}_t)$. By Theorem~\ref{thm:infchain-Lua}, $L^\ua$ is the maximal logic with infinite depth. Thus $L^\ua$ is the only pretabular logic with infinite depth. To characterize pretabular logics with finite depth, we introduce some auxiliary definitions.

\begin{definition}
    Let $\mathsf{tp}=\cset{\pm,+,-,\circ}$. For all $\lambda$ such that $\lambda\leq\omega$, we define $\C^\circ_\lambda=(\C_1)^0_\lambda$, $\C^+_\lambda=(\C_2)^1_\lambda$, $\C^-_\lambda=(\C_2)^0_\lambda$ and $\C^\pm_\lambda=(\C_3)^1_\lambda$.
    For all $t\in\mathsf{tp}$, let $L^t=\mathsf{Log}(\C^t_\omega)$.
\end{definition}

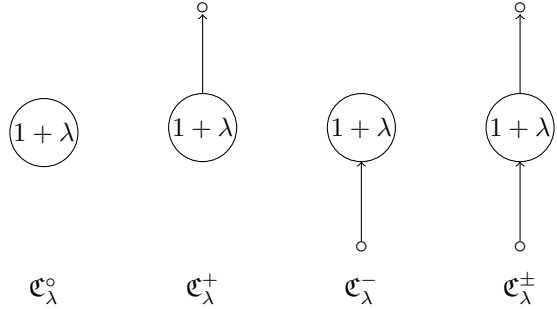
\begin{figure}[ht]
\[
\begin{tikzpicture}[scale=1.5]
    \draw (0,-1.2) node[below]{$\C^\circ_\lambda$};
    \draw (0,0) node{$1+\lambda$};
    \draw (0,0) ellipse (.3 and .3);
\end{tikzpicture}
\quad\quad\quad
\begin{tikzpicture}[scale=1.5]
    \draw (0,-1.2) node[below]{$\C^+_\lambda$};
    \draw (0,0) node{$1+\lambda$};
    \draw (0,0) ellipse (.3 and .3);
    \draw [->] (0,.3) -- (0,1);
    \draw (0,1.05) node{$\circ$};
\end{tikzpicture}
\quad\quad\quad
\begin{tikzpicture}[scale=1.5]
    \draw (0,-1.2) node[below]{$\C^-_\lambda$};
    \draw (0,0) node{$1+\lambda$};
    \draw (0,0) ellipse (.3 and .3);
    \draw [->] (0,-1) -- (0,-.3);
    \draw (0,-1.05) node{$\circ$};
\end{tikzpicture}
\quad\quad\quad
\begin{tikzpicture}[scale=1.5]
    \draw (0,-1.2) node[below]{$\C^\pm_\lambda$};
    \draw (0,0) node{$1+\lambda$};
    \draw (0,0) ellipse (.3 and .3);
    \draw [->] (0,.3) -- (0,1);
    \draw (0,1.05) node{$\circ$};
    \draw [->] (0,-1) -- (0,-.3);
    \draw (0,-1.05) node{$\circ$};
\end{tikzpicture}
\]
\caption{Frames ${\C^\circ_\lambda}$, $\C^+_\lambda$, $\C^-_\lambda$ and $\C^\pm_\lambda$}
\label{fig:pretabular-S4.3}
\end{figure}
The frames ${\C^\circ_\lambda}$, $\C^+_\lambda$, $\C^-_\lambda$ and $\C^\pm_\lambda$ are depicted in Figure.\ref{fig:pretabular-S4.3}.

\begin{proposition}\label{prop:fap-ptab-S4.3}
    $L^\circ=\mathsf{S5}_t=\mathsf{S4}_t\oplus(\D p\to\B\D p)$, $L^+=\mathsf{S4}_t\oplus\cset{\mathbf{bd}_2,\B\D p\to\D\B p}$, $L^-=\mathsf{S4}_t\oplus\cset{\mathbf{bd}_2,\bb\bd p\to\bd\bb p}$ and $L^\pm=\mathsf{S4}_t\oplus\cset{\mathbf{bd}_3,\B\D p\to\D\B p,\bb\bd p\to\bd\bb p}$.
\end{proposition}
\begin{proof}
    The proof is standard.
\end{proof}

In what follows, we show that the logics $\mathsf{PTAB}(\mathsf{S4.3}_t)=\cset{L^\pm,L^+,L^-,L^\circ,L^\ua}$.

\begin{lemma}\label{lem:pretab-S43-fin-depth}
    For all $t\in\mathsf{tp}$, $L^t$ is pretabular.
\end{lemma}
\begin{proof}
    Since $\C^t_1$ is c-irreducible, by Lemma~\ref{lem:pretabularity-fin-skeleton}, $L^t=\mathsf{Log}(\C^t_\omega)$ is pretabular.
\end{proof}

\begin{lemma}\label{lem:4-pretab-S43-fin-depth}
    Let $t,s\in\mathsf{tp}$. Then $t=s$ if and only if $L^t=L^s$.
\end{lemma}
\begin{proof}
    Follows from Lemma~\ref{lem:irreducible-frame-same-logic-implies-iso} immediately.
\end{proof}

\begin{lemma}\label{lem:c-irreducible-S4.3-finite-depth}
    Let $\F=(X,R)\in\mathsf{Fr}_r(\mathsf{S4.3}_t)$ be a skeleton of finite depth. Then 
    \begin{center}
        for all $x\in X$, $\F^x_\omega$ is c-irreducible if and only if $\F^x_\omega\in\mathsf{IM}(\cset{\C^t_\omega:t\in\mathsf{tp}})$.
    \end{center}
\end{lemma}
\begin{proof}
    The right-to-left direction is trivial. For the other direction, suppose $\F^x_\omega$ is c-irreducible. Let $Y_0=R[x]\setminus C(x)$ and $Y_1=\breve{R}[x]\setminus C(x)$. We claim that $|Y_0|<2$ and $|Y_1|<2$. Suppose $|Y_0|\geq 2$. Then let $\G=(X_0,R_0)$ where $X_0=Y_1\cup C(x)\cup\cset{y}$ and $R_0$ is the transitive-reflexive closure of $(Z\times C(x))\cup (C(x)\times\cset{y})$. It is obvious that $\G\in\mathsf{TM}(\F^x_\omega)\setminus\mathsf{IM}(\cset{\F^x_n:n\leq\omega})$, which contradicts $\F^x_\omega$ is c-irreducible. Similarly, $|Y_1|\geq 2$ implies $\F^x_\omega$ is c-reducible, which is impossible. Thus $|Y_0|<2$ and $|Y_1|<2$. The reader can verify that $\tup{|Y_0|,|Y_1|} = \tup{0,0}$, $\tup{0,1}$, $\tup{1,0}$, and $\tup{1,1}$ correspond to $\F^x_\omega \iso \C^\circ_\omega$, $\C^+_\omega$, $\C^-_\omega$, and $\C^\pm_\omega$, respectively.
\end{proof}

\begin{theorem}\label{thm:pretab-S4.3t}
    $\mathsf{PTAB}(\mathsf{S4.3}_t)=\cset{L^\pm,L^+,L^-,L^\circ,L^\ua}$ and $|\mathsf{PTAB}(\mathsf{S4.3}_t)|=5$.
\end{theorem}
\begin{proof}
    By Lemmas~\ref{lem:grz.3-pretabular}, \ref{lem:pretab-S43-fin-depth} and \ref{lem:4-pretab-S43-fin-depth}, we see that $\cset{L^\pm,L^+,L^-,L^\circ,L^\ua}\sub\mathsf{PTAB}(\mathsf{S4.3}_t)$ and $|\cset{L^\pm,L^+,L^-,L^\circ,L^\ua}|=5$. Take any $L\in\mathsf{PTAB}(\mathsf{S4.3}_t)$. Suppose $\mathrm{dep}(L)=\aleph_0$. By Theorem~\ref{thm:infchain-Lua} and Lemma~\ref{lem:grz.3-pretabular}, $L=L^\ua$. Suppose $\mathrm{dep}(L)=n<\aleph_0$. By Fact~\ref{fact:bounds}, $L\in\NExt(\mathsf{S4BP}^{n,1}_{1,1})$. By Theorem~\ref{thm:pretabular-S4BP-ch}, $L=\mathsf{Log}(\F^x_\omega)$ for some c-irreducible rooted $\F^x_\omega$. By Lemma~\ref{lem:c-irreducible-S4.3-finite-depth}, $\F^x_\omega\in\mathsf{IM}(\cset{\C^t_\omega:t\in\mathsf{tp}})$ and so $L\in\cset{L^t:t\in\mathsf{tp}}$.
\end{proof}

Note that $L^\ua=\mathsf{S4.3}_t\oplus\mathbf{grz}$, where $\mathbf{grz}=\B(\B(p\to\B p)\to p)\to p$. By Proposition~\ref{prop:fap-ptab-S4.3} and Theorems~\ref{thm:BP-KC-FMP} and \ref{thm:pretab-S4.3t}, we have

\begin{theorem}\label{thm:pretab-S4.3t-FMP}
    Let $L\in\mathsf{PTAB}(\mathsf{S4.3}_t)$. Then $L$ is finitely axiomatizable and has the FMP. As a corollary, $L$ is decidable.
\end{theorem}

By Theorems~\ref{thm:pretab-S4.3t} and \ref{thm:pretab-S4.3t-FMP}, we have

\begin{theorem}
    Tabularity in $\NExt(\mathsf{S4.3}_t)$ is decidable. That is, given any formula $\phi$, it is decidable whether $\mathsf{S4.3}_t\oplus\phi$ is tabular.
\end{theorem}

\section{Pretabular tense logics in $\NExt(\mathsf{S4BP}^{2,\omega}_{2,2})$}

In this section, we study pretabular logics over $\mathsf{S4BP}_{2,2}^{2,\omega}$. Comparing to $\mathsf{S4.3}_t$, the tense logic $\mathsf{S4BP}^{2,\omega}_{2,2}$ has weaker constraints on the width of logics, but a stronger constraint on the depth. As we will show in this section later, rooted frames of $\mathsf{S4BP}^{2,\omega}_{2,2}$ are garlands and hoops. We provide a full characterization of $\mathsf{PTAB}(\mathsf{S4BP}^{2,\omega}_{2,2})$ and it turns out that the cardinality of $\mathsf{PTAB}(\mathsf{S4BP}^{2,\omega}_{2,2})$ is $\aleph_0$. Thus we have a constructive proof for the claim in \cite{Rautenberg1979} that $\mathsf{PTAB}(\mathsf{S4}_t)$ is infinite.

Before characterizing the pretabular logics in $\NExt(\mathsf{S4BP}^{2,\omega}_{2,2})$, we need to define some finite skeletons of $\mathsf{S4BP}^{2,\omega}_{2,2}$, which plays an important role in our proof.

\begin{definition}
    Let $\G_\mathbb{Z}$ denote the frame $(\mathbb{Z},R_z)$ where 
    \begin{center}
        $R_z=\cset{\tup{i,i}:i\in\mathbb{Z}}\cup\cset{\tup{2i+1,2i}:i\in\mathbb{Z}}\cup\cset{\tup{2i+1,2i+2}:i\in\mathbb{Z}}$.
    \end{center}
    For all $\lambda\leq\omega$, we define $\G_\lambda$ as the frame $\G_\mathbb{Z}\rsto(1+\lambda)$. For each $n\in\mathbb{O}$, we define $\H_n$ as the frame $(1+n,(R_z\rsto{(1+n)})\cup\cset{\tup{n,0}})$. We call $\G_n$ an $n$-garland and $\H_n$ an $n$-hoop. Let $\mathcal{G}=\cset{\G_n:n\in\omega}$, $\breve{\mathcal{G}}=\cset{\breve{G_n}:n\in\omega}$ and $\mathcal{H}=\cset{\H_n:n\in\mathbb{O}}$.
\end{definition}
The frames $\G_\mathbb{Z}$, $\G_n$ and $\H_n$ are depicted in Figure.\ref{fig:pretabular-S4.3}. 

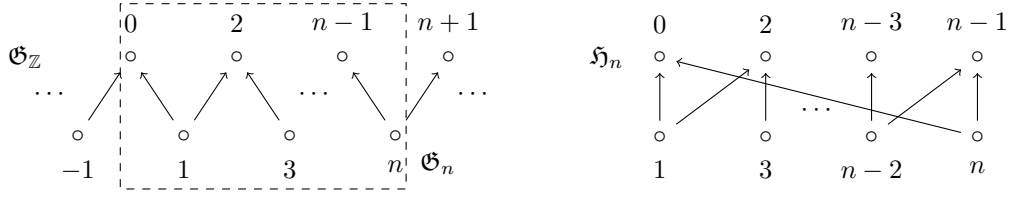
\begin{figure}[ht]
  \[
    \begin{tikzpicture}[scale=.7]
        \def\ptRad{.2pt}
        \draw (-2,1.5) node{$\G_\mathbb{Z}$};
        \draw (5.8,-.5) node{$\G_n$};
        \draw (3.5,.8) node{$\cdots$};
        \draw (-1.5,.8) node{$\cdots$};
        \draw (6.5,.8) node{$\cdots$};

        \draw [dashed] (-.2,2.5) rectangle (5.2,-1);

        \node (-1) at (-1,0)[label=below:$-1$]{$\circ$};
        \node (0) at (0,1.5)[label=above:$0$]{$\circ$};
        \node (1) at (1,0)[label=below:$1$]{$\circ$};
        \node (2) at (2,1.5)[label=above:$2$]{$\circ$};
        \node (3) at (3,0)[label=below:$3$]{$\circ$};
        \node (4) at (4,1.5)[label=above:$n-1$]{$\circ$};
        \node (5) at (5,0)[label=below:$n$]{$\circ$};
        \node (6) at (6,1.5)[label=above:$n+1$]{$\circ$};

        \draw [->,shorten <>=\ptRad] (-1) -- (0);
        \draw [->,shorten <>=\ptRad] (1) -- (0);
        \draw [->,shorten <>=\ptRad] (1) -- (2);
        \draw [->,shorten <>=\ptRad] (3) -- (2);
        \draw [->,shorten <>=\ptRad] (5) -- (4);
        \draw [->,shorten <>=\ptRad] (5) -- (6);
    \end{tikzpicture}
    \qquad\quad
    \begin{tikzpicture}[scale=.7]
        \def\ptRad{.2pt}
        \draw (-1,1.5) node{$\H_n$};
        \draw (3,.5) node{$\cdots$};

        \node (0) at (0,1.5)[label=above:$0$]{$\circ$};
        \node (1) at (0,0)[label=below:$1$]{$\circ$};
        \node (2) at (2,1.5)[label=above:$2$]{$\circ$};
        \node (3) at (2,0)[label=below:$3$]{$\circ$};
        \node (4) at (4,1.5)[label=above:$n-3$]{$\circ$};
        \node (5) at (4,0)[label=below:$n-2$]{$\circ$};
        \node (6) at (6,1.5)[label=above:$n-1$]{$\circ$};
        \node (7) at (6,0)[label=below:$n$]{$\circ$};

        \draw [->,shorten <>=\ptRad] (1) -- (0);
        \draw [->,shorten <>=\ptRad] (1) -- (2);
        \draw [->,shorten <>=\ptRad] (3) -- (2);
        \draw [->,shorten <>=\ptRad] (5) -- (4);
        \draw [->,shorten <>=\ptRad] (5) -- (6);
        \draw [->,shorten <>=\ptRad] (7) -- (6);
        \draw [->,shorten <>=\ptRad] (7) -- (0);
    \end{tikzpicture}
  \]
  \caption{The garlands $\G_n$ and hoops $\H_n$ for some $n\in\mathbb{O}$}
  \label{fig:BS222}
\end{figure}

\begin{lemma}
    $\mathcal{G}\cup\breve{\mathcal{G}}\cup\mathcal{H}\cup\cset{\G_\omega,\breve{\G_\omega}}\sub\mathsf{TM}(\G_\mathbb{Z})$.
\end{lemma}
\begin{proof}
    It is easy to see $\mathrm{abs}(\cdot):\G_\mathbb{Z}\twoheadrightarrow\G_\omega$, where $\mathrm{abs}(\cdot)$ is the absolute value function defined as follows:
    \begin{align*}
        \mathrm{abs}(x)=
        \begin{cases}
            x &\text{ if }x\geq 0,\\
            -x &\text{ otherwise.}
        \end{cases}
    \end{align*}
    Take any $n\in\omega$. Then the function $h_n:x\mapsto x\mod(1+n)$ is a t-morphism from $\G_\mathbb{Z}$ to $\H_n$. Moreover, $g_n:\G_\omega\twoheadrightarrow\G_n$, where
    \begin{align*}
        g(x)=
        \begin{cases}
            x\mod 2n &\text{ if }(x\mod 2n)\leq n;\\
            2n-(x\mod 2n) &\text{ otherwise.}
        \end{cases}
    \end{align*}
    Thus $\mathcal{G}\cup\mathcal{H}\cup\cset{\G_\omega}\sub\mathsf{TM}(\G_\mathbb{Z})$ and so $\breve{\mathcal{G}}\cup\cset{\breve{\G_\omega}}\sub\mathsf{TM}(\breve{\G_\mathbb{Z}})$. Note that $s:x\mapsto x+1$ is an isomorphism between $\G_\mathbb{Z}$ and $\breve{\G_\mathbb{Z}}$, we see $\mathcal{G}\cup\breve{\mathcal{G}}\cup\mathcal{H}\cup\cset{\G_\omega,\breve{\G_\omega}}\sub\mathsf{TM}(\G_\mathbb{Z})$.
\end{proof}

\begin{fact}\label{fact:finite-frame-BS222}
    Let $\F=(X,R)\in\mathsf{Fr}_r(\mathsf{S4}_t)$ be a skeleton. Then the following holds: 
    \begin{enumerate}[(1)]
        \item If $\F\md\mathbf{bd}_2$, then for all $x\in X$, $R[x]=\R[x]$ or $\breve{R}[x]=\R[x]$.
        \item If $\F\md\mathbf{bd}_2\wedge\mathbf{bw}^+_2$, then for all $x\in X$, $|R[x]|\leq 3$.
        \item If $\F\md\mathbf{bd}_2\wedge\mathbf{bw}^-_2$, then for all $x\in X$, $|\breve{R}[x]|\leq 3$.
    \end{enumerate}
\end{fact}
\begin{proof}
    The proof is standard, we omit it here.
\end{proof}

\begin{lemma}\label{lem:rooted-skeleton-BS222}
    Let $\F=(X,R)\in\mathsf{Fr}_r(\mathsf{S4BP}^{2,\omega}_{2,2})$ be a skeleton and $x\in X$. Then 
    \begin{enumerate}[(1)]
        \item for all $n\in\omega$, $\F\rsto\R^n[x]\in\mathsf{IM}(\mathcal{G}\cup\breve{\mathcal{G}}\cup\mathcal{H})$.
        \item $\F\rsto\R^\omega[x]\in\mathsf{IM}(\mathcal{G}\cup\breve{\mathcal{G}}\cup\mathcal{H}\cup\cset{\G_\omega,\breve{\G_\omega},\G_\mathbb{Z}})$.
    \end{enumerate}
\end{lemma}
\begin{proof}
    For (1), the proof proceeds by induction on $n$. When $n=0$, we see $\F\rsto\R^0[x]\iso\G_0$. Let $n>0$. By induction hypothesis, $\F\rsto\R^{n-1}[x]\in\mathsf{IM}(\mathcal{G}\cup\breve{\mathcal{G}}\cup\mathcal{H})$. Suppose $\F\rsto\R^{n-1}[x]\in\mathsf{IM}(\mathcal{H})$. Then $|\R[y]\cap\R^{n-1}[x]|=3$ for all $y\in X$. By Fact~\ref{fact:finite-frame-BS222}, $\R^n[x]=\R^{n-1}[x]$ and so $\F\rsto\R^{n}[x]\in\mathsf{IM}(\mathcal{H})$. Suppose $\F\rsto\R^{n-1}[x]\in\mathsf{IM}(\mathcal{G})$. Then there exists $f:\G_m\iso\F\rsto\R^{n-1}[x]$ for some $m\in\omega$. Since $\R^{n-1}[x]\sub\R^n[x]\neq X$, by Fact~\ref{fact:tabn}(1), $\R[f(0)]\nsubseteq\R^{n-1}[x]$ or $\R[f(m)]\nsubseteq\R^{n-1}[x]$. Then we have three cases:
    \begin{itemize}
        \item $\R[f(0)]\sub\R^{n-1}[x]$ and $\R[f(m)]\nsubseteq\R^{n-1}[x]$. Note that $\cset{f(m-1),f(m)}\sub\R[f(m)]$, by Fact~\ref{fact:finite-frame-BS222}, there exists a unique point $y\in\R[f(m)]\setminus\R^{n-1}[x]$. Clearly, $\R[y]\cap\R^{n-1}[x]=\cset{f(m)}$. Thus $f\cup\cset{\tup{m+1},y}:\G_{m+1}\iso\F\rsto\R^n[x]$.
        \item $\R[f(0)]\nsubseteq\R^{n-1}[x]$ and $\R[f(m)]\sub\R^{n-1}[x]$. By a similar argument, by Fact~\ref{fact:finite-frame-BS222}, there is a unique point $y\in\breve{R}[f(0)]\setminus\R^{n-1}[x]$ and we see that $g:\breve{\G_{m+1}}\iso\F\rsto\R^n[x]$, where $g(0)=y$ and $g(1+k)=f(k)$ for all $0<k\leq m$.
        \item $\R[f(0)]\nsubseteq\R^{n-1}[x]$ and $\R[f(m)]\nsubseteq\R^{n-1}[x]$. Again, by Fact~\ref{fact:finite-frame-BS222}, there exists a unique point $y\in\R[f(m)]\setminus\R^{n-1}[x]$ and a unique point $z\in\breve{R}[f(0)]\setminus\R^{n-1}[x]$. Suppose $y\neq z$. Then $g:\breve{\G_{m+2}}\iso\F\rsto\R^n[x]$, where $g(0)=z$, $g(m+1)=y$ and $g(1+k)=f(k)$ for all $0<k\leq m$. Otherwise, $y=z$ and we see that $f\cup\cset{\tup{m+1,y}}:\H_{m+1}\iso\F\rsto\R^n[x]$.
    \end{itemize}
    Suppose $\F\rsto\R^{n-1}[x]\in\mathsf{IM}(\breve{\mathcal{G}})$. Then $\breve{\F}\rsto\R^{n-1}[x]\in\mathsf{IM}(\mathcal{G})$, which entails $\breve{\F}\rsto\R^{n}[x]\in\mathsf{IM}(\mathcal{G})$. Thus $\F\rsto\R^{n}[x]\in\mathsf{IM}(\breve{\mathcal{G}})$. Since $\F$ is rooted, (2) follows from (1).
\end{proof}

Let $\mathsf{Ga}=\mathsf{Log}(\G_\mathbb{Z})$. The following theorem shows that $\mathsf{Ga}$ is the maximal logic in $\NExt(\mathsf{S4BP}^{2,\omega}_{2,2})$ with infinite z-degree.

\begin{theorem}\label{thm:infzigzag-Ga}
    Let $L\in\NExt(\mathsf{S4BP}^{2,\omega}_{2,2})$. Then $L\sub\mathsf{Ga}$ iff $\mathbf{bz}_n\not\in L$ for any $n\in\mathbb{Z}^+$.
\end{theorem}
\begin{proof}
    The left-to-right direction is trivial. Suppose $\mathbf{bz}_n\not\in L$ for any $n\in\mathbb{Z}^+$. Take any $\phi\not\in\mathsf{Ga}$ with $md(\phi)=k$. Then $\G_\mathbb{Z},m\not\md\phi$ for some $m\in\mathbb{Z}$. Since $\mathbf{bz}_{4k+2}\not\in L$ and $L$ is Kripke complete, there exists a frame $\F=(X,R)\in\mathsf{Fr}_r(L)$ and $x\in X$ such that $\F,x\not\md\mathbf{bz}_{4k+2}$. Then $\F^S,C(x)\not\md\mathbf{bz}_{4k+2}$ and so $\mathsf{zdg}(C(x))\geq 4k+2$. Let $\G=\F^S\rsto(R^S)_\sharp^{4k+1}[C(x)]$. By Lemma~\ref{lem:rooted-skeleton-BS222}, $\G\in\mathsf{IM}(\mathcal{G}\cup\breve{\mathcal{G}}\cup\mathcal{H})$. Since $\mathsf{zdg}(C(x))\geq 4k+2$, $\G\not\in\mathsf{IM}(\mathcal{H})$. Suppose $g:\G\iso\G_l$ for some $l\in\omega$. Since $\mathsf{zdg}(C(x))\geq 4k+2$, we see $l\geq 4k+1$. Let $f:1+l\to\mathbb{Z}$ be the function such that $f(x)=x+m+(m\mod 2)-2k$. Then we see $f\circ g:(\F,y)\to^{l-1}(\G_\mathbb{Z},m)$, where $y$ is such that $g(y)=2k-(m\mod 2)$. By Lemma~\ref{lem:k-t-morphism}, $\F,x\not\md\phi$ and so $\phi\not\in L$.
\end{proof}

\begin{corollary}\label{coro:extensions-Ga}
    Let $\mathcal{F}\sub\mathcal{G}\cup\breve{\mathcal{G}}\cup\mathcal{H}$ be infinite. Then $\mathsf{Ga}=\mathsf{Log}(\mathcal{F})$.
\end{corollary}

\begin{theorem}\label{thm:Ga-is-pretabular}
    $\mathsf{Ga}\in\mathsf{PTAB}(\mathsf{S4BP}^{2,\omega}_{2,2})$.
\end{theorem}
\begin{proof}
    Clearly, $\mathsf{Ga}$ is non-tabular. Take any $L\supsetneq\mathsf{Ga}$. By Theorem~\ref{thm:infzigzag-Ga}, $\mathbf{bz}_n\in L$ for some $n\in\mathbb{Z}^+$. By Lemma~\ref{lem:rooted-skeleton-BS222}, $\mathsf{Fr}_r(L)\sub\mathsf{IM}(\mathcal{G}\cup\breve{\mathcal{G}}\cup\mathcal{H})$. By Corollary~\ref{coro:extensions-Ga}, $\mathsf{Fr}_r(L)\sub\mathsf{IM}(\mathcal{F})$ for some finite $\mathcal{F}\sub\mathcal{G}\cup\breve{\mathcal{G}}\cup\mathcal{H}$. By Theorem~\ref{thm:BP-KC-FMP}, we see that $L=\mathsf{Log}(\bigoplus\mathcal{F})$ is tabular.
\end{proof}

In what follows, we give a characterization for pretabular logics over $\mathsf{S4BP}^{2,\omega}_{2,2}$ with finite z-degree.

\begin{lemma}\label{lem:finite-BS222-frame-irreducible}
    Let $m\leq n\in\omega$, $\F^x_1$ be a pre-skeleton and $f:(\G_n)^m_1\twoheadrightarrow\F^x_1$. Then
    \begin{enumerate}[(1)]
        \item $f^{-1}[f(m)]=\cset{m}$ and $f^{-1}[f(m_1)]=\cset{m_1}$.
        \item for all different $k,l\in 1+n$, $f(k)=f(l)$ implies $k+l=2m$.
        \item if $2m<n$, then $f:(\G_n)^m_1\iso\F^x_1$.
    \end{enumerate}
\end{lemma}
\begin{proof}
    Recall that the domain of $(\G_n)^m_1$ is $(1+n)\cup\cset{m_1}$. Since $f:(\G_n)^m_1\twoheadrightarrow\F^x_1$, by Lemma~\ref{lem:finiteK4-ontopmorphism-cluster}, $f[\cset{m,m_1}]=\cset{x,x_1}$ is the unique proper cluster in $\F^x_1$. For (1), take any $k\in(\G_n)^m_1$. Suppose $f(k)=f(m)$. Then $C(f(k))$ is a proper cluster. By Lemma~\ref{lem:finiteK4-ontopmorphism-cluster}, $C(k)$ is proper. Since $f(m)\neq f(m_1)$ and $\cset{m,m_1}$ is the unique proper cluster in $(\G_n)^m_1$, we see $k=m$. Hence $f^{-1}[f(m)]=\cset{m}$. Similarly, $f^{-1}[f(m_1)]=\cset{m_1}$.

    For (2), take any different $k,l\in 1+n$ such that $f(k)=f(l)$. By (1), $k\neq m\neq l$. Let $\F^x_1=(X,R)$. Then there exists $a\in\mathbb{Z}^+$ such that $f(m)\in\R^{a}[f(k)]\setminus\R^{a-1}[f(k)]$. Thus $f(m)\in f[\R^{a}[k]]\setminus f[\R^{a-1}[k]]$ and $f(m)\in f[\R^{a}[l]]\setminus f[\R^{a-1}[l]]$. By (1), $m\in\R^{a}[k]\setminus\R^{a-1}[k]$ and $m\in\R^{a}[l]\setminus\R^{a-1}[l]$. Since $k\neq l$, we see $\cset{k,l}=\cset{m-a,m+a}$. Thus $k+l=2m$.

    For (3), suppose $2m<n$. By (1) and (2), $f^{-1}[f[y]]=\cset{y}$ for all $y\in\cset{m,m_1}\cup\cset{k\in 1+n:k>2m}$. Since $2m+1\in\R^m[m+a]\setminus\R^m[m-a]$ for all $a\leq m$, by (2), we see $f^{-1}[f[k]]=\cset{k}$ for all $k\leq 2m$. Thus $f$ is injective, which entails $f:(\G_n)^m_1\iso\F^x_1$.
\end{proof}

\begin{lemma}\label{lem:Gn-irreducible}
    Let $m,n\in\omega$. Suppose $n>0$ and $2m\leq n$. Then 
    \begin{enumerate}[(1)]
        \item $(\G_n)^m_1$ is c-irreducible if and only if $2m\neq n$.
        \item $(\breve{\G_n})^m_1$ is c-irreducible if and only if $2m\neq n$.
    \end{enumerate}
\end{lemma}
\begin{proof}
    Note that $\mathsf{TM}((\breve{\G_n})^m_1)=\cset{\breve{\F}:\F\in\mathsf{TM}((\G_n)^m_1)}$, (2) follows from (1) immediately. For (1), suppose $2m\neq n$. By Lemma~\ref{lem:finite-BS222-frame-irreducible}(3), $(\G_n)^m_1$ is c-irreducible. For the other direction, suppose $2m=n$. Let $f:X\to(1+m)\cup\cset{m_1}$ be defined by:
    \begin{align*}
        f(x)=
        \begin{cases}
            x &\text{ if }x\leq m \text{ or }x=m_1;\\
            n-x &\text{ otherwise.}
        \end{cases}
    \end{align*}
    Thus $f:(\G_n)^m_1\twoheadrightarrow(\G_m)^m_1$. Since $n>0$, we see $n>m$ and so $(\G_m)^m_1\not\iso(\G_n)^m_1$. By Lemma~\ref{lem:finiteK4-ontopmorphism-cluster}, $(\G_m)^m_1\not\in\mathsf{TM}(\G_n)$. Thus $(\G_n)^m_1$ is c-reducible.
\end{proof}

\begin{lemma}\label{lem:Hn-reducible}
    For all $n\in\mathbb{O}$ and $m\leq n$, $(\H_n)^m_1$ is c-reducible.
\end{lemma}
\begin{proof}
    Suppose $n=4k+1$ for some $k\in\omega$ and $m\in\mathbb{E}$. Let $f:(4k+2)\cup\cset{m_1}\to(4k+2)\cup\cset{(2k)_1}$ be the map defined as:
    \begin{align*}
        f(x)=
        \begin{cases}
            (x+2k-m)\mod(4k+1) &\text{ if } x\in(4k+2);\\
            (2k)_1 &\text{ otherwise.}
        \end{cases}
    \end{align*}
    Then $f:(\H_n)^m_1\iso(\H_n)^{2k}_1$. We define $g:(4k+2)\cup\cset{(2k)_1}\to(2k+1)\cup\cset{0_1}$ by
    \begin{align*}
        g(x)=
        \begin{cases}
            l &\text{ if } \mathrm{abs}(x-2m)=l \text{ and }x\in(n+1);\\
            0_1 &\text{ otherwise.}
        \end{cases}
    \end{align*}
    Then we see that $g:(\H_n)^{2k}_1\twoheadrightarrow(\G_{2k+2})^0_1$, which entails $g\circ f:(\H_n)^{m}_1\twoheadrightarrow(\G_{2k+2})^0_1$. Note that $(\G_{2k+2})^0_1\not\in\mathsf{TM}(\H_n)$ and $(\G_{2k+2})^0_1\not\iso(\H_n)^{m}_1$, $(\H_n)^{m}_1$ is c-reducible. Assume $m\in\mathbb{O}$. Then we can construct maps $f':(\H_n)^m_1\iso(\H_n)^{2k+1}_1$ and $g':(\H_n)^{2k+1}_1\twoheadrightarrow(\breve{\G_{2k+1}})^0_1$ in a similar way, which also implies that $(\H_n)^{m}_1$ is c-reducible. Similar arguments work for the case when $n=4k+3$ for some $k\in\omega$.
\end{proof}

Now we are ready to prove the main theorem of this section:

\begin{theorem}\label{thm:pretabular-BS222-ch}
  Let $\mathsf{S5}_t=\mathsf{S4}_t\oplus(\D p\to\B\D p)$. Then
    \begin{center}
      $\mathsf{PTAB}(\mathsf{S4BP}^{2,\omega}_{2,2})=\cset{\mathsf{Ga},\mathsf{S5}_t}\cup\cset{\mathsf{Log}((\G_n)^m_\omega),\mathsf{Log}((\breve{\G_n})^m_\omega):2m<n\in\mathbb{Z}^+}$.
    \end{center}
\end{theorem}
\begin{proof}
    Note that $\mathsf{S5}_t=L^\circ$ and $\mathsf{S5}_t\supseteq\mathsf{S4BP}^{2,\omega}_{2,2}$, by Theorem~\ref{thm:pretab-S4.3t}, $\mathsf{S5}_t$ is pretabular. By Theorem~\ref{thm:Ga-is-pretabular}, $\mathsf{Ga}\in\mathsf{PTAB}(\mathsf{S4BP}^{2,\omega}_{2,2})$. Take any $2m<n\in\omega$. By Lemma~\ref{lem:Gn-irreducible} and Lemma~\ref{lem:c-irreducible-2-omega}, $(\G_n)^m_\omega$ is c-irreducible. By Theorem~\ref{thm:pretabular-S4BP-ch}, $\mathsf{Log}((\G_n)^m_\omega)$ is pretabular.

    Take any $L\in\mathsf{PTAB}(\mathsf{S4BP}^{2,\omega}_{2,2})$. Suppose $\mathsf{bz}_n\not\in L$ for any $n\in\mathbb{Z}^+$. Then by Theorem~\ref{thm:infzigzag-Ga}, $L\sub\mathsf{Ga}$. Since $L$ is pretabular, $L=\mathsf{Ga}$. Suppose $\mathsf{bz}_n\in L$ for some $n\in\mathbb{Z}^+$. Then $L\in\mathsf{PTAB}(\mathsf{S4BP}^{2,n}_{2,2})$. By Theorem~\ref{thm:pretabular-S4BP-ch}, $L=\mathsf{Log}(\F^x_\omega)$ for some c-irreducible rooted finite pre-skeleton. By Lemma~\ref{lem:rooted-skeleton-BS222}, $\F\in\mathsf{IM}(\mathcal{G}\cup\breve{\mathcal{G}}\cup\mathcal{H})$. If $|\F|=1$, then $\F^x_\omega\iso(\G_0)^0_\omega$ and so $L=\mathsf{S5}_t$. Suppose $|\F|\neq 1$, then $|\F|=1+n$ for some $n\in\mathbb{Z}^+$. By Lemma~\ref{lem:Gn-irreducible} and Lemma~\ref{lem:Hn-reducible}, $\F^x_\omega\iso(\G_n)^m_\omega$ or $\F^x_\omega\iso(\breve{\G_n})^m_\omega$ for some $m\in\omega$ such that $2m\neq n$. If $2m<n$, then we are done. If $2m>n$, then we see that $\F^x_\omega\iso(\G_n)^{n-m}_\omega$ or $\F^x_\omega\iso(\breve{\G_n})^{n-m}_\omega$. Thus $L\in\cset{\mathsf{Log}((\G_n)^m_\omega):2m<n\in\omega}$.
\end{proof}

\begin{example}
    Note that the frame $\F^{x_4}_\omega$ given in Figure~\ref{fig:ske-preske} is isomorphic to $(\breve{G_3})^0_\omega$, by Theorem~\ref{thm:pretabular-BS222-ch}, we see that $\mathsf{Log}(\F^{x_4}_\omega)$ is pretabular.
\end{example}

\begin{theorem}\label{thm:pretab-BS222}
    $|\mathsf{PTAB}(\mathsf{S4BP}^{2,\omega}_{2,2})|=\aleph_0$.
\end{theorem}
\begin{proof}
    By Theorem~\ref{thm:pretabular-BS222-ch}, $|\mathsf{PTAB}(\mathsf{S4BP}^{2,\omega}_{2,2})|\leq\aleph_0$. Take any $n,m,k,l\in\omega$ such that $2m<n$ and $2l<k$. It is clear that $(\G_n)^m_\omega$, $(\breve{\G_n})^m_\omega$, $(\G_k)^l_\omega$ and $(\breve{\G_k})^l_\omega$ are pairwise non-isomorphic. By Lemma~\ref{lem:irreducible-frame-same-logic-implies-iso}, their logics are pairwise different. Thus by Theorem~\ref{thm:pretabular-BS222-ch}, $|\mathsf{PTAB}(\mathsf{S4BP}^{2,\omega}_{2,2})|=\aleph_0$.
\end{proof}

Moreover, we show the following anti-dichotomy theorem for cardinality of pretabular extensions in $\NExt(\mathsf{S4BP}_{2,2}^{2,\omega})$:

\begin{theorem}\label{thm:anti-dichotomy-BS222}
    For all $\kappa\leq{\aleph_0}$, there exists $L\in\NExt(\mathsf{S4}_t)$ with $|\mathsf{PTAB}(L)|=\kappa$.
\end{theorem}
\begin{proof}
    By Theorem~\ref{thm:pretab-BS222}, $|\mathsf{PTAB}(\mathsf{S4BP}^{2,\omega}_{2,2})|=\aleph_0$. Obviously, $\mathsf{Log}(\C_1))\supseteq\mathsf{S4BP}_{2,2}^{2,\omega}$ and $|\mathsf{PTAB}(\mathsf{Log}(\C_1))|=0$. Take any $\kappa\in\mathbb{Z}^+$. Let $\mathcal{F}=\cset{(\G_{n+1})^0_\omega:n<\kappa}$ and $L=\bigcap_{\F\in\mathcal{F}}\mathsf{Log}(\F)$. Then $|\mathcal{F}|=\kappa$. Note that $\mathbf{bd}_1\not\in L$ and $\mathbf{bz}_{\kappa+1}\in L$, $\cset{\mathsf{Ga},\mathsf{S5}_t}\cap\mathsf{PTAB}(L)=\ve$. Take any $\F^x_\omega\in\cset{(\G_n)^m_\omega,(\breve{\G_n})^m_\omega:2m<n\in\mathbb{Z}^+}$. Suppose $\F^x_\omega\not\in\mathcal{F}$. Then similar to the proof of Lemma~\ref{lem:irreducible-frame-same-logic-implies-iso}, we see $\J^{\mathsf{zdg}(\F)+\kappa}(\F^x_1,x)\not\in L$. Thus $\F\in\mathcal{F}$ if and only if $\mathsf{Log}(\F)\supseteq L$. By Theorem~\ref{thm:pretabular-BS222-ch}, $\mathsf{PTAB}(L)=\cset{\mathsf{Log}(\F):\F\in\mathcal{F}}$. By Lemma~\ref{lem:irreducible-frame-same-logic-implies-iso}, $|\mathsf{PTAB}(L)|=\kappa$. 
\end{proof}

\begin{theorem}\label{thm:BS222-pretab-fmp}
    For all $L\in\mathsf{PTAB}(\mathsf{S4BP}^{2,\omega}_{2,2})$, $L$ has the FMP.
\end{theorem}
\begin{proof}
    Follows from Lemma~\ref{lem:ThFomega=CapThFn}, Corollary~\ref{coro:extensions-Ga} and Theorem~\ref{thm:pretabular-BS222-ch} immediately.
\end{proof}

\begin{remark}\label{rem:Kracht-Ga}
    The results obtained in this section is closely related to the ones in \cite[Section 4]{Kracht1992}. The logic $\mathsf{Ga}$ was defined to be $\mathsf{S4}_t\oplus\cset{\mathbf{grz},\mathbf{alt}^+_3,\mathbf{alt}^-_3,\mathbf{bd}_2}$. Garlands $\G_n$ were also defined there. It was proved that $\mathsf{Ga}=\mathsf{Log}(\G_\omega)=\bigcap_{n\in\omega}\mathsf{Log}(\G_n)$ is pretabular. 
    
    However, there are some problematic claims in \cite[Section 4]{Kracht1992}, which makes the characterization of $\NExt(\mathsf{Ga})$ given there incomplete. It was claimed that a rooted frame $\F$ validates $\mathsf{Ga}$ if and only if $\F\iso\G_n$ for some $n\in\omega$. It follows that $\cset{\mathsf{Log}(\G_n):n\text{ is prime}}$ is the set of all logics in $\NExt(\mathsf{Ga})$ which are of co-dimension 3. But as Lemma~\ref{lem:rooted-skeleton-BS222} shows, the class of rooted frame for $\mathsf{Ga}$ is $\mathsf{IM}(\mathcal{G}\cup\breve{\mathcal{G}}\cup\mathcal{H}\cup\cset{\G_\omega,\breve{\G_\omega},\G_\mathbb{Z}})$, but not $\mathsf{IM}(\mathcal{G})$. Consider the frames $\H_3$ and $\breve{G_2}$, as shown in Figure~\ref{fig:H3G2}. Then we see that $\H_3\md\mathsf{Ga}$ but $\H_3\ncong\G_n$ for any $n\in\omega$. Moreover, note that $\mathbf{tab}^T_4\wedge\neg\J^4(\G_3,0)\in\mathsf{Log}(\H_3)$ and $\mathbf{tab}^T_3\not\in\mathsf{Log}(\H_3)$, we see that $\mathsf{Log}(\H_3)\not\in\cset{\bigcap_{i\in I}\mathsf{Log}(\G_i):I\sub\omega}$. Thus $\mathsf{Log}(\H_3)$ is missing in the characterization given by Kracht \cite{Kracht1992}. It is also straightforward to show that $\mathsf{Log}(\breve{\G_2})\not\in\cset{\bigcap_{i\in I}\mathsf{Log}(\G_i):I\sub\omega}$.
    
    With the results obtained in this section, we can even give a full characterization of $\NExt(\mathsf{Ga})$. It can be shown that $\NExt(\mathsf{Ga})$ is dually isomorphic to the distributive lattice freely $\bigcup$-generated by $(\mathcal{G}\cup\breve{\mathcal{G}}\cup\mathcal{H},\twoheadrightarrow)$. But we do not go into the details now and leave this work for future research.
\begin{figure}[ht]
    \[
      \begin{tikzpicture}[scale=.7]
          \def\ptRad{.2pt}
          \draw (.5,1.5) node{$\breve{\G_2}$};
          \node (1) at (1,0)[label=below:$0$]{$\circ$};
          \node (2) at (2,1.5)[label=above:$1$]{$\circ$};
          \node (3) at (3,0)[label=below:$2$]{$\circ$};
  
          \draw [->,shorten <>=\ptRad] (1) -- (2);
          \draw [->,shorten <>=\ptRad] (3) -- (2);
      \end{tikzpicture}
      \qquad\qquad
      \begin{tikzpicture}[scale=.7]
          \def\ptRad{.2pt}
          \draw (-1,1.5) node{$\H_3$};
  
          \node (0) at (0,1.5)[label=above:$0$]{$\circ$};
          \node (1) at (0,0)[label=below:$1$]{$\circ$};
          \node (2) at (2,1.5)[label=above:$2$]{$\circ$};
          \node (3) at (2,0)[label=below:$3$]{$\circ$};
  
          \draw [->,shorten <>=\ptRad] (1) -- (0);
          \draw [->,shorten <>=\ptRad] (1) -- (2);
          \draw [->,shorten <>=\ptRad] (3) -- (2);
          \draw [->,shorten <>=\ptRad] (3) -- (0);
      \end{tikzpicture}
    \]
    \caption{The garland $\breve{G_2}$ and hoop $\H_3$}
    \label{fig:H3G2}
  \end{figure}
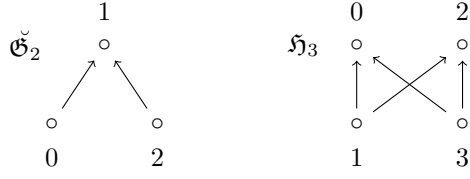
\end{remark}

\section{Pretabular tense logics in $\NExt(\mathsf{S4BP}^{2,\omega}_{2,3})$}\label{sec:BS223}

So far, we studied several families of tense logics with bounded parameters and constructed countably many pretabular logics in $\NExt(\mathsf{S4}_t)$. In this section, we study the pretabular tense logics extending $\mathsf{S4BP}^{2,\omega}_{2,3}$, which has a weaker constraint on back-width than $\mathsf{S4BP}^{2,\omega}_{2,2}$. As we will see soon, the class of rooted frames for $\mathsf{S4BP}^{2,\omega}_{2,3}$ is much more complicated than the one for $\mathsf{S4BP}^{2,\omega}_{2,2}$. In fact, there are continuum many frames for $\mathsf{S4BP}^{2,\omega}_{2,3}$ whose logics are pairwise different.

The aim of this section is to show that $|\mathsf{PTAB}(\mathsf{S4BP}^{2,\omega}_{2,3})|=2^{\aleph_0}$. The strategy is as follows: We first construct a continual family of sequences by generalizing the Thue-Morse sequences. Based on these sequences, we construct corresponding `umbrella-like' frames and show that their logics are pairwise different and pretabular. As a corollary, $|\mathsf{PTAB}(\mathsf{S4}_t)|=2^{\aleph_0}$, which answers the open problem given in \cite{Rautenberg1979}. 

\subsection{Preliminaries of sequences}

For all $i,j\in\mathbb{Z}$, we write $[i,j]$ for $\cset{k\in\mathbb{Z}:i\leq k\leq j}$. A subset $I$ of $\mathbb{Z}$ is said to be an {\em interval} in $\mathbb{Z}$ if for all $i,j\in I$, $[i,j]\sub I$. A map $t:\mathbb{Z}\to\mathbb{Z}$ is called a {\em translation} if there exists $k\in\mathbb{Z}$ such that $t(i)=i+k$ for all $i\in\mathbb{Z}$. We write $t:a\mapsto b$ for the translation $t$ such that $t:i\mapsto(i+b-a)$ for all $i\in\mathbb{Z}$. Let $X$ be a non-empty set. An {\em $X$-sequence} is a partial function $f:\mathbb{Z}\to X$ where $\mathsf{dom}(f)$ is an interval in $\mathbb{Z}$. Let $\mathsf{Seq}(X)$ and $\mathsf{Seq}^{<\aleph_0}(X)$ denote the sets of all $X$-sequences and finite $X$-sequences, respectively.

Let $\alpha:[a,b]\to X$ be a nonempty finite sequence such that $\alpha(i)=x_{i-a}$ for all $i\in[a,b]$. Then we write $\tup{a,\tup{x_0,\cdots,x_b}}$ for $\alpha$. We write $\tup{x_0,\cdots,x_b}$ for $\alpha$ if $a=0$.

\begin{definition}
    Let $X\neq\ve$ and $\alpha,\beta\in\mathsf{Seq}(X)$. Then we say (i) {\em $\alpha$ is embedded into $\beta$ (notation: $\alpha\trianglelefteq\beta$)}, if $\alpha\circ s\sub\beta$ for some translation $s$. (ii) {\em $\alpha$ is finitely covered by $\beta$ (notation: $\alpha\preceq\beta$)}, if $\gamma\trianglelefteq\beta$ for all finite sequence $\gamma\sub\alpha$. (iii) {\em $\alpha$ and $\beta$ are similar (notation: $\alpha\approx\beta$)}, if $\alpha\preceq\beta$ and $\beta\preceq\alpha$. $\alpha$ and $\beta$ are {\em dissimilar} if $\alpha\not\approx\beta$.
\end{definition}

\begin{definition}[Concatenation]\label{def:concatenation}
    Let $\alpha:[a,b]\to X$ and $\beta:[c,d]\to X$ be finite $X$-sequences for some nonempty set $X$. Then we define the sequence $\alpha\ast\beta$ by
    \begin{center}
        $\alpha\ast\beta={\tup{\alpha(a),\cdots,\alpha(b),\beta(c),\cdots,\beta(d)}}$. 
    \end{center}
    The sequences $\alpha^\dagger\ast\beta$ and $\alpha\ast\beta^\dagger$ are defined as follows:
    \begin{center}
        $\alpha^\dagger\ast\beta=\tup{a,\alpha\ast\beta}$ and $\alpha\ast\beta^\dagger=\tup{c+a-(b+1),\alpha\ast\beta}$.
    \end{center}
    Let $\tup{\alpha_i:i\in n}$ be a finite tuple of finite $X$-sequences. Then we write $\alpha_0\ast\alpha_1\ast\cdots\ast\alpha_{n-1}$ or $\alpha_0\alpha_1\cdots\alpha_{n-1}$ for $(\cdots(\alpha_0\ast\alpha_1)\ast\cdots\ast\alpha_{n-2})\ast\alpha_{n-1}$. Moreover, we define 
    \begin{center}
        $\alpha_0\ast\cdots\ast{\alpha_m}^\dagger\ast\cdots\ast\alpha_{n-1}=((\alpha_0\ast\cdots\ast{\alpha_{m-1}})\ast{\alpha_m}^\dagger)^\dagger\ast({\alpha_{m+1}}\ast\cdots\ast\alpha_{n-1})$.
    \end{center}
\end{definition}

The notation in Definition~\ref{def:concatenation} looks complicated, but the idea behind is simple. Given any finite tuple $A=\tup{\alpha_i:i\in n}$ of finite $X$-sequences, the sequence $\alpha_0\ast\cdots\ast{\alpha_m}^\dagger\ast\cdots\ast\alpha_{n-1}$ is designed to be the concatenation of $A$ which always preserves the index of $\alpha_m$. An example is given in Table~\ref{tab:example-concatenation}.

\begin{table}[htbp]\label{tab:example-concatenation}
    \begin{center}
        \begin{tabular}{l|cccccccccccccccccc}
            $\mathbb{Z}$&$\cdots$&-5&-4&-3&-2&-1&0&1&2&3&4&5&6&$\cdots$\\
            \\
            $\alpha$&&&&&&a&b&c&d&&&&\\
            $\beta$&&&&&x&y&z&&&&&&\\
            $\gamma$&&&&&u&v&w&&&&&&\\
            $\alpha\ast\beta$&&&&&&&a&b&c&d&x&y&z\\
            $\alpha^\dagger\ast\beta$&&&&&&a&b&c&d&x&y&z&\\
            $\alpha\ast\beta\ast\gamma^\dagger\ast\alpha$&$\cdots$&x&y&z&u&v&w&a&b&c&d&&\\
        \end{tabular}
    \end{center}
    \caption{Example of concatenations of sequences}
\end{table}

\begin{definition}
    Let $X$ be a non-empty set, $A\sub\mathsf{Seq}(X)$ and $n\in\omega$. Then the set $\mathsf{Con}(A,n)$ of all $n$-concatenations of $A$ is defined as follows:
    \begin{center}
        $\mathsf{Con}(A,n)=\cset{\alpha_0\cdots\alpha_{n-1}:\forall{i<n}(\alpha_i\in A)}$.
    \end{center}
    Let $\mathsf{Con}(A)=\bigcup_{i\in\omega}\mathsf{Con}(A,n)$. We call $\mathsf{Con}(A)$ the set of all concatenations of $A$.
\end{definition}

\begin{definition}
    Let $\alpha:\mathbb{Z}\to 2$ be a map. We say that $\alpha$ is {\em finitely perfect} if for all finite subsequence $\beta$ of $\alpha$, there is $n\in\omega$ such that $\beta\triangleleft\zeta$ for all $\zeta\triangleleft\alpha$ with $|\zeta|>n$.
\end{definition}

\subsection{Generalized Thue-Morse sequences}

In this subsection, our goal is to construct a continuum of pairwise dissimilar finitely perfect $2$-sequences. For every function $f:\omega\to 2$, we define the generalized Thue-Morse sequence $\chi^f$ generated by $f$. 

The Thus-Morse sequence $\alpha^t$ is defined by $\alpha^t=\bigcup_{i\in\omega}\alpha^t_i$, where $\alpha^t_0=\tup{0}$ and $\alpha^t_{i}=\alpha_i\ast\ol{\alpha_{i-1}}$ for all $i>0$. The sequence $\alpha^t$ has many nice properties. For example, $\alpha^t$ is shown to be overlap-free, i.e., $x\beta x\beta x\ntrianglelefteq\alpha^t$ for any $2$-sequence $\beta$ and $x<2$, see \cite[Proposition 5.1.6]{Fogg.Berthe.ea2002}. The readers can also check that $\alpha^t$ is finitely perfect. To obtain a continual family of finitely perfect $2$-sequences, we generalize the Thue-Morse sequence as follows:

\begin{definition}[Generalized Thus-Morse sequence]
    Let $f:\omega\to 2$. For each $i\in\omega$, the finite binary sequence $\chi_i$ is defined as follows:
    \begin{itemize}
        \item $\chi^f_0=\tup{0,0,1}$;
        \item $\chi^f_{2i+1}={\chi_{2i}}^\dagger\ast\tup{f(2i)}\ast\ol{\chi_{2i}}$;
        \item $\chi^f_{2i+2}=\ol{\chi_{2i+1}}\ast\tup{f(2i+1)}\ast{\chi_{2i+1}}^\dagger$.
    \end{itemize}
    Let $\chi^f=\bigcup_{i\in\omega}\chi^f_i$. Then we see $\chi^f$ is a function from $\mathbb{Z}$ to $2$. 
    The sequence $\chi^f$ is called the generalized Thus-Morse sequence generated by $f$.
\end{definition}

\begin{example}
    Consider the maps $f:\omega\to\cset{0}$ and $g:\omega\to\cset{1}$. Then the sequences $\chi^f$ and $\chi^g$ are constructed as follows:
    \begin{center}
        \begin{tabular}{l|cccccccccccccccccc}
            $\mathbb{Z}$&$\cdots$&-8&-7&-6&-5&-4&-3&-2&-1&0&1&2&3&4&5&6&7&$\cdots$\\
            \\
            $\chi^f_0$&&&&&&&&&&0&0&1&&&&&\\
            $\chi^g_1$&&&&&&&&&&0&0&1&1&1&1&0&\\
            $\chi^f_2$&&1&1&0&1&0&0&1&0&0&0&1&0&1&1&0&&\\
            $\chi^g_2$&&1&1&0&0&0&0&1&1&0&0&1&1&1&1&0&&\\
            &&&&&&&&&$\vdots$&&&&&&&\\
            $\chi^f$&$\cdots$&1&1&0&1&0&0&1&0&0&0&1&0&1&1&0&0&$\cdots$\\
            $\chi^g$&$\cdots$&1&1&0&0&0&0&1&1&0&0&1&1&1&1&0&1&$\cdots$\\
        \end{tabular}
    \end{center}
\end{example}

\begin{lemma}
    Let $f\in 2^\omega$, $j\leq i\in\omega$ and $\nu(j)=2^{i-j}$. Then there are $\alpha_1,\cdots,\alpha_{\nu(j)}\in\cset{\chi^f_j,\ol{\chi^f_j}}$ and $n_1,\cdots,n_{{\nu(j)}-1}\in 2$ such that $\chi^f_i\approx\alpha_0n_0\cdots\alpha_{{\nu(j)}-1}n_{{\nu(j)}-1}\alpha_{{\nu(j)}}$.
\end{lemma}
\begin{proof}
    The proof proceeds by induction on $i-j$. The case $i-j=0$ is trivial. Suppose $i-j>0$. By induction hypothesis, $\chi^f_i\approx\alpha_0n_0\cdots\alpha_{\nu(j+1)}n_{\nu(j+1)-1}\alpha_{{\nu(j+1)}}$ for some $\alpha_1,\cdots,\alpha_{\nu(j)}\in\cset{\chi^f_{j+1},\ol{\chi^f_{j+1}}}$ and $n_0,\cdots,n_{{\nu(j+1)}-1}\in 2$. Note that for all $l\leq\nu(j+1)$, we have $\alpha_l\approx\beta_l^1m_l\beta_l^2$ for some $\beta_l^1,\beta_l^2\in\cset{\chi^f_j,\ol{\chi^f_j}}$ and $m_l\in 2$. Thus we have $\chi^f_i\approx\beta_0^1m_0\beta_0^2n_0\cdots n_{\nu(j+1)-1}\beta_{\nu(j+1)}^1m_{\nu(j+1)}\beta_{\nu(j+1)}^2$, which concludes the proof.
\end{proof}

\begin{corollary}\label{coro:construction-chif}
    Let $f\in 2^\omega$. Then for all $j\in\omega$,
    \[
    \chi^f\preceq\mathsf{Con}(\cset{\chi^f_j0,\chi^f_j1,\ol{\chi^f_j}0,\ol{\chi^f_j}1}).
    \]
\end{corollary}

Intuitively, Corollary~\ref{coro:construction-chif} says that the $2$-sequence is builded up by iterating elements in $\cset{\chi^f_j0,\chi^f_j1,\ol{\chi^f_j}0,\ol{\chi^f_j}1}$. This makes the structure of $\chi^f$ clear.

\begin{lemma}\label{lem:chif-finitely-perfect}
    Let $f\in 2^\omega$. Then $\chi^f$ is finitely perfect.
\end{lemma}
\begin{proof}
    Take any finite subsequence $\alpha\triangleleft\chi^f$. Then $\alpha\triangleleft\chi^f_{i-1}$ for some $i\in\omega$. Let $n=2(|\chi^f_i|+1)$. Take any $\gamma\triangleleft\chi^f$ with $|\gamma|>n$. By Corollary~\ref{coro:construction-chif}, $\beta\triangleleft\gamma$ for some $\beta\in\cset{\chi^f_i0,\chi^f_i1,\ol{\chi^f_i}0,\ol{\chi^f_i}1}$. Note that $\chi^f_{i-1}\triangleleft\chi^f_i$ and $\chi^f_{i-1}\triangleleft\ol{\chi^f_i}$, we see $\alpha\triangleleft\beta\triangleleft\gamma$.
\end{proof}

\begin{lemma}\label{lem:chifi-emb-chifi1-unique}
    Let $f\in 2^\omega$ and $i\in\mathbb{Z}^+$. Let $\alpha:[a,b]\to 2$ and $\beta:[c,d]\to 2$ be 2-sequences such that $\alpha\approx\chi^f_i$ and $\beta\approx\alpha x\ol{\alpha}$. Let $t:\mathbb{Z}\to\mathbb{Z}$ be a translation. Then (1) $\alpha\circ t\sub\beta$ implies $t:a\mapsto c$; (2) $\alpha\circ t\sub\ol{\beta}$ implies $t:b\mapsto d$.
\end{lemma}
\begin{proof}
    We prove (1) and (2) together by induction on $i$. Suppose $i=1$. Then we see that $\alpha\approx 001f(0)110$, $\beta\approx 110\ol{f(0)}001x001f(0)110$ and $\ol{\beta}\approx 001{f(0)}110\ol{x}110\ol{f(0)}001$. It is not hard to verify that both (1) and (2) hold. Let $i>1$. Assume $i\in\mathbb{O}$. Then $\alpha\approx\chi^f_{i-1}f(i-1)\ol{\chi^f_{i-1}}$ and $\beta\approx{\chi^f_{i-1}}{f(i-1)}\ol{\chi^f_{i-1}}x\ol{\chi^f_{i-1}}\ol{f(i-1)}{\chi^f_{i-1}}$. Suppose $\alpha\circ t\sub\beta$. Let $\gamma=\alpha\rsto[a,a+|\chi^f_i|-1]$, $\gamma'=\alpha\rsto[b+1-|\chi^f_i|,b]$, $\lambda=\beta\rsto[c,c+|\alpha|-1]$ and $\lambda'=\beta\rsto[d+1-|\alpha|,d]$. Then either $\gamma\circ t\sub\lambda$ or $\gamma'\circ t\sub\lambda'$. By induction hypothesis, we see $t:a\mapsto c$. The case for $\ol{\beta}$ is similar.
\end{proof}

\begin{lemma}\label{lem:TM-sequence-emb}
    Let $f,g\in 2^\omega$ be distinct. Then $\chi^f\npreceq\chi^g$.
\end{lemma}
\begin{proof}
    Since $f\neq g$, there exists $i\in\omega$ such that $f(i)\neq g(i)$ and $f(j)=g(j)$ for all $j<i$. Without loss of generality, assume $f(i)=0$ and $i\in\mathbb{E}$. Then $\chi^f_2\approx\ol{\chi^f_i}1{\chi^f_i}f(i+1){\chi^f_i}0\ol{\chi^f_i}$. By Corollary~\ref{coro:construction-chif}, $\chi^g\preceq\mathsf{Con}(\cset{\chi^g_{i+1}0,\chi^g_{i+1}1,\ol{\chi^g_{i+1}}0,\ol{\chi^g_{i+1}}1})$. Since $g(i)=1$ and $f(j)=g(j)$ for all $j<i$, $\chi^g_{i+1}=\chi^f_i1\ol{\chi^f_i}$. Then we see that $\chi^g\preceq\mathsf{Con}(\cset{\chi^f_i1\ol{\chi^f_i}0,\chi^f_i1\ol{\chi^f_i}1,\ol{\chi^f_i}0{\chi^f_i}0,\ol{\chi^f_i}0{\chi^f_i}1})$. However, by Lemma~\ref{lem:chifi-emb-chifi1-unique}, neither $\chi^f_i1\ol{\chi^f_i}\triangleleft\chi^f_2\approx\ol{\chi^f_i}1{\chi^f_i}f(i+1){\chi^f_i}0\ol{\chi^f_i}$ nor $\ol{\chi^f_i}0\chi^f_i\triangleleft\chi^f_2$ holds. Thus $\chi^f_{i+2}\ntrianglelefteq\chi^g$ and so $\chi^f\npreceq\chi^g$, which concludes the proof.
\end{proof}

\subsection{Umbrellas and their properties}

We are now ready to construct the umbrellas, which are rooted frames for $\mathsf{S4BP}_{2,3}^{2,\omega}$. For each sequence $\alpha:\mathbb{Z}\to 2$, we define the umbrella $\Z_\alpha$ generated by $\alpha$. A key property of umbrellas is that their structures are locally preserved by $k$-t-morphisms, provided $k$ is sufficiently large. The precise statement of this key property is given in Lemma~\ref{lem:injective-map}. With the help of this property, we show that $\mathsf{Log}(\Z_\alpha)$ is pretabular for any finitely perfect $\alpha:\mathbb{Z}\to 2$.

Let us begin with the definition of umbrellas.

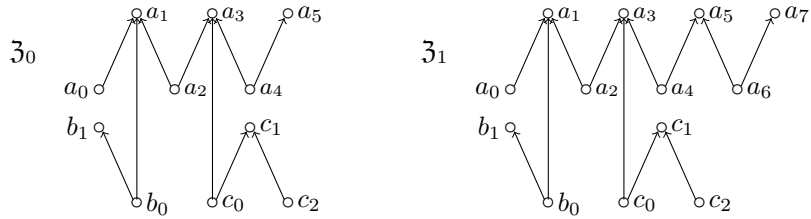
\begin{figure}[ht]
  \[
  \begin{tikzpicture}
    \draw (-2,0.5) node{$\Z_0$};
    \draw (0,0) node{$\circ$};
    \draw (.5,1) node{$\circ$};
    \draw [->] (.05,.05) -- (0.45,.95);
    \draw (1,0) node{$\circ$};
    \draw [->] (.95,.05) -- (0.55,.95);
    \draw (1.5,1) node{$\circ$};
    \draw (1.5,1) node[right]{$a_5$};
    \draw (-1,0) node{$\circ$};
    \draw (-1,0) node[left]{$a_0$};
    \draw [->] (.05-1,.05) -- (0.45-1,.95);
    
    \draw (-.5,1) node{$\circ$};
    \draw [->] (-.05,.05) -- (-0.45,.95);
    \draw [->] (1.05,0.05) -- (1.45,.95);
    \draw (0,0) node[right]{\small $a_2$};
    \draw (-0.5,1) node[right]{\small $a_1$};
    \draw (1,0) node[right]{\small $a_4$};
    \draw (0.5,1) node[right]{\small $a_3$};

    \draw (-.5,-1.5) node{$\circ$};
    \draw (-.5,-1.5) node[right]{$b_0$};
    \draw [->] (-.5,-1.45) -- (-.5,.95);
    \draw [->] (-.55,-1.45) -- (-0.95,-.55);
    \draw (-1,-.5) node{$\circ$};
    \draw (-1,-.5) node[left]{$b_1$};

    \draw (.5,-1.5) node{$\circ$};
    \draw (.5,-1.5) node[right]{$c_0$};
    \draw [->] (.5,-1.45) -- (.5,.95);
    \draw [->] (.55,-1.45) -- (0.95,-.55);
    \draw (1,-.5) node{$\circ$};
    \draw (1,-.5) node[right]{$c_1$};
    \draw (1.5,-1.5) node{$\circ$};
    \draw (1.5,-1.5) node[right]{$c_2$};
    \draw [->] (1.45,-1.45) -- (1.05,-.55);
  \end{tikzpicture}
  \quad\quad\quad
  \begin{tikzpicture}
    \draw (-2,0.5) node{$\Z_1$};
    \draw (0,0) node{$\circ$};
    \draw (.5,1) node{$\circ$};
    \draw [->] (.05,.05) -- (0.45,.95);
    \draw (1,0) node{$\circ$};
    \draw [->] (.95,.05) -- (0.55,.95);
    \draw (1.5,1) node{$\circ$};
    \draw (1.5,1) node[right]{$a_5$};
    \draw (2,0) node{$\circ$};
    \draw (2,0) node[right]{$a_6$};
    \draw [->] (1.95,.05) -- (1.55,.95);
    
    \draw (-.5,1) node{$\circ$};
    \draw [->] (-.05,.05) -- (-0.45,.95);
    \draw (-1,0) node{$\circ$};
    \draw (-1,0) node[left]{$a_0$};
    \draw [->] (-.95,.05) -- (-0.55,.95);
    \draw (2.5,1) node{$\circ$};
    \draw (2.5,1) node[right]{$a_7$};
    \draw [->] (1.05,0.05) -- (1.45,.95);
    \draw [->] (1.05+1,0.05) -- (1.45+1,.95);
    \draw (0,0) node[right]{\small $a_2$};
    \draw (-0.5,1) node[right]{\small $a_1$};
    \draw (1,0) node[right]{\small $a_4$};
    \draw (0.5,1) node[right]{\small $a_3$};

    \draw (-.5,-1.5) node{$\circ$};
    \draw (-.5,-1.5) node[right]{$b_0$};
    \draw [->] (-.5,-1.45) -- (-.5,.95);
    \draw [->] (-.55,-1.45) -- (-0.95,-.55);
    \draw (-1,-.5) node{$\circ$};
    \draw (-1,-.5) node[left]{$b_1$};

    \draw (.5,-1.5) node{$\circ$};
    \draw (.5,-1.5) node[right]{$c_0$};
    \draw [->] (.5,-1.45) -- (.5,.95);
    \draw [->] (.55,-1.45) -- (0.95,-.55);
    \draw (1,-.5) node{$\circ$};
    \draw (1,-.5) node[right]{$c_1$};
    \draw (1.5,-1.5) node{$\circ$};
    \draw (1.5,-1.5) node[right]{$c_2$};
    \draw [->] (1.45,-1.45) -- (1.05,-.55);
  \end{tikzpicture}
  \]
  \caption{The frames $\Z_0$ and $\Z_1$}
  \label{fig:BulidingBlocks}
\end{figure}

\begin{definition}
    Let $\Z_0=(Z_0,R_0)$ and $\Z_1=(Z_1,R_1)$ be frames as depicted in Figure~\ref{fig:BulidingBlocks}. To be precise, we define $Z_0=\cset{a_i:i<6}\cup\cset{b_0,b_1}\cup\cset{c_i:i<2}$, $Z_1=Z_0\cup\cset{a_6,a_7}$, $R_0=R_1\rsto Z_0$ and $R_1$ to be the transitive-reflexive closure of the set
    \begin{center}
        $\cset{\tup{a_{2i},a_{2i+1}}:i<4}\cup\cset{\tup{a_{2i+2},a_{2i+1}}:i<3}\cup\cset{\tup{b_0,b_1},\tup{b_0,a_1},\tup{c_0,c_1},\tup{c_2,c_1},\tup{c_0,a_3}}$.
    \end{center}

    For each $\alpha:\mathbb{Z}\to 2$, we define the frame $\Z_\alpha=(Z_\alpha,R_\alpha)$ as follows:
    \begin{itemize}
        \item $Z_\alpha=\biguplus_{i\in\mathbb{Z}}Z_{\alpha(i)}=\cset{\tup{x,i}:x\in Z_{\alpha(i)}\text{ and }i\in\mathbb{Z}}$
        \item $\tup{x,i}R_{\alpha}\tup{y,j}$ if and only if one of the following holds:
        \begin{itemize}
            \item $i=j$ and $R_{\alpha(i)}xy$;
            \item $j=i+1$, $f(i)=0$, $x=a_5$ and $y=a_0$;
            \item $j=i+1$, $f(i)=1$, $x=a_7$ and $y=a_0$.
        \end{itemize}
    \end{itemize}    
\end{definition}

\begin{figure}[ht]
  \[
  \begin{tikzpicture}
    \draw (2.5,0.5) node{$\cdots$};
    \draw (0,0) node{$\circ$};
    \draw (.5,1) node{$\circ$};
    \draw [->] (.05,.05) -- (0.45,.95);
    \draw (1,0) node{$\circ$};
    \draw [->] (.95,.05) -- (0.55,.95);
    \draw (1.5,1) node{$\circ$};
    \draw (1.5,1) node[right]{$a^5_5$};
    \draw (2,0) node{$\circ$};
    \draw (2,0) node[right]{$a^6_6$};
    \draw [->] (1.95,.05) -- (1.55,.95);
    
    \draw (-.5,1) node{$\circ$};
    \draw [->] (-.05,.05) -- (-0.45,.95);
    \draw (-1,0) node{$\circ$};
    \draw (-1,0) node[right]{$a^0_0$};
    \draw [->] (-.95,.05) -- (-0.55,.95);
    \draw [->] (1.05,0.05) -- (1.45,.95);
    \draw (0,0) node[right]{$a^2_2$};
    \draw (-0.5,1) node[right]{$a^1_1$};
    \draw (1,0) node[right]{$a^4_4$};
    \draw (0.5,1) node[right]{$a^3_3$};

    \draw (-.5,-1.5) node{$\circ$};
    \draw (-.5,-1.5) node[right]{$b^0_1$};
    \draw [->] (-.5,-1.45) -- (-.5,.95);
    \draw [->] (-.55,-1.45) -- (-0.95,-.55);
    \draw (-1,-.5) node{$\circ$};
    \draw (-1,-.5) node[left]{$b^1_1$};

    \draw (.5,-1.5) node{$\circ$};
    \draw (.5,-1.5) node[right]{$c^0_3$};
    \draw [->] (.5,-1.45) -- (.5,.95);
    \draw [->] (.55,-1.45) -- (0.95,-.55);
    \draw (1,-.5) node{$\circ$};
    \draw (1,-.5) node[right]{$c^1_3$};
    \draw (1.5,-1.5) node{$\circ$};
    \draw (1.5,-1.5) node[right]{$c^2_3$};
    \draw [->] (1.45,-1.45) -- (1.05,-.55);

    \draw (-4+-1.5,0.5) node{$\cdots$};
    \draw (-4+0,0) node{$\circ$};
    \draw (-4+.5,1) node{$\circ$};
    \draw [->] (-4+.05,.05) -- (-4+0.45,.95);
    \draw (-4+1,0) node{$\circ$};
    \draw [->] (-4+.95,.05) -- (-4+0.55,.95);
    \draw (-4+1.5,1) node{$\circ$};
    \draw (-4+1.5,1) node[right]{$a^5_{-3}$};
    \draw (-3+1.5,1) node{$\circ$};
    \draw (-3+1.5,1) node[right]{$a^7_{-1}$};
    \draw (-4+2,0) node{$\circ$};
    \draw [->] (-4+1.95,.05) -- (-4+1.55,.95);
    \draw [->] (-3+1.95,.05) -- (-3+1.55,.95);
    
    \draw (-4+-.5,1) node{$\circ$};
    \draw [->] (-4+-.05,.05) -- (-4+-0.45,.95);
    \draw (-4+-1,0) node{$\circ$};
    \draw [->] (-4+-.95,.05) -- (-4+-0.55,.95);
    \draw [->] (-4+1.05,0.05) -- (-4+1.45,.95);
    \draw [->] (-3+1.05,0.05) -- (-3+1.45,.95);
    \draw (-4+0,0) node[below]{$a^2_{-6}$};
    \draw (-4+-0.5,1) node[right]{$a^1_{-7}$};
    \draw (-4+-1,0) node[left]{$a^0_{-8}$};
    \draw (-3+1,0) node[right]{$a^6_{-2}$};
    \draw (-4+1,0) node[right]{$a^4_{-4}$};
    \draw (-4+0.5,1) node[right]{$a^3_{-5}$};

    \draw (-4+-.5,-1.5) node{$\circ$};
    \draw (-4+-.5,-1.5) node[right]{$b^0_{-7}$};
    \draw [->] (-4+-.5,-1.45) -- (-4+-.5,.95);
    \draw [->] (-4+-.55,-1.45) -- (-4+-0.95,-.55);
    \draw (-4+-1,-.5) node{$\circ$};
    \draw (-4+-1,-.5) node[left]{$b^1_{-7}$};

    \draw (-4+.5,-1.5) node{$\circ$};
    \draw (-4+.5,-1.5) node[right]{$c^0_{-5}$};
    \draw [->] (-4+.5,-1.45) -- (-4+.5,.95);
    \draw [->] (-4+.55,-1.45) -- (-4+0.95,-.55);
    \draw (-4+1,-.5) node{$\circ$};
    \draw (-4+1,-.5) node[right]{$c^1_{-5}$};
    \draw (-4+1.5,-1.5) node{$\circ$};
    \draw (-4+1.5,-1.5) node[right]{$c^2_{-5}$};
    \draw [->] (-4+1.45,-1.45) -- (-4+1.05,-.55);
  \end{tikzpicture}
  \]
  \caption{Fragment of relabelled $\Z_\alpha$}
  \label{fig:relabel}
\end{figure}
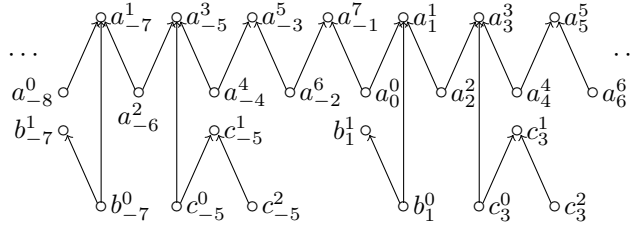

In what follows, let $\alpha:\mathbb{Z}\to 2$ be an arbitrarily fixed sequence on $\mathbb{Z}$ and $\Z_\alpha=(Z,R)$. Our aim is to show that $\mathsf{Log}(\Z_\alpha)$ is pretabular. To make the proofs below easier to read, we re-indexed the elements in $Z$ by a onto map $l:Z\to\mathbb{Z}$ where:
\begin{itemize}
    \item $l(\tup{a_0,0})=0$;
    \item for all $\tup{a_i,j},\tup{a_{i'},j'}\in Z$, $l(\tup{a_i,j})<l(\tup{a_{i'},j'})$ iff $j<j'$ or, $j=j'$ and $i<i'$;
    \item for all $i<2$ and $j\in\mathbb{Z}$, $l(\tup{b_i,j})=l(\tup{a_1,j})$;
    \item for all $i<3$ and $j\in\mathbb{Z}$, $l(\tup{c_i,j})=l(\tup{a_3,j})$.
\end{itemize}
It is easy to see such map $f$ exists and is unique. To simplify notations, we write $a^i_{l(\tup{x_i,j})}$, $b^i_{l(\tup{y_i,j})}$ and $c^i_{l(\tup{z_i,j})}$ for $\tup{a_i,j}$, $\tup{b_i,j}$ and $\tup{c_i,j}$, respectively. We also write $a_j$ for $a^i_j$. A fragment of the re-indexed frame $\Z_\alpha$ is shown in Figure~\ref{fig:relabel}. For all $i,j\in\mathbb{Z}$ such that $i\leq j$, we define $Z[i,j]=\cset{z\in Z: i\leq l(z)\leq j}$.

Note that for all $a^{i_1}_{j_1},b^{i_2}_{j_2},c^{i_3}_{j_3}\in Z$, we have $i_1+j_1\in\mathbb{E}$, and $j_2,j_3\in\mathbb{O}$.

\begin{lemma}\label{lem:BoundaryOfIntervalsInZ-alpha}
    Let $\G=(Y,S)\in\mathsf{Fr}_r$ and $f:(\Z_\alpha,z_0)\to^{n}(\G,y_0)$. Then for all $i\leq j$ such that $Z[i,j]\sub\R^{n}[z_0]$, $Z[i,j]$ is sufficient if $f[\cset{x_i,x_j}]\sub f[Z[i,j]\setminus\cset{x_i,x_j}]$.
\end{lemma}
\begin{proof}
    Note that $\R[Z[i,j]\setminus\cset{x_i,x_j}]\sub Z[i,j]$. By Fact~\ref{fact:boundary-sufficient}, $f[\cset{x_i,x_j}]\sub f[Z[i,j]\setminus\cset{x_i,x_j}]$ implies that $Z[i,j]$ is sufficient.
\end{proof}

To simplify proofs of the following lemmas, we present the following proposition. 

\begin{proposition}\label{prop:k-t-morphism}
    Let $\F=(X,R)$, $\F'=(X',R')$ be frames, $x\in X$, $x'\in X'$ and $f:(\F,x)\to^{k}(\F',x')$. Let $m<k$. Suppose $z,z'\in \R^{m}[x]$ and $f(z)=f(z')$. Then 
    \begin{center}
        $f[\R^n[z]]=f[{R'}_\sharp^n[z']]$ for all $n\in\omega$ such that $n+m<k$.
    \end{center}
\end{proposition}
\begin{proof}
    The proof proceeds by induction on $n$. The case $n=0$ is trivial. Let $n>0$. Take any $w\in\R^{n}[z]$. Then there exists $u\in\R^{n-1}[z]$ such that $w\in\R[u]$. By induction hypothesis, $f(u)\in f[\R^{n-1}[z]]=f[\R^{n-1}[z']]$ and so $f(u)=f(u')$ for some $u'\in\R^{n-1}[z']$. Note that $u,u'\in\R^{k-1}[x]$, we have $f[\R[u]]={R'}_\sharp[f(u)]={R'}_\sharp[f(u')]=f[\R[u']]$. Thus $f(w)\in f[\R[u]]=f[\R[u']]\sub f[\R^{n}[z']]$. By arbitrariness of $w$, $f[\R^n[z]]\sub f[\R^n[z']]$.
    Symmetrically, we see $f[\R^n[z']]\sub f[\R^n[z]]$. Hence $f[\R^n[z]]=f[\R^n[z']]$.
\end{proof}

Intuitively, Proposition~\ref{prop:k-t-morphism} shows that if two points have the same image under a local t-morphism $f$, then their neighborhoods are `equivalent' with respect to $f$.

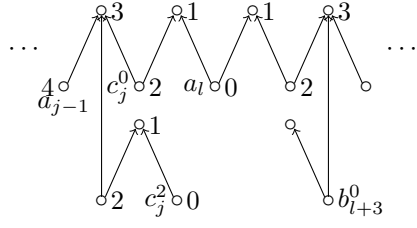
\begin{figure}[ht]
  \[
  \begin{tikzpicture}
    \draw (-2,0) node[left]{$a_l$};
    \draw (-4+1.5,-1.5) node[left]{$c^2_j$};
    \draw (.5,0.5) node{$\cdots$};
    \draw (0,0) node{$\circ$};
    
    \draw (-.5,1) node{$\circ$};
    \draw [->] (-.05,.05) -- (-0.45,.95);
    \draw (-1,0) node{$\circ$};
    \draw (-1,0) node[right]{$2$};
    \draw [->] (-.95,.05) -- (-0.55,.95);
    \draw (-0.5,1) node[right]{$3$};
    \draw (-0.5,-1.5) node[right]{$b^0_{l+3}$};

    \draw (-.5,-1.5) node{$\circ$};
    \draw [->] (-.5,-1.45) -- (-.5,.95);
    \draw [->] (-.55,-1.45) -- (-0.95,-.55);
    \draw (-1,-.5) node{$\circ$};

    \draw (-4.5,0.5) node{$\cdots$};
    \draw (-4+0,0) node{$\circ$};
    \draw (-4+.5,1) node{$\circ$};
    \draw [->] (-4+.05,.05) -- (-4+0.45,.95);
    \draw (-4+1,0) node{$\circ$};
    \draw [->] (-4+.95,.05) -- (-4+0.55,.95);
    \draw (-4+1.5,1) node{$\circ$};
    \draw (-4+1.5,1) node[right]{$1$};
    \draw (-3+1.5,1) node{$\circ$};
    \draw (-3+1.5,1) node[right]{$1$};
    \draw (-4+2,0) node{$\circ$};
    \draw [->] (-4+1.95,.05) -- (-4+1.55,.95);
    \draw [->] (-3+1.95,.05) -- (-3+1.55,.95);
    \draw [->] (-4+1.05,0.05) -- (-4+1.45,.95);
    \draw [->] (-3+1.05,0.05) -- (-3+1.45,.95);
    \draw (-4+0,0) node[left]{$4$};
    \draw (-4+0,0) node[below]{$a_{j-1}$};
    \draw (-3+1,0) node[right]{$0$};
    \draw (-4+1,0) node[right]{$2$};
    \draw (-4+1,0) node[left]{$c^0_j$};
    \draw (-4+0.5,1) node[right]{$3$};

    \draw (-4+.5,-1.5) node{$\circ$};
    \draw (-4+.5,-1.5) node[right]{$2$};
    \draw [->] (-4+.5,-1.45) -- (-4+.5,.95);
    \draw [->] (-4+.55,-1.45) -- (-4+0.95,-.55);
    \draw (-4+1,-.5) node{$\circ$};
    \draw (-4+1,-.5) node[right]{$1$};
    \draw (-4+1.5,-1.5) node{$\circ$};
    \draw (-4+1.5,-1.5) node[right]{$0$};
    \draw [->] (-4+1.45,-1.45) -- (-4+1.05,-.55);
  \end{tikzpicture}
  \]
  \caption{For Claim 3 in the proof of Lemma~\ref{lem:injective-map}}
  \label{fig:aux-proof-3}
\end{figure}

\begin{lemma}\label{lem:Zalpha-k-morphism-same-type}
    Let $\G=(Y,S)\in\mathsf{Fr}_r(\mathsf{S4}_t)$, $z_0\in Z$, $y_0\in Y$ and $n\in\omega$. Suppose $f:(\Z_\alpha,z_0)\to^{n+19}(\G,y_0)$ and $f$ is not sufficient. Then for all $x,y\in\cset{a,b,c}$ and $i,j,k,l\in\mathbb{Z}$ such that $f(x^i_j)=f(y^k_l)$, the following holds:
    \begin{enumerate}[(1)]
        \item if $x^i_j,y^k_l\in\R^{n+18}[z_0]$, then $i+k\in\mathbb{E}$;
        \item if $x^i_j,y^k_l\in\R^{n+16}[z_0]$, then $x=b$ implies $y\neq c$;
        \item if $x^i_j,y^k_l\in\R^{n+6}[z_0]$, then $x=b$ implies $y\neq a$;
        \item if $x^i_j,y^k_l\in\R^{n}[z_0]$, then $x=c$ implies $y\neq a$.
    \end{enumerate}
\end{lemma}
\begin{proof}
    For (1), suppose $x^i_j,y^k_l\in\R^{n+18}[z_0]$ and $i+k\not\in\mathbb{E}$. Assume $i\in\mathbb{E}$. Then for each $z\in\cset{x^i_j,y^k_l}$, we have $f(x)=f(x^i_j)=f(y^k_l)$ and $R[y^k_l]\cup\breve{R}[x^i_j]\sub\cset{x^i_j,y^k_l}$. By Fact~\ref{fact:boundary-sufficient}(1), $\cset{x^i_j,y^k_l}$ is sufficient. Assume $i\in\mathbb{O}$. Then $\breve{R}[y^k_l]\cup R[x^i_j]\sub\cset{x^i_j,y^k_l}$, which also implies that $\cset{x^i_j,y^k_l}$ is sufficient.
    \vspace{.5em}

    \noindent For (2), suppose $x^i_j,y^k_l\in\R^{n+16}[z_0]$, $x=b$ and $y=c$. Then we have two cases:

    (2.1) $i=1$. By (1), $k=1$. By Proposition~\ref{prop:k-t-morphism}, we see $\cset{f(b^0_j),f(b^1_j)}=f[\R[b^1_j]]=f[\R[c^1_l]]=\cset{f(c^0_l),f(c^1_l),f(c^2_l)}$, which entails $f(c^0_l)=f(c^2_l)=f(b^0_j)$. Let $X=\cset{b^0_j,b^1_j,c^2_l,c^1_l}$. Then $\R[b^1_j]\cup\R[c^2_l]\sub X$. Note that $X\sub\R^{n+18}[z_0]$, $f(c^2_l)=f(b^0_j)$ and $f(c^1_l)=f(b^1_j)$, by Fact~\ref{fact:boundary-sufficient}, we see $X$ is sufficient. 
    
    (2.2) $i\neq 1$. Then $i=0$ and so $k\in\cset{0,2}$. By (2.1), $f(b^1_j)\neq f(c^1_j)$. If $k=2$, then by Proposition~\ref{prop:k-t-morphism} and (1), $f(b^1_j)=f(c^1_j)$, which is impossible. Thus $k=0$. By Proposition~\ref{prop:k-t-morphism} and (1), we see $f(a^1_j)=f(c^1_l)$ and $f(a^1_l)=f(b^1_j)$. Consider the set $X'=Z[j,j]\cup Z[l,l]$. Clearly, $X'\sub\R^2[x^i_j]\cup\R^2[y^k_l]\sub\R^{n+18}[z_0]$. Note that $\R[X'\setminus\cset{a^1_j,a^1_l}]\sub X'$ and $f[\cset{a^1_j,a^1_l}]\sub f[X'\setminus\cset{a^1_j,a^1_l}]$, we see $X'$ is sufficient. Thus this case is impossible.
    \vspace{.5em}
    
    \noindent For (3), let $x=b$ and $y=a$. We first prove the following two claims: 
    
    \noindent\textbf{Claim 1:} Suppose $i=1$. If $x^i_j,y^k_l\in\R^{n+15}[z_0]$, then $j=l+2$.
    
    \noindent\textbf{Proof of Claim 1:} By $x^i_j,y^k_l\in\R^{n+15}[z_0]$, we see $Z[\min(j,l),\max(j,l)]\sub\R^{n+18}[z_0]$. Suppose $l=j$. Then $f(b^1_j)=f(a^1_j)$. By Lemma~\ref{lem:BoundaryOfIntervalsInZ-alpha}, $Z[j,j]$ is sufficient. Suppose $l=j+2$. Then by Proposition~\ref{prop:k-t-morphism} and (1), we have $f(b^0_j)=f(c^0_l)$, which contradicts (2.1). Suppose $l>j+2$. By (1) and Proposition~\ref{prop:k-t-morphism}, we see $f(b^0_j)=f(a_{l-1})$ and $f(a_j)=f(a_{l-2})$. Note that $a_{l-2}\neq a_j$, we see $f[\cset{a_j,a_l}]\sub f[Z[j,l]\setminus\cset{a_j,a_l}]$. By Lemma~\ref{lem:BoundaryOfIntervalsInZ-alpha}, $Z[j,l]$ is sufficient. Suppose $l+2<j$. By (1) and Proposition~\ref{prop:k-t-morphism}, we see $f(b^0_j)=f(a_{l+1})$ and $f(a_j)=f(a_{l+2})$. Note that $a_{l+2}\neq a_j$, we see $f[\cset{a_j,a_l}]\sub f[Z[l,j]\setminus\cset{a_j,a_l}]$. By Lemma~\ref{lem:BoundaryOfIntervalsInZ-alpha}, $Z[l,j]$ is sufficient.\hfill$\dashv$

    \vspace{.5em}
    \noindent\textbf{Claim 2:} Suppose $i=1$. If $x^i_j,y^k_l\in\R^{n+7}[z_0]$, then $j\neq l+2$.
    
    \noindent\textbf{Proof of Claim 2:} Suppose $j=l+2$. Then $k\in\cset{5,7}$. 
    
    (a) $k=5$. Consider the frame $\Z_\alpha$ with labels in Figure~\ref{fig:aux-proof-1}. By Proposition~\ref{prop:k-t-morphism} and (1), points with same label have the same $f$-image. Then $f(c^0_{l-2})\in f[\cset{a_{j+1},b^0_j}]$. By (2), $f(c^0_{l-2})\neq f(b^0_j)$ and so $f(a_{j+1})=f(c^0_{l-2})$. Note that $f(a_{l-2})=f(a_j)$ and $Z[l-2,j+1]\sub\R^{n+18}[z_0]$, by Lemma~\ref{lem:BoundaryOfIntervalsInZ-alpha}, $Z[l-2,j+1]$ is sufficient, which is impossible.

    (b) $k=7$. Consider the frame $\Z_\alpha$ with labels given in Figure~\ref{fig:aux-proof-2}(a). By Proposition~\ref{prop:k-t-morphism} and (1), we see that points with same label have the same $f$-image. Moreover, we see that $f(c^0_{l-4})\in f[\cset{a_{j+3},c^0_{j+2}}]$.
    Suppose $f(c^0_{l-4})=f(a_{j+3})$. Note that $Z[l-4,j+3]\sub\R^{n+18}[z_0]$, by Lemma~\ref{lem:BoundaryOfIntervalsInZ-alpha}, $Z[l-4,j+3]$ is sufficient. Suppose $f(c^0_{l-4})=f(c^0_{j+2})$. By (1) and Proposition~\ref{prop:k-t-morphism}, we can verify that in Figure~\ref{fig:aux-proof-2}(b), points with same label have the same $f$-image. Then $f(b^1_{l-6})\in f[\cset{{a_{j+2},a_{j+6}}}]$. Note that $\cset{a_{j+2},a_{j+6},b^1_{l-6}}\sub\R^{n+15}[z_0]$, by Claim 1, $f(b^1_{l-6})\not\in f[\cset{{a_{j+2},a_{j+6}}}]$.\hfill$\dashv$
    \vspace{.5em}

    Suppose $x^i_j,y^k_l\in\R^{n+6}[z_0]$. By Claim 1 and Claim 2, $i\neq 1$. Then $i=0$. By (1) and Proposition~\ref{prop:k-t-morphism}, $f(b^1_j)\in\cset{f(a_{l-1}),f(a_{l+1})}$. Note that $\cset{b^1_j,a_{l-1},a_{l+1}}\sub\R^{n+7}[z_0]$, by Claim 1 and Claim 2, $f(b^1_j)\not\in\cset{f(a_{l-1}),f(a_{l+1})}$, which is impossible.
    \vspace{.5em}

    \noindent For (4), let $x=c$ and $y=a$. We first prove the following claims:

    \noindent\textbf{Claim 3:} Suppose $x^i_j,y^k_l\in\R^{n+2}[z_0]$. then $i\neq 2$.
    
    \noindent\textbf{Proof of Claim 3:} suppose $x^i_j,y^k_l\in\R^{n+2}[z_0]$ and $i=2$. By (1), $j+l\in\mathbb{O}$. Suppose $l=j+1$. By Proposition~\ref{prop:k-t-morphism} and (1), $f(a_l)=f(c^1_l)$. By Lemma~\ref{lem:BoundaryOfIntervalsInZ-alpha}, $Z[j,l]$ is sufficient. 
    Similarly, $l=j-1$ implies $Z[l,j]$ is sufficient. 
    Suppose $l>j+3$. By Proposition~\ref{prop:k-t-morphism} and (1), $f(a_l)=f(a_{j-3})$. By Lemma~\ref{lem:BoundaryOfIntervalsInZ-alpha}, $Z[j,l]$ is sufficient. Similarly, $l<j-3$ implies $Z[l,j]$ is sufficient. 
    Suppose $l=j-3$. By Proposition~\ref{prop:k-t-morphism} and (1), $f(b^0_{l-1})\in f[\R^2[a_l]]=f[\R^2[c^2_j]]=f[\cset{c^0_j,c^1_j,c^2_j}]$. Since $\cset{b^0_{l-1},c^0_j,c^1_j,c^2_j}\sub\R^{n+16}[z_0]$, by (2), $f(b^0_{l-1})\not\in f[\cset{c^0_j,c^1_j,c^2_j}]$.
    Suppose $l=j+3$. Then $k=0$ or $k=6$. If $k=0$, then by (1) and Proposition~\ref{prop:k-t-morphism}, we see $f(b^0_{l+1})\in f[\R^2[a_l]]=f[\R^2[c^2_j]]=f[\cset{c^0_j,c^1_j,c^2_j}]$, which contradicts (2). Suppose $k=6$. Consider the relabelled frame in Figure~\ref{fig:aux-proof-3}. By Proposition~\ref{prop:k-t-morphism} and (1), points with same label have the same $f$-image. Thus $f(b^0_{l+3})\in\cset{f(c^0_j),f(a_{j-1})}$. Note that $\cset{b^0_{l+3},c^0_j,a_{j-1}}\sub\R^{n+6}[z_0]$, by (2) and (3), $f(b^0_{l+3})\not\in\cset{f(c^0_j),f(a_{j-1})}$.\hfill$\dashv$

    \noindent\textbf{Claim 4:} Suppose $x^i_j,y^k_l\in\R^{n+1}[z_0]$. then $i\neq 1$.
    
    \noindent\textbf{Proof of Claim 4:} suppose $x^i_j,y^k_l\in\R^{n+1}[z_0]$ and $i=1$. Note that $\R[a^k_l]\sub\R^{n+2}[z_0]$, by (1) and Proposition~\ref{prop:k-t-morphism}, $f(c^2_j)\in f[\R[a^k_l]]$. By (2) and (3), $k=3$ and $f(c^0_l)=f(c^2_j)$. By Proposition~\ref{prop:k-t-morphism}, $f(c^1_j)=f(c^1_l)$, which entails that $\cset{c^0_l,c^1_l,c^2_l,c^1_j,c^2_j}$ is sufficient.\hfill$\dashv$
    \vspace{.5em}

    Suppose $x^i_j,y^k_l\in\R^{n}[z_0]$. By Claim 3 and Claim 4, $i\not\in\cset{1,2}$. Then $i=0$. By (1) and Proposition~\ref{prop:k-t-morphism}, $f(c^1_j)\in\cset{f(a_{l-1}),f(a_{l+1})}$. Note that $\cset{c^1_j,a_{l-1},a_{l+1}}\sub\R^{n+1}[z_0]$, by Claim 3 and Claim 4, $f(c^1_j)\not\in\cset{f(a_{l-1}),f(a_{l+1})}$, which is impossible.
\end{proof}

\begin{lemma}\label{lem:Zalpha-k-morphism-same-label}
    Let $\G=(Y,S)\in\mathsf{Fr}_r(\mathsf{S4}_t)$, $z_0\in Z$, $y_0\in Y$ and $n\in\omega$. Suppose $f:(\Z_\alpha,z_0)\to^{n+26}(\G,y_0)$ and $f$ is not sufficient. Then for all $x,y\in\cset{a,b,c}$ and $i,j,k,l\in\mathbb{Z}$ such that $f(x^i_j)=f(y^k_l)$, the following holds:
    \begin{enumerate}[(1)]
        \item if $x^i_j,y^k_l\in\R^{n+4}[z_0]$, then $i=k$;
        \item if $x^i_j,y^k_l\in\R^{n}[z_0]$, then $j=l$.
    \end{enumerate}
\end{lemma}
\begin{proof}
    Take any $x^i_j,y^k_l\in\R^n[z_0]$ such that $f(x^i_j)=f(y^k_l)$. By Lemma~\ref{lem:Zalpha-k-morphism-same-type}, $x=y$ and $i+k\in\mathbb{E}$.
    For (1), consider the following three cases:

    (1.1) $x=b$. Then $i,k<2$. Since $i+k\in\mathbb{E}$, we see $i=k$.

    (1.2) $x=c$. Then $i,k\in\cset{0,1,2}$. Suppose $i\neq k$. Then $\cset{i,k}=\cset{0,2}$. Suppose $i=0$. By Proposition~\ref{prop:k-t-morphism}, $f(a_i)\in f[\R[c^0_j]]=f[\R[c^2_l]]=\cset{f(c^2_l),f(c^1_l)}$. Since $f:(\Z_\alpha,z_0)\to^{n+26}(\G,y_0)$ and $\cset{a_i,c^2_l,c^1_l}\sub\R^{n+7}[z_0]$, by Lemma~\ref{lem:Zalpha-k-morphism-same-type}(4), $f(a_i)\not\in\cset{f(c^0_l),f(c^1_l)}$, which leads to a contradiction. Symmetrically, $i=0$ is also impossible.

    (1.3) $x=a$. Then $i,k\leq 7$. We now show $i=k$ by showing the following claims:
    
    \noindent\textbf{Claim 1:} If $x^i_j,y^k_l\in\R^{n+5}[z_0]$, then $i\in\mathbb{E}$ implies $i=k$. 
    
    \noindent\textbf{Proof of Claim 1:} Consider the following cases: 
    
    (a) $2\in\cset{i,k}$. Suppose $i=2$. By Proposition~\ref{prop:k-t-morphism}, $f(b^0_{j-1}),f(c^0_{j+1})\in f[\R^2[a^k_l]]$. Since $\R^2[a^k_l]\sub\R^{n+7}[z_0]$, by Lemma~\ref{lem:Zalpha-k-morphism-same-type}, $k\in\mathbb{E}$ and there are $b^{k_1}_{l_1},c^{k_2}_{l_2}\in\R^2[a^k_l]$, which entails $k=2$. Symmetrically, $k=2$ implies $i=2$. 
    
    (b) $0\in\cset{i,k}$. Suppose $i=0$. By Proposition~\ref{prop:k-t-morphism}, $f(b^0_{j+1})\in f[\R^2[a^k_l]]$, which entails $k\in\cset{0,1,2}$. By Lemma~\ref{lem:Zalpha-k-morphism-same-type}(1) and (a), $k=0$. Symmetrically, $k=0$ implies $i=k$. 
    
    (c) $4\in\cset{i,k}$. Suppose $i=4$. Similar to the argument for (b).
    
    Note that $i+k\in\mathbb{E}$ and $i,k\leq 7$, by (a)-(c), we see $i=6$ if and only if $k=6$.\hfill$\dashv$

    \noindent\textbf{Claim 2:} If $x^i_j,y^k_l\in\R^{n+4}[z_0]$, then $i=k$.
    
    \noindent\textbf{Proof of Claim 2:} The case $i\in\mathbb{E}$ follows from Claim 1 immediately. Suppose $i\in\mathbb{O}$. By Proposition~\ref{prop:k-t-morphism} and the definition of $Z$, $f[\cset{a^{i-1}_{j-1},a^{i'}_{j+1}}]\sub f[\cset{a^{k-1}_{l-1},a^{k'}_{l+1}}]$. Note that $\cset{a_{j-1},a_{j+1}}\sub\R^{n+5}[z_0]$, by Claim 1, $\cset{i-1,i'}=\cset{k-1,k'}$. Suppose $i\neq k$. Then $i-1=k'\in\cset{k+1,0}$ and $k-1=i'\in\cset{i+1,0}$. If $i-1=0$, then $k\in\cset{5,7}$, which contradicts $k-1\in\cset{i+1,0}=\cset{2,0}$. Thus $i-1=k+1$ and so $k-1=0$. Then $i\in\cset{5,7}$, which contradicts $i-1\in\cset{k+1,0}=\cset{2,0}$. Hence $i=k$.\hfill$\dashv$

    For (2), suppose $j<l$. By Lemma~\ref{lem:Zalpha-k-morphism-same-type} and (1), $x=y$ and $i=k$. By Proposition~\ref{prop:k-t-morphism}, $f(a_j)=f(a_l)$. Note that $a_j,a_l\in\R^{3}[x^i_j]\cup\R^3[y^k_l]$, we see $\R[a_{j}]\cup\R[a_{l}]\sub\R^{n+4}[z_0]$. Then by (1), $f(a_{j-1})=f(a_{l-1})$. Thus $Z[j-1,l]$ is sufficient. Symmetrically, $j>l$ implies that $Z[l-1,j]$ is sufficient. By assumption, $f$ is not sufficient, which leads to a contradiction. Hence $j=l$.
\end{proof}

\begin{lemma}\label{lem:injective-map}
    Let $\G=(Y,S)\in\mathsf{Fr}_r$ be an infinite frame. Let $z_0\in Z$, $y_0\in Y$ and $n\in\omega$. Suppose $f:(\Z_\alpha,z_0)\to^{n+26}(\G,y_0)$. Then $g=f\rsto\R^n[z_0]$ is injective.
\end{lemma}
\begin{proof}
    Take any $z_1,z_2\in\mathsf{dom}(g)$. Then there exists  $x,y\in\cset{a,b,c}$ and $i,j,k,l\in\mathbb{Z}$ such that $z_1=x^i_j$ and $z_2=y^k_l$. Suppose $f(z_1)=f(z_2)$. Then by Lemmas~\ref{lem:Zalpha-k-morphism-same-type} and \ref{lem:Zalpha-k-morphism-same-label}, $x=y$, $i=k$ and $j=l$, which entails $z_1=z_2$ immediately.
\end{proof}

\begin{figure}[ht]
  \[
  \begin{tikzpicture}
    \draw (2.5,0.5) node{$\cdots$};
    \draw (0,0) node{$\circ$};
    \draw (.5,1) node{$\circ$};
    \draw [->] (.05,.05) -- (0.45,.95);
    \draw (1,0) node{$\circ$};
    \draw [->] (.95,.05) -- (0.55,.95);
    \draw (1.5,1) node{$\circ$};
    \draw (2,0) node{$\circ$};
    \draw [->] (1.95,.05) -- (1.55,.95);
    
    \draw (-.5,1) node{$\circ$};
    \draw [->] (-.05,.05) -- (-0.45,.95);
    \draw (-1,0) node{$\circ$};
    \draw (-1,0) node[right]{$1$};
    \draw [->] (-.95,.05) -- (-0.55,.95);
    \draw [->] (1.05,0.05) -- (1.45,.95);
    \draw (0,0) node[right]{$3$};
    \draw (-0.5,1) node[right]{$2$};

    \draw (-.5,-1.5) node{$\circ$};
    \draw (-.5,-1.5) node[right]{$1$};
    \draw [->] (-.5,-1.45) -- (-.5,.95);
    \draw [->] (-.55,-1.45) -- (-0.95,-.55);
    \draw (-1,-.5) node{$\circ$};
    \draw (-1,-.5) node[below]{$0$};
    \draw (-1,-.5) node[right]{$b^1_j$};

    \draw (.5,-1.5) node{$\circ$};
    \draw [->] (.5,-1.45) -- (.5,.95);
    \draw [->] (.55,-1.45) -- (0.95,-.55);
    \draw (1,-.5) node{$\circ$};
    \draw (1.5,-1.5) node{$\circ$};
    \draw [->] (1.45,-1.45) -- (1.05,-.55);

    \draw (-3+-1.5,0.5) node{$\cdots$};
    \draw (-3+0,0) node{$\circ$};
    \draw (-3+.5,1) node{$\circ$};
    \draw [->] (-3+.05,.05) -- (-3+0.45,.95);
    \draw (-3+1,0) node{$\circ$};
    \draw [->] (-3+.95,.05) -- (-3+0.55,.95);
    \draw (-3+1.5,1) node{$\circ$};
    \draw (-3+1.5,1) node[right]{$0$};
    \draw (-3+1.5,1) node[above]{$a_l$};
    \draw (-3+2,0) node{$\circ$};
    \draw [->] (-3+1.95,.05) -- (-3+1.55,.95);
    
    \draw (-3+-.5,1) node{$\circ$};
    \draw [->] (-3+-.05,.05) -- (-3+-0.45,.95);
    \draw (-3+-1,0) node{$\circ$};
    \draw [->] (-3+-.95,.05) -- (-3+-0.55,.95);
    \draw [->] (-3+1.05,0.05) -- (-3+1.45,.95);
    \draw (-3+0,0) node[below]{$a_{l-3}$};
    \draw (0,0) node[below]{$a_{j+1}$};
    \draw (-3+1,0) node[right]{$1$};
    \draw (-3+0.5,1) node[right]{$2$};

    \draw (-3+-.5,-1.5) node{$\circ$};
    \draw [->] (-3+-.5,-1.45) -- (-3+-.5,.95);
    \draw [->] (-3+-.55,-1.45) -- (-3+-0.95,-.55);
    \draw (-3+-1,-.5) node{$\circ$};

    \draw (-3+.5,-1.5) node{$\circ$};
    \draw (-3+.5,-1.5) node[right]{$c^0_{l-2}$};
    \draw (-1+.5,-1.5) node[left]{$b^0_{j}$};
    \draw [->] (-3+.5,-1.45) -- (-3+.5,.95);
    \draw [->] (-3+.55,-1.45) -- (-3+0.95,-.55);
    \draw (-3+1,-.5) node{$\circ$};
    \draw (-3+1.5,-1.5) node{$\circ$};
    \draw [->] (-3+1.45,-1.45) -- (-3+1.05,-.55);
  \end{tikzpicture}
  \]
  \caption{For Claim 2(a) in the proof of Lemma~\ref{lem:injective-map}}
  \label{fig:aux-proof-1}
\end{figure}

\begin{figure}[ht]
    \[
    \begin{tikzpicture}
        \draw (-5+-1,0) node[left]{(a)};
      \draw (2.5,0.5) node{$\cdots$};
      \draw (0,0) node{$\circ$};
      \draw (.5,1) node{$\circ$};
      \draw [->] (.05,.05) -- (0.45,.95);
      \draw (1,0) node{$\circ$};
      \draw [->] (.95,.05) -- (0.55,.95);
      \draw (1.5,1) node{$\circ$};
      \draw (1.5,1) node[right]{$6$};
      \draw (2,0) node{$\circ$};
      \draw (2,0) node[right]{$7$};
      \draw [->] (1.95,.05) -- (1.55,.95);
      
      \draw (-.5,1) node{$\circ$};
      \draw [->] (-.05,.05) -- (-0.45,.95);
      \draw (-1,0) node{$\circ$};
      \draw (-1,0) node[right]{$1$};
      \draw [->] (-.95,.05) -- (-0.55,.95);
      \draw [->] (1.05,0.05) -- (1.45,.95);
      \draw (0,0) node[right]{$3$};
      \draw (0,0) node[below]{$a_{j+1}$};
      \draw (-0.5,1) node[right]{$2$};
      \draw (1,0) node[left]{$5$};
      \draw (1,0) node[right]{$a_{j+3}$};
      \draw (0.5,1) node[right]{$4$};
  
      \draw (-.5,-1.5) node{$\circ$};
      \draw (-.5,-1.5) node[right]{$1$};
      \draw [->] (-.5,-1.45) -- (-.5,.95);
      \draw [->] (-.55,-1.45) -- (-0.95,-.55);
      \draw (-1,-.5) node{$\circ$};
      \draw (-1,-.5) node[below]{$0$};
      \draw (-1,-.5) node[right]{$b^1_j$};
  
      \draw (.5,-1.5) node{$\circ$};
      \draw (.5,-1.5) node[right]{$4'$};
      \draw (.5,-1.5) node[below]{$c^0_{j+2}$};
      \draw [->] (.5,-1.45) -- (.5,.95);
      \draw [->] (.55,-1.45) -- (0.95,-.55);
      \draw (1,-.5) node{$\circ$};
      \draw (1,-.5) node[right]{$4''$};
      \draw (1.5,-1.5) node{$\circ$};
      \draw (1.5,-1.5) node[right]{$4'''$};
      \draw [->] (1.45,-1.45) -- (1.05,-.55);

      \draw (-4+-1.5,0.5) node{$\cdots$};
      \draw (-4+0,0) node{$\circ$};
      \draw (-4+.5,1) node{$\circ$};
      \draw [->] (-4+.05,.05) -- (-4+0.45,.95);
      \draw (-4+1,0) node{$\circ$};
      \draw [->] (-4+.95,.05) -- (-4+0.55,.95);
      \draw (-4+1.5,1) node{$\circ$};
      \draw (-4+1.5,1) node[right]{$2$};
      \draw (-3+1.5,1) node{$\circ$};
      \draw (-3+1.5,1) node[right]{$0$};
      \draw (-3+1.5,1) node[above]{$a_l$};
      \draw (-4+2,0) node{$\circ$};
      \draw [->] (-4+1.95,.05) -- (-4+1.55,.95);
      \draw [->] (-3+1.95,.05) -- (-3+1.55,.95);
      
      \draw (-4+-.5,1) node{$\circ$};
      \draw [->] (-4+-.05,.05) -- (-4+-0.45,.95);
      \draw (-4+-1,0) node{$\circ$};
      \draw [->] (-4+-.95,.05) -- (-4+-0.55,.95);
      \draw [->] (-4+1.05,0.05) -- (-4+1.45,.95);
      \draw [->] (-3+1.05,0.05) -- (-3+1.45,.95);
      \draw (-3+1,0) node[right]{$1$};
      \draw (-4+1,0) node[right]{$3$};
      \draw (-4+0.5,1) node[right]{$4$};
      \draw (-4+0.5,1) node[above]{$a_{l-4}$};
  
      \draw (-4+-.5,-1.5) node{$\circ$};
      \draw [->] (-4+-.5,-1.45) -- (-4+-.5,.95);
      \draw [->] (-4+-.55,-1.45) -- (-4+-0.95,-.55);
      \draw (-4+-1,-.5) node{$\circ$};
  
      \draw (-4+.5,-1.5) node{$\circ$};
      \draw (-4+.5,-1.5) node[right]{$c^0_{l-4}$};
      \draw [->] (-4+.5,-1.45) -- (-4+.5,.95);
      \draw [->] (-4+.55,-1.45) -- (-4+0.95,-.55);
      \draw (-4+1,-.5) node{$\circ$};
      \draw (-4+1.5,-1.5) node{$\circ$};
      \draw [->] (-4+1.45,-1.45) -- (-4+1.05,-.55);
    \end{tikzpicture}
    \]
  \[
  \begin{tikzpicture}
    \draw (-5+-1,0) node[left]{(b)};
    \draw (3,0.5) node{$\cdots$};
    \draw (0,0) node{$\circ$};
    \draw (.5,1) node{$\circ$};
    \draw [->] (.05,.05) -- (0.45,.95);
    \draw (1,0) node{$\circ$};
    \draw [->] (.95,.05) -- (0.55,.95);
    \draw (1.5,1) node{$\circ$};
    \draw (1.5,1) node[right]{$6$};
    \draw (2,0) node{$\circ$};
    \draw (2,0) node[right]{$7$};
    \draw [->] (1.95,.05) -- (1.55,.95);
    \draw (2.5,1) node{$\circ$};
    \draw (2.5,1) node[right]{$8$};
    \draw (2.5,1) node[above]{$a_{j+6}$};
    \draw [->] (2.05,.05) -- (2.45,.95);
    
    \draw (-.5,1) node{$\circ$};
    \draw [->] (-.05,.05) -- (-0.45,.95);
    \draw (-1,0) node{$\circ$};
    \draw (-1,0) node[right]{$1$};
    \draw [->] (-.95,.05) -- (-0.55,.95);
    \draw [->] (1.05,0.05) -- (1.45,.95);
    \draw (0,0) node[right]{$3$};
    \draw (-0.5,1) node[right]{$2$};
    \draw (1,0) node[right]{$5$};
    \draw (.5,1) node[above]{$a_{j+2}$};
    \draw (0.5,1) node[right]{$4$};

    \draw (-.5,-1.5) node{$\circ$};
    \draw (-.5,-1.5) node[right]{$1$};
    \draw [->] (-.5,-1.45) -- (-.5,.95);
    \draw [->] (-.55,-1.45) -- (-0.95,-.55);
    \draw (-1,-.5) node{$\circ$};
    \draw (-1,-.5) node[below]{$0$};
    \draw (-1,-.5) node[right]{$b^1_j$};

    \draw (.5,-1.5) node{$\circ$};
    \draw (.5,-1.5) node[right]{$4'$};
    \draw [->] (.5,-1.45) -- (.5,.95);
    \draw [->] (.55,-1.45) -- (0.95,-.55);
    \draw (1,-.5) node{$\circ$};
    \draw (1,-.5) node[right]{$4''$};
    \draw (1.5,-1.5) node{$\circ$};
    \draw (1.5,-1.5) node[right]{$4'''$};
    \draw [->] (1.45,-1.45) -- (1.05,-.55);

    \draw (-4+-1.5,0.5) node{$\cdots$};
    \draw (-4+0,0) node{$\circ$};
    \draw (-4+.5,1) node{$\circ$};
    \draw [->] (-4+.05,.05) -- (-4+0.45,.95);
    \draw (-4+1,0) node{$\circ$};
    \draw [->] (-4+.95,.05) -- (-4+0.55,.95);
    \draw (-4+1.5,1) node{$\circ$};
    \draw (-4+1.5,1) node[right]{$2$};
    \draw (-3+1.5,1) node{$\circ$};
    \draw (-3+1.5,1) node[right]{$0$};
    \draw (-3+1.5,1) node[above]{$a_l$};
    \draw (-4+2,0) node{$\circ$};
    \draw [->] (-4+1.95,.05) -- (-4+1.55,.95);
    \draw [->] (-3+1.95,.05) -- (-3+1.55,.95);
    
    \draw (-4+-.5,1) node{$\circ$};
    \draw [->] (-4+-.05,.05) -- (-4+-0.45,.95);
    \draw (-4+-1,0) node{$\circ$};
    \draw [->] (-4+-.95,.05) -- (-4+-0.55,.95);
    \draw [->] (-4+1.05,0.05) -- (-4+1.45,.95);
    \draw [->] (-3+1.05,0.05) -- (-3+1.45,.95);
    \draw (-4+0,0) node[right]{$5$};
    \draw (-4+-0.5,1) node[right]{$6$};
    \draw (-3+1,0) node[right]{$1$};
    \draw (-4+1,0) node[right]{$3$};
    \draw (-4+0.5,1) node[right]{$4$};

    \draw (-4+-.5,-1.5) node{$\circ$};
    \draw (-4+-1,-.5) node[left]{$b^1_{l-6}$};
    \draw [->] (-4+-.5,-1.45) -- (-4+-.5,.95);
    \draw [->] (-4+-.55,-1.45) -- (-4+-0.95,-.55);
    \draw (-4+-1,-.5) node{$\circ$};

    \draw (-4+.5,-1.5) node{$\circ$};
    \draw (-4+.5,-1.5) node[right]{$4'$};
    \draw [->] (-4+.5,-1.45) -- (-4+.5,.95);
    \draw [->] (-4+.55,-1.45) -- (-4+0.95,-.55);
    \draw (-4+1,-.5) node{$\circ$};
    \draw (-4+1,-.5) node[right]{$4''$};
    \draw (-4+1.5,-1.5) node{$\circ$};
    \draw (-4+1.5,-1.5) node[right]{$4'''$};
    \draw [->] (-4+1.45,-1.45) -- (-4+1.05,-.55);
  \end{tikzpicture}
  \]
  \caption{For Claim 2(b) in the proof of Lemma~\ref{lem:injective-map}}
  \label{fig:aux-proof-2}
\end{figure}

\begin{lemma}\label{lem:partial-iso-Zalpha}
    Let $\G=(Y,S)\in\mathsf{Fr}_r(\mathsf{Log}(\Z_\alpha))$ be an infinite frame. Then for all $k\in\omega$ and $y\in Y$, $\G\rsto S_\sharp^k[y]\rightarrowtail\Z_\alpha$.
\end{lemma}
\begin{proof}
    Since $\G\md\mathsf{Log}(\Z_\alpha)$ and $\G,y\not\md\neg\J^{k+26}(\G,y)$, there exists $x\in Z$ such that $\Z_\alpha,x\not\md\neg\J^{k+26}(\G,y)$. By Proposition~\ref{lem:JankovLemma-k}, there exists a map $f:(\Z_\alpha,x)\to^{k+26}(\F,w)$. By Lemma~\ref{lem:injective-map}, $g=f\rsto \R^k[x]$ is injective. By Fact~\ref{fact:inj-tmorphism}, $g:\R^k[x]\iso\R^k[y]$. 
\end{proof}

\begin{lemma}
    If $\alpha:\mathbb{Z}\to 2$ is finitely perfect, then $\mathsf{Log}(\Z_\alpha)$ is pretabular.
\end{lemma}
\begin{proof}
    Let $L\supseteq\mathsf{Log}(\Z_\alpha)$ be non-tabular. By Theorem~\ref{thm:nontabular-infrootedframe}, $L\sub\mathsf{Log}(\gf)$ for some rooted refined frame $\gf$. Note that $\gf\md\mathbf{alt}^+_3\wedge\mathbf{alt}^-_4$, we see that $\gf$ is image-finite. By Lemma~\ref{lem:pointwise-finite-kappaF}, $\mathsf{Log}(\gf)=\mathsf{Log}(\kappa\gf)$. Let $\kappa\gf=G=(Y,S)$. It suffices to show that $\mathsf{Log}(\Z_\alpha)\supseteq\mathsf{Log}(\G)$. Take any $\phi\not\in\mathsf{Log}(\Z_\alpha)$. Then $\Z_\alpha,z\not\md\phi$ for some $z\in Z$ and there exists a finite subsequence $\beta$ of $\alpha$ such that $\Z_\alpha\rsto\R^{\mathsf{md}(\phi)}[z]\rightarrowtail\Z_\beta$. Recall that $\alpha$ is finitely perfect, there exists $n\in\omega$ such that $\beta\triangleleft\gamma$ for all $\gamma\triangleleft\alpha$ such that $|\gamma|>n$. Let $m=8(n+3)$. Take any $x\in Y$. By Lemma~\ref{lem:partial-iso-Zalpha}, $\G\rsto S_\sharp^m[y]\cong\Z$ for some $\Z\rightarrowtail\Z_\alpha$. Then $\mathsf{zdg}(Z)\geq m$. By the construction of $\Z_\alpha$, $\Z_\gamma\rightarrowtail\Z$ for some $\gamma\triangleleft\alpha$ with $|\gamma|\geq n$. Thus $\Z_\alpha\rsto\R^{\mathsf{md}(\phi)}[z]\rightarrowtail\Z_\beta\rightarrowtail\Z_\gamma\rightarrowtail\Z\iso\G\rsto S_\sharp^m[y]$, which implies $\phi\not\in\mathsf{Log}(\G)$. Hence $\mathsf{Log}(\Z_\alpha)\supseteq\mathsf{Log}(\G)\supseteq L\supseteq\mathsf{Log}(\Z_\alpha)$ and so $\mathsf{Log}(\Z_\alpha)$ is pretabular.
\end{proof}

\begin{corollary}\label{coro:LogZchif-is-pretabular}
    For all $f:\omega\to 2$, the logic $\mathsf{Log}(\Z_{\chi^f})$ is pretabular.
\end{corollary}

\begin{lemma}\label{lem:distinctMap-distinctLogic}
    For all sequences $\alpha,\beta:\mathbb{Z}\to 2$, $\beta\npreceq\alpha$ implies $\mathsf{Log}(\Z_{\alpha})\nsubseteq\mathsf{Log}(\Z_{\beta})$.
\end{lemma}
\begin{proof}
    Let $\beta\npreceq\alpha$. Then $\gamma\ntriangleleft\alpha$ for some $\gamma\triangleleft\beta$. Suppose $\mathsf{Log}(\Z_{\alpha})\subseteq\mathsf{Log}(\Z_{\beta})$. Then $\Z_{\beta}\in\mathsf{Fr}_r(\mathsf{Log}(\Z_{\alpha}))$. By Lemma~\ref{lem:partial-iso-Zalpha}, $\Z_\gamma\rightarrowtail\Z_{\alpha}$. By the construction of $\Z_\gamma$ and $\Z_{\alpha}$, it is clear that $\gamma\triangleleft\alpha$, which is impossible. Thus $\mathsf{Log}(\Z_{\alpha})\nsubseteq\mathsf{Log}(\Z_{\beta})$.
\end{proof}

As consequences, the following theorems hold: 

\begin{theorem}\label{thm:pretab-BS223}
    $|\mathsf{PTAB}(\mathsf{S4}_t)|\geq|\mathsf{PTAB}(\mathsf{S4BP}^{2,\omega}_{2,3})|=2^{\aleph_0}$.
\end{theorem}
\begin{proof}
    By Lemmas~\ref{lem:TM-sequence-emb} and \ref{lem:distinctMap-distinctLogic}, $|\cset{\mathsf{Log}(\Z_{\chi^f}):f\in 2^\omega}|=2^{\aleph_0}$. By Corollary~\ref{coro:LogZchif-is-pretabular}, we see that $\cset{\mathsf{Log}(\Z_{\chi^f}):f\in 2^\omega}\sub\mathsf{PTAB}(\mathsf{S4BP}^{2,\omega}_{2,3})\sub\mathsf{PTAB}(\mathsf{S4}_t)$.
\end{proof}

\begin{theorem}\label{thm:pretab-S4t}
    For all $\kappa\leq{\aleph_0}$ or $\kappa=2^{\aleph_0}$, $|\mathsf{PTAB}(L)|=\kappa$ for some $L\in\NExt(\mathsf{S4}_t)$.
\end{theorem}

\begin{remark}\label{rem:fmp-pretab}
    As it is shown in this section, the logics $\mathsf{Log}(\Z_{\chi^f})$ are pretabular, Kripke complete and of finite depth. There exist modal logics in $\NExt(\mathsf{S4})$, say $\mathsf{S5}$, which also satisfy these properties. 
    
    However, there are substantial differences between the lattices $\NExt(\mathsf{S4})$ and $\NExt(\mathsf{S4}_t)$. For example, by \cite[Theorems 12.7 and 12.11]{Chagrov.Zakharyaschev1997}, the following claims are true in $\NExt(\mathsf{S4})$, even in $\NExt(\mathsf{K4})$:
    \begin{enumerate}[(i)]
        \item every tabular logic has finitely many immediate predecessors;
        \item every pretabular logic enjoys the finite model property.
    \end{enumerate}
    On the other hand, we have the following conjecture: $\mathsf{Fin}_r(\mathsf{Log}(\Z_{\chi^f}))=\mathsf{IM}(\cset{\C_1,\C_2})$ for any $f\in 2^\omega$. This will give that $\mathsf{Log}(\C_2)$ has a continuum of immediate predecessors and $\mathsf{Log}(\Z_{\chi^f})$ is pretabular but lacks the FMP. Thus neither (i) nor (ii) may hold in $\NExt(\mathsf{S4}_t)$. In order to prove this conjecture, we would need to show that the critical exponent of $\chi^f$ is always finite. We leave this to future research.
\end{remark}

\section{Conclusions}

The present work studies pretabularity in tense logics above $\mathsf{S4}_t$. We started with tense logics $\mathsf{S4BP}_{n,m}^{k,l}$ with bounded parameters. We gave a full characterization of $\mathsf{PTAB}(\mathsf{S4BP}_{n,m}^{k,l})$ for cases where all parameters $k,l,n,m$ are finite. We also investigated some concrete tense logics where some of the parameters are infinite. Full characterizations for pretabular logics extending $\mathsf{S4.3}_t$ and $\mathsf{S4BP}_{2,2}^{2,\omega}$ were provided, where $\mathsf{S4.3}_t=\mathsf{S4BP}_{1,1}^{\omega,\omega}$ and $\mathsf{S4BP}_{2,2}^{2,\omega}$ is closely related to the logic $\mathsf{Ga}$ studied in \cite{Kracht1992}. We showed that $|\mathsf{PTAB}(\mathsf{S4.3}_t)|=5$ and $|\mathsf{PTAB}(\mathsf{S4BP}_{2,2}^{2,\omega})|=\aleph_0$. Finally, we studied pretabular tense logics in $\NExt(\mathsf{S4BP}^{2,\omega}_{2,3})$. By Theorem~\ref{thm:pretab-BS223}, the cardinality of $\mathsf{PTAB}(\mathsf{S4BP}^{2,\omega}_{2,3})$ is $2^{\aleph_0}$. Theorem~\ref{thm:pretab-S4t} answers the open problem raised in \cite{Rautenberg1979}.

In fact, Theorem~\ref{thm:pretab-BS223} suggests that a full characterization of $\mathsf{PTAB}(\mathsf{S4BP}^{2,\omega}_{2,3})$ or $\mathsf{PTAB}(\mathsf{S4}_t)$ is unattainable. Likewise, the decidability of tabularity in $\NExt(\mathsf{S4}_t)$ cannot be obtained via pretabular logics. However, this does not mean that research on pretabular tense logics in $\NExt(\mathsf{S4}_t)$ is exhausted; much remains to be explored. Some open problems have already been mentioned in the remarks (see Remark~\ref{rem:K4BP} and Remark~\ref{rem:Kracht-Ga}) and we outline a few further topics below.

One possible direction for future work is to investigate pretabular logics in other sublattices of $\NExt(\mathsf{S4}_t)$. For example, consider the tense logic $\mathsf{S4BP}_{2,2}^{3,\omega}$, which has the forth-width and back-width 2 and the depth 3. The cardinality of $\mathsf{PTAB}(\mathsf{S4BP}_{2,2}^{3,\omega})$ remains unknown. Similar to pretabular pre-transitive modal logics, the pretabular logics of finite z-degrees are not well-understood.

Another direction for future work is to investigate pretabular logics in $\NExt(\mathsf{S4}_t)$ with the FMP. Pretabular logics can be viewed as boundaries of tabular logics. It is natural to consider that pretabular logics with the FMP act as the limit of certain set of tabular logics. As it is shown in \cite[Theorem 12.11]{Chagrov.Zakharyaschev1997}, every pretabular modal logic in $\NExt(\mathsf{K4})$ has the FMP. By Theorem~\ref{thm:BS222-pretab-fmp}, every pretabular tense logic in $\NExt(\mathsf{S4BP}_{2,2}^{2,\omega})$ has the FMP. However, if our conjecture in Remark~\ref{rem:fmp-pretab} is proved to be correct, then there exists a continuum-sized family of Kripke complete pretabular tense logics lacking the FMP in $\NExt(\mathsf{S4BP}_{2,3}^{2,\omega})$. This raises at least two natural questions: (i) When does $\NExt(L)$ contain pretabular logics lacking the FMP? (ii) How many pretabular logics with the FMP exist in $\NExt(\mathsf{S4}_t)$? Exploring these questions will deepen our understanding of the lattices of tense logics.

\vspace{2em}

\noindent\textbf{Acknowledgement.} The author is grateful to Nick Bezhanishvili for his very helpful and insightful comments, which significantly improved the manuscript. The author would like to thank the referee for their insightful and constructive comments, which have helped improve the clarity and quality of this paper. The author also thanks Tenyo Takahashi for his inspiring suggestions, contributing to the main proof in Section~\ref{sec:BS223}. Finally, the author is indebted to Minghui Ma and Rodrigo Nicolau Almeida for the discussions that helped shape the ideas presented here. 
The author is supported by Chinese Scholarship Council.

\bibliographystyle{siam}

\bibliography{../../0_Bibliography/References}

\end{document}